\documentclass[3p]{elsarticle}

\usepackage{amssymb}
\usepackage{amsmath}
\usepackage{multirow}
\usepackage{subcaption}
\usepackage{booktabs}
\usepackage{rotating}
\usepackage{float}
\usepackage{siunitx}
\usepackage{ulem}
\newcommand{\p}[2]{\partial_{x_#1}#2}
\newcommand{\px}[1]{\partial_{x}#1}
\DeclareMathOperator{\diag}{diag}

\DeclareMathOperator{\Ha}{Ha}
\DeclareMathOperator{\Ma}{Ma}
\DeclareMathOperator{\Mm}{M_m}
\DeclareMathOperator{\Rh}{Re}
\DeclareMathOperator{\Rm}{R_m}

\numberwithin{equation}{section}
\usepackage{comment}
\usepackage{bm}
\usepackage{xcolor}
\usepackage{hyperref}
\hypersetup{colorlinks=true,
}
\newcommand{\review}[1]{{#1}}

\journal{Journal of Computational Physics}

\begin{document}

\begin{frontmatter}



\title{Entropy conservative and entropy stable solid wall boundary conditions
for the resistive magnetohydrodynamic equations}


\author[label1]{Vladimir Pimanov}
\ead{vladimir.pimanov@kaust.edu.sa}
\cortext[cor1]{Corresponding author.}
\author[label1]{Lisandro Dalcin\corref{cor1}}
\ead{dalcinl@gmail.com}
\author[label1,label2]{Matteo Parsani}
\ead{matteo.parsani@kaust.edu.sa}

\affiliation[label1]{organization={Computer Electrical and Mathematical Science and Engineering Division (CEMSE), 
Extreme Computing Research Center (ECRC), King Abdullah University of Science
    and Technology (KAUST)},%Department and Organization
            %addressline={}, 
            city={Thuwal},
            postcode={}23955-6900, 
            %state={},
            country={Saudi Arabia}}

\affiliation[label2]{organization={Physical Science and Engineering Division (PSE), 
King Abdullah University of Science and Technology (KAUST)},%Department and Organization
            %addressline={}, 
            city={Thuwal},
            postcode={}23955-6900, 
            %state={},
            country={Saudi Arabia}}
\begin{abstract}
We present a novel technique for imposing non-linear entropy conservative and 
entropy stable wall boundary conditions for the resistive magnetohydrodynamic 
equations in the presence of an adiabatic wall or a wall with a prescribed
 heat entropy flow, addressing three scenarios: electrically 
insulating walls, thin walls with finite conductivity, 
and perfectly conducting walls. The procedure relies on the formalism and mimetic 
properties of diagonal-norm, summation-by-parts, and simultaneous-approximation-term 
operators. Using the method of lines, a semi-discrete entropy estimate for the entire 
domain is obtained when the proposed numerical imposition of boundary conditions 
is coupled with an entropy-conservative or entropy-stable discrete interior operator. 
The resulting estimate mimics the global entropy estimate obtained at the continuous level. 
The boundary data at the wall are weakly imposed using a penalty flux approach and a 
simultaneous-approximation-term technique for both the conservative variables and the
gradient of the entropy variables. Discontinuous spectral collocation operators
(mass lumped nodal discontinuous Galerkin operators) on high-order unstructured 
grids are used to demonstrate the new procedure's accuracy, robustness, and efficacy
for weakly enforcing boundary conditions. Numerical simulations confirm the non-linear 
stability of the proposed technique, with applications to three-dimensional flows. 
The procedure described is compatible with any diagonal-norm summation-by-parts spatial
operator, including finite element, finite difference, finite volume, nodal and modal 
discontinuous Galerkin, and flux reconstruction schemes.
\end{abstract}

\begin{graphicalabstract}
\end{graphicalabstract}


\begin{keyword}
Magnetohydrodynamic equations \sep Electrically conducting walls \sep 
Electrically insulating walls \sep
Entropy conservation \sep
Entropy stability \sep
Summation-by-parts operators \sep 
Simultaneous-approximation-terms



\end{keyword}

\end{frontmatter}


%\documentclass{article}
%\usepackage{graphicx} % Required for inserting images
%\usepackage{biblatex}
%\addbibresource{mhd.bib}
%\usepackage{amsmath}
%\usepackage{amsfonts}
%\usepackage{amsthm}
%\usepackage{comment}
\newtheorem{thm}{Theorem}
\newproof{pf}{Proof}
\newtheorem{remark}{Remark}
\newcommand{\Mnode}{(-)}
\newcommand{\Pnode}{(+)}
\newcommand{\nodeA}{\mathcal{A}}
\newcommand{\nodeB}{\mathcal{B}}
\newcommand{\Bextvec}{\vec B^0}
\newcommand{\Bext}{B^0}
\newcommand{\jump}[1]{[\![#1]\!]}
\newcommand{\jumpDB}[1]{[\![#1]\!]_{(D,B)}}
\newcommand{\meanDB}[1]{\{\!\{#1\}\!\}_{(D,B)}}
\newcommand{\logavr}[1]{{#1}^\text{ln}}
\newcommand{\avr}[1]{\{\!\!\{#1\}\!\!\}}
%\let\oldref\ref
%\renewcommand{\ref}[1]{(\oldref{#1})}
%\newtheorem{corollary}{Corollary}[theorem]
%\newtheorem{lemma}[theorem]{Lemma}
%\usepackage{graphicx}
%\graphicspath{ {./images/} }

%\title{Entropy stable and entropy conservative wall boundary conditions for resistive MHD equations}
%\author{Vladimir Pimanov, ...}
%\date{September 2023}

%\begin{document}

%\maketitle

%\section{}
%\label{}





\section{Introduction}
\label{intro_section_label}

The resistive magnetohydrodynamics (MHD) equations
\cite{godunov1972symmetricMHD,davidson2002resistbook,Davidson_2001,hosking2016fundamental} govern the motion of ionized
plasma and electrically conducting fluids and are of great interest in space
physics, astrophysics, and many other technical applications. One of these critical
applications is nuclear fusion energy generation \cite{harms2000principles,marcus2022systems}. 

The resistive MHD equations
for compressible flows explicitly establish the conservation of mass, momentum,
total energy, and magnetic fields. These equations also encapsulate two
implicit principles. The first principle demands a divergence-free constraint
on the magnetic field. At the continuum level, a divergence-free initial state
guarantees the fulfillment of this principle at any subsequent time. However, numerical 
schemes
often struggle to 
maintain this property due to accumulating numerical errors. Thus, special attention 
is required to fulfill this principle at the semi-discrete and fully-discrete levels. 
The second fundamental property of the resistive MHD equations is they satisfy
the second law of thermodynamics. The total 
thermodynamic entropy can not decrease in time if the physical system is closed. 
Its variation is associated with 
two mechanisms: the production of entropy on discontinuities of the solution and
diffusion mechanisms. The MHD equations meet this second principle because they 
are equipped with an additional non-conservative term proportional to the 
magnetic field divergence \cite{godunov1972symmetricMHD}. This non-conservative 
term is known as the Powell term \cite{PowellTerm1999}. Also, in this case, numerical 
schemes often struggle to 
maintain this property due to accumulating numerical errors.
%Very often, this
%thermodynamic principle is at the core of the construction of semi-discrete and
%fully-discrete algorithms that are non-linearly stable (depending on the field, 
%which is also called entropy stable or Lyapunov stable algorithms).

%The compatibility of this principle with resistive MHD equations is proven in \cite{DGSEM_MHD_Part1}. 

%An accurate representation of the conservation properties of continuous 
%equations in numerical schemes is crucial for the reliability of solutions. 
Having the
second law of thermodynamics built-in plays a unique role.
Very often, this
thermodynamic principle is at the core of the construction of semi-discrete and
fully-discrete algorithms that are non-linearly stable (depending on the field, 
this property is also called entropy stability or Lyapunov stability) --
provided that the positivity of the solution is preserved.
%it allows to proof 
%nonlinear stability of the continuous or discrete system in some cases. 
In particular, in closed systems, under the physically informed assumptions of 
density and pressure positivity, the fulfillment of the second law of
thermodynamic establishes an upper bound on the
$L_2$ norm of the conservative variables (see, for instance,
\cite{Dafermos2005book})
which can be translated into an upper bound of the primitive variables. 
Mimicking the continuous analyses at the discrete level is a rigorous way to 
develop robust, non-linearly stable numerical discretizations.

At the semi-discrete level, this work focuses on collocated discontinuous Galerkin (DG)
discretizations for
unstructured grids. For a high-order accurate representation of the solution, these methods deliver a high
computation-to-communication ratio, which makes them a natural choice for the 
current and future generation of heterogeneous data-centric computing hardware 
\cite{fischer2020scalability, isaac2015support, rojas2021robustness}. In particular, they deliver
better accuracy per degree of freedom than low-order counterparts, which leads
to much smaller numerical errors in terms of dispersion and dissipation 
\cite{hesthaven2017numerical, wang2013high, reynanolasco2022eigen}. However, high-order accurate schemes
generally suffer from instabilities more often than low-order schemes due to discontinuities
in the solution or under-resolved physical features.  

In the last two decades, entropy stable discontinuous Galerkin methods have been developed
to solve several systems of hyperbolic equations, including the compressible Euler equations and 
the compressible Navier--Stokes equations; see, for instance, \cite{fisher2013high,carpenter2014entropy}. An entropy stable discontinuous
Galerkin spectral element method (DGSEM) for the resistive MHD equations is presented for spatially periodic 
solutions in domains without boundaries in \cite{DGSEM_MHD_Part1}. 
In \cite{DGSEM_MHD_Part2}, the solver above is coupled with a classical finite volume 
discretization for shockwave capturing. To deal with divergence-free 
constrain on the magnetic field, the numerical schemes in \cite{DGSEM_MHD_Part1, DGSEM_MHD_Part2}
are designed using a generalized Lagrange
multiplier, providing divergence-cleaning properties to the spatial
discretizations \cite{IdealGLMMHD}.  

Issues remain on the path toward complete entropy stability for the MHD 
equations, \emph{e.g.}, bound-preserving limiter for high-order accurate discretizations.
In particular, one major obstacle is the need for boundary conditions that 
preserve the entropy
conservation or stability property of the interior operators. Practical 
experience indicates that
numerical instabilities frequently originate at domain boundaries, where
``the action takes place". For instance, the interaction of shocks with these physical 
boundaries is particularly challenging for high-order formulations. In this 
article, we expand upon the method presented in \cite{DGSEM_MHD_Part1} to
introduce entropy conservative and entropy stable solid wall boundary conditions. 
The starting points are the entropy stable solid wall boundary conditions developed 
for the compressible Navier--Stokes 
equations in \cite{dalcin2019conservative,parsani2015entropy} for collocated 
discontinuous Galerkin methods with the
summation-by-parts property (DG-SBP) and simultaneous-approximation-terms (SATs) 
operators
\cite{carpenter2014entropy,carpenter2016towards,fernandez2020entropy}. To develop the new
boundary conditions, we
consider a general conducting wall with conductance parameter $c$ and address three scenarios: electrically 
insulating walls ($c=0$), thin walls with finite conductivity ($0<c<\infty$), 
and perfectly conducting walls ($c=\infty$). 
%These boundary conditions, 
%particularly in the incompressible limit, enable the derivation of %analytical 
%solutions for flows in pipes and channels, including plane channel flow 
%\%cite{hartmann1937hg}, flow in rectangular ducts
%\cite{davidson2002resistbook}, and circular pipe flow 
%\cite{ihara1967flow}. 

We demonstrate the accuracy of the proposed boundary 
conditions through a convergence study towards an analytical solution in a circular pipe under a 
varying wall conductance parameter. The choice of pipe 
flow as a case study is strategic, showcasing the method's efficacy in curved 
geometries, where meshes with curved elements are often used in conjunction with
high-order accurate schemes. To verify the entropy conservative properties of the spatial 
discretization algorithm and our proposed hydrodynamic, thermal, and magnetic boundary conditions, 
we study the time rate of change of the entropy function of a flow \review{around} a stationary rotating spheroid. Finally, we test
the new boundary conditions in two applications: fluid flow in a microchannel and a self-field magnetoplasmadynamic thruster.

The manuscript is organized as follows. 
Section \ref{equations_section_label} introduces the continuous formulation of the problem, 
including the notation and the system of governing equations, boundary conditions, 
and a thermodynamic analysis of the system. 
Section \ref{discret_section_label} details the one-dimensional DG-SBP-SAT discretization 
of the system and dives into the development of the new entropy-stable solid-wall boundary conditions. 
Section \ref{results_section_label} discusses the numerical results that support
the primary goal of this study.
%, including a convergence study 
%for a circular pipe flow with different wall conductance parameter values. 
The concluding remarks are presented in Section \ref{conclusions_section_label}.
%\ref{iflux_appendix_label} defines the entropy
%conservative
%two-point flux function used in the discretization.
%\ref{entropy_proff_appendix_label} presents an analysis of the numerical entropy production 
%within a single mesh element, as utilized in Section \ref{discret_section_label}.

\section{Resistive compressible MHD equations}
\label{equations_section_label}
The compressible resistive magnetohydrodynamics (MHD) equations \cite{godunov1972symmetricMHD} and
the formulation of an entropy stable discontinuous Galerkin method with the
summation-by-part property (DG-SBP) are described in detail in
\cite{DGSEM_MHD_Part1}. Here, we summarize the continuous equations
and their key properties, keeping the notation close to that of
\cite{DGSEM_MHD_Part1}. We also extend the continuous analysis to
include entropy conservative and entropy stable solid wall boundary conditions.

\subsection{Notation}

We assume to work in a three-dimensional space, \emph{i.e.}, in $\mathbb{R}^3$.
Spatial vectors are denoted with an arrow (\emph{e.g.}, $\vec v = (v_1, v_2, v_3) \in \mathbb{R}^{1 \times 3}$).
State vectors are denoted in bold
(\emph{e.g.}, $\mathbf u = (\rho, \rho \vec v, E, \vec B, \psi)^T \in \mathbb{R}^{9 \times 1}$).
Block vectors are vectors that contain state vectors in each spatial direction, they are denoted in bold with double arrow
(\emph{e.g.}, $\overset{\leftrightarrow}{\mathbf f} = (\mathbf f_1, \mathbf f_2, \mathbf f_3) \in \mathbb R^{9 \times 3}$).

In Cartesian coordinates $(x_1,x_2,x_3)$, the gradient of a state vector is a block vector,
\begin{equation}
    \vec \nabla \mathbf u = (\p{1}{\mathbf u}, \p{2}{\mathbf u}, \p{3}{\mathbf u}).
\end{equation}
The gradient of a spatial vector is a second order tensor and in matrix form is given by
\begin{equation}
    \vec \nabla \vec v =
    \begin{pmatrix}
        \p{1}{v_1} & \p{2}{v_1} & \p{3}{v_1} \\
        \p{1}{v_2} & \p{2}{v_2} & \p{3}{v_2} \\
        \p{1}{v_3} & \p{2}{v_3} & \p{3}{v_3}
    \end{pmatrix}.
\end{equation}
The divergence of a block vector is defined as
\begin{equation}
    \vec \nabla \cdot \overset{\leftrightarrow}{\mathbf f} = \p{1}{\mathbf f_1} +
    \p{2}{\mathbf f_2} + \p{3}{\mathbf f_3}.
\end{equation}
We additionally define the jump operator, arithmetic, and logarithmic means between two values $a_L$ and $a_R$ as
\begin{equation}
    [\![a]\!]_{(L,R)} := a_R - a_L, \quad
    \{\!\{a\}\!\}_{(L,R)} := \frac{1}{2}(a_R + a_L), \quad
    a^\text{ln}_{(L,R)} := [\![a]\!]_{(L,R)} / [\![\text{ln}(a)]\!]_{(L,R)}.
\end{equation}
A numerically stable procedure to evaluate the logarithmic mean is given in \cite{IsmailRoe2009}.

\subsection{The system of PDEs}
We consider the compressible visco-resistive MHD equations compatible with the 
continuous entropy analysis presented by Godunov in 1972 \cite{godunov1972symmetricMHD}. 
This system of PDEs reads 
\begin{equation}
    \label{conservation_law_eq}
\begin{split}
    & \partial_t \mathbf{u} + 
    \vec \nabla \cdot \overset{\leftrightarrow}{\mathbf{f}}\!^a(\mathbf{u}) +
    \mathbf{\Upsilon}(\mathbf{u}, \vec \nabla \mathbf{u}) = 
    \vec \nabla \cdot \overset{\leftrightarrow}{\mathbf{f}}\!^\nu(\mathbf{u}, \vec \nabla \mathbf{u}) + \mathbf r(\mathbf u), 
    \quad \forall \vec x \in \Omega, \quad t > 0, \\
    &\Pi(\mathbf u, \vec \nabla \mathbf u) = 0, \quad \forall \vec x \in \Gamma, \quad t > 0, \\
    %& c_1 \mathbf u(\vec x, t) + c_2\ \partial_n \mathbf u(\vec x, t) = \mathbf u^\text{wall}(\vec x, t), \quad \forall \vec x \in \Gamma, \quad t \ge 0, \\
    & \mathbf u = \mathbf u^\text{init}, \quad \forall \vec x \in \Omega, \quad t=0,
\end{split}
\end{equation}
where $\mathbf{u}$ is the state vector, 
$\overset{\leftrightarrow}{\mathbf{f}^a}$ and $\overset{\leftrightarrow}{\mathbf{f}^\nu}$ 
are advective and diffusive fluxes, and $\mathbf{\Upsilon}$ is the
non-conservative term which
enhances the numerical properties of system \eqref{conservation_law_eq} and ensures it is symmetrizable,
extending the original ideas from \cite{godunov1972symmetricMHD,PowellTerm1999} to the 
ideal GLM-MHD system reported in \cite{IdealGLMMHD}.
%that appears by rewriting the ideal MHD equations in conservative form \cite{IdealGLMMHD}.
%\textcolor{red}{Read part I and II of Gassner}
%\textcolor{red}{Say when fv = 0}
The vector $\mathbf r$ is a purely algebraic source term introduced to dampen the 
magnetic field's divergence-free numerical error. The Greek
letter $\Omega$ indicates the spatial domain with boundary $\Gamma$. 
Boundary conditions are expressed in a general form through an implicit function $\Pi$.
In the generalized Lagrangian multiplier (GLM) formulation of the MHD equations
(GLM-MHD), the state vector is given by $\mathbf u = (\rho, \rho \vec v, E, \vec B, \psi)^T$, 
where $\rho$ is the density, $\vec v = (v_1, v_2, v_3)$ is the velocity, 
$E$ is the total energy, $\vec B = (B_1, B_2, B_3)$ is the magnetic field, 
$\psi$ is the generalized Lagrange multiplier introduced in \cite{IdealGLMMHD} for the purpose of
magnetic field divergence cleaning. The introduction of the generalized Lagrange
multiplier builds on the original ideas from \cite{MUNZ2000484, DEDNER2002645}.
Note that if the viscous fluxes are set to zero, the ideal MHD for compressible
flows are recovered. 
In this work, we consider the perfect gas model.

The advective flux, split into the Euler (\emph{i.e.}, the inviscid part of the compressible
Navier--Stokes equations), ideal MHD, and GLM parts, is
\begin{equation}
    \overset{\leftrightarrow}{\mathbf{f}}\!^a(\mathbf u) = 
    \begin{pmatrix}
        \rho \vec v \\
        \rho \vec v^T \vec v + p \underline{I} \\
        \left( \frac{1}{2} \rho  |\!|\vec v|\!|^2 + \frac{\gamma p}{\gamma-1}\right) \vec v \\
        \underline{0} \\
        \vec 0
    \end{pmatrix} + 
    \begin{pmatrix}
        \vec 0 \\
        \frac{1}{2\mu_0}|\!|\vec B|\!|^2 \underline{I} - \frac{1}{\mu_0} \vec B^T \vec B \\
        \frac{1}{\mu_0} \left( |\!|\vec B|\!|^2 \vec v  - \left( \vec v \cdot \vec B \right) \vec B \right) \\
        \vec B^T \vec v - \vec v^T \vec B \\
        \vec 0
    \end{pmatrix} + 
    \begin{pmatrix}
        \vec 0 \\
        \underline{0} \\
        \frac{c_h}{\mu_0} \psi \vec B \\
        c_h \psi \underline{I} \\
        c_h \vec B
    \end{pmatrix},
\end{equation}
where $p$ is the gas pressure, $\underline I$ is the $3\times3$ identity matrix, $\mu_0$ is the permeability of the medium, 
$\gamma$ is the heat capacity ratio, and $c_h$ is the hyperbolic divergence cleaning speed.

The non-conservative term, $\mathbf{\Upsilon}$, consists of two components, namely
$\mathbf{\Upsilon} = \mathbf{\Upsilon}^\text{MHD} + \mathbf{\Upsilon}^\text{GLM}$.
The first non-conservative term, $\mathbf{\Upsilon}^\text{MHD}$, also known as
the Powell term \cite{PowellTerm1999}, reads
\begin{equation}
    \mathbf{\Upsilon}^\text{MHD} = \left(\vec \nabla \cdot \vec B \right) \pmb{\phi}^\text{MHD} = 
    \left(\vec \nabla \cdot \vec B \right) \left(0, \frac{1}{\mu_0} \vec B,
    \frac{1}{\mu_0} \vec v \cdot \vec B, \vec v, 0\right)^T.
\end{equation}
The second non-conservative term, $\mathbf{\Upsilon}^\text{GLM}$, guarantees 
Galilean invariance of the GLM-MHD equations and is given by
\begin{equation}
    \mathbf{\Upsilon}^\text{GLM} = \overset{\leftrightarrow}{\pmb{\phi}}\!^\text{GLM} \cdot \vec \nabla \psi = 
    \pmb{\phi}_1\!^\text{GLM} \p{1}{\psi} + \pmb{\phi}_2\!^\text{GLM} \p{2}{\psi} + \pmb{\phi}_3\!^\text{GLM} \p{3}{\psi},
\end{equation}
where $\pmb{\phi}_l\!^\text{GLM} = \left(0,\vec 0,\mu_0^{-1} v_l \psi,\vec 0,v_l\right)^T$.

%\textcolor{red}{MP: To discusse here}
At the continuous level, $\vec \nabla \cdot \vec B = 0$. Consequently, the
non-conservative term,
$\mathbf{\Upsilon}^\text{MHD}$, is also zero. 
%However, at the discrete level, a nonzero divergence of the magnetic field is allowed.
%The addition of the non-conservative term, $\mathbf{\Upsilon}$, is know to have some favorable advantages in numerical simulations respect to the ideal MHD equations.
%\textcolor{red}{MP: To discusse here}
%Note that this extra term accommodates for the conservation of entropy in the advective part of the system. 
For a more detailed discussion of the model and its derivation, we refer the reader to~\cite{IdealGLMMHD}.

The diffusive flux includes electroresistive, viscous, and thermoconductive effects and reads
\begin{equation}\label{eq:diffusive-flux-mhd}
    \overset{\leftrightarrow}{\mathbf{f}}\!^\nu(\mathbf u, \vec \nabla \mathbf u) = 
    \begin{pmatrix}
        \vec 0 \\
        \underline{\tau} \\
        \vec v \underline{\tau} - \vec q - \frac{\mu_R}{\mu_0^2} (\vec \nabla \times \vec B) \times \vec B  \\
        \frac{\mu_R}{\mu_0} \left( \vec \nabla \vec B - (\vec \nabla \vec B)^T \right) \\
        \vec 0
    \end{pmatrix},
\end{equation}
where $\underline \tau$ is the stress tensor, %{\color{red} Sign was changed to meet 1D formulation and entropy flux condition}
\begin{equation}
    \underline{\tau} = \mu_{NS} \left( \vec \nabla \vec v + 
    ( \vec \nabla \vec v)^T  \right) - \frac{2}{3} \mu_{NS} \left(\vec \nabla \cdot \vec v \right) \underline{I},
\end{equation}
and $\vec q$ is the heat flux, linearly proportional to the gradient of temperature, $T$, via the Fourier's law, \emph{i.e.},
\begin{equation}
    \vec q = -\kappa \vec \nabla T. %\left( \frac{p}{R \rho} \right),
\end{equation}
The physical parameters $\mu_{NS}$, $\mu_R$, and $\kappa$ are the fluid dynamic viscosity, 
the resistivity of the plasma, and the thermal conductivity, respectively. The symbol $R$ indicates the specific gas constant
and is given by $R = R_u/M_w$, where $R_u$ is the universal gas constant, and $M_w$ is the molecular weight of the gas.
Here, because we consider the perfect gas model, the equation for the temperature  is given by $T = p/\rho R$ whereas, 
the pressure, extended to account for the GLM term, reads 
\begin{equation}
    p = (\gamma - 1) \left( 
        E - \frac{1}{2}\rho |\!|\vec v|\!|^2 - \frac{1}{2\mu_0}|\!|\vec B|\!|^2 - \frac{1}{2 \mu_0} \psi^2
    \right).
\end{equation}

If otherwise stated, the algebraic source term $\mathbf r$ is given by
\begin{equation}
    \label{r_assumption}
    \mathbf r = \left(0,\vec 0,0,\vec 0,-\alpha\psi\right)^T,
\end{equation}
with $\alpha \ge 0$. We reiterate that introducing this algebraic source term 
enhances the cleaning of the $\vec \nabla \cdot \vec B$ at the semi-discrete 
level \cite{DEDNER2002645,TRICCO2016326} (\emph{i.e.}, of the spatial discretization). Specifically, an additional damping on 
the numerical error of the divergence-free magnetic field is attained if $\alpha > 0$.
%According to~\cite{IdealGLMMHD}, $\alpha = c_h / 0.18$ provides optimal divergence cleaning properties of the scheme. 
% The rational choose of $\alpha$ is discussed in \cite{IdealGLMMHD}.
%\textcolor{red}{Cite the original paper and not Gassner}

\subsection{Thermodynamic properties of the system}\label{sec:eq-thermo}
As shown in \cite{DGSEM_MHD_Part1,IdealGLMMHD}, 
the GLM-MHD system \eqref{conservation_law_eq} adheres to the second law of thermodynamics.
%The GLM-MHD system \eqref{conservation_law_eq} 
In fact, it possesses a strictly convex entropy function, $S(\mathbf u)$,
\begin{equation}
\label{math_entropy_eq}
    S(\mathbf u) = -\rho s,
\end{equation}
which is proportional (with the negative sign) to the specific thermodynamic entropy 
\begin{equation}
  \label{specific_entropy}
    s = \frac{R}{\gamma-1} \ln\left(\frac{T}{T^*}\right) - R \ln\left(\frac{\rho}{\rho^*}\right),
\end{equation}
where $T^*$ and $\rho^*$ are the reference values for temperature and pressure, respectively.
The entropy function $S(\mathbf u)$ is strictly convex on its arguments in a bounded domain $\Omega$ 
under the physically informed assumption that 
\begin{equation}\label{eq:positivity_assumption}
    T > 0, \quad \rho > 0, \quad \forall \vec x \in \Omega, \quad t \ge 0.
\end{equation}
In this work, we assume \eqref{eq:positivity_assumption} to be always true. 
The convexity of the entropy function is a valuable tool for providing stability for the PDE system in the $L^2$ norm \cite{Dafermos2005book, svard2015weak}.

Differentiating the entropy function with respect to the conservative variables leads to the definition of the entropy variables $\mathbf w$,
\begin{equation}
\label{entropy_var_eq}
    \mathbf w^T = \frac{\partial S}{\partial \mathbf u} =
    \left( \frac{\gamma R}{\gamma-1} - s - \frac{1}{2}\frac{|\!|\vec v|\!|^2}{T}, \ \frac{\vec v}{T},
     \ -\frac{1}{T}, \ \frac{1}{\mu_0} \frac{\vec B}{T}, \ \frac{1}{\mu_0} \frac{\psi}{T} \right).
\end{equation}
The convexity of the entropy function guarantees the invertibility of the mapping between the entropy and the conservative variables. 

If the solution is smooth, contracting the entropy variables with the \review{balance} equations \eqref{conservation_law_eq} yields the (scalar) \review{balance} equation for the entropy function,
\begin{equation}
\label{entropy_cons_law}
    \frac{\partial S}{\partial t} + \vec \nabla \cdot \vec f^S =
    \mathbf w^T \left( \vec \nabla \cdot \overset{\leftrightarrow}{\mathbf{f}}\!^\nu \right) + \mathbf w^T \mathbf r,
\end{equation}
where $\vec f^S = S \, \vec v$ is the entropy flux which, as shown in \cite{IdealGLMMHD} for GLM-MHD system, satisfies the identity
\begin{equation}
  \vec \nabla \cdot \vec f^S =
  \mathbf w^T \vec \nabla \cdot \overset{\leftrightarrow}{\mathbf{f}}\!^a + \mathbf w^T \mathbf{\Upsilon}.
\end{equation}

In entropy variables, the diffusive flux
$\overset{\leftrightarrow}{\mathbf f}\!^\nu = (\mathbf f_1^\nu, \mathbf f_2^\nu, \mathbf f_3^\nu)$
takes the following (matrix-vector multiplication) form:
\begin{equation}
\label{entropy_diff_term}
    \mathbf f_i^\nu = \sum_{j=1}^3\mathrm C_{ij}^\nu \frac{\partial \mathbf w}{\partial x_j}.
\end{equation}
The $3\times3$ block matrix, compiled from 9 matrices $C_{ij}^\nu$ of size $9\times9$ each, appearing in \eqref{entropy_diff_term}, is symmetric positive semi-definite.
With the definition \eqref{r_assumption} and the positivity assumption \eqref{eq:positivity_assumption}, the last term in \eqref{entropy_cons_law} equals
\begin{equation}
\label{r_effect_equation}
    \mathbf w^T \mathbf r = -\frac{\alpha \psi^2}{\mu_0 T}
\end{equation}
and is non-positive.

Considering \eqref{entropy_diff_term} and \eqref{r_effect_equation}, integrating \eqref{entropy_cons_law}
over the domain $\Omega$ with boundary $\Gamma$, and invoking the divergence theorem leads to the following scalar global equation for the entropy function:
%\textcolor{red}{refer to divergemnce theorem}
\begin{equation}
\label{entropy_int_eq}
    \frac{d}{dt} \int_\Omega S d\Omega \le
    \oint_{\Gamma} \left(\mathbf w^T \mathbf f^\nu_n - f^S_n \right) d\Gamma -
    \int_\Omega \left(\frac{\partial \mathbf w}{\partial x_i} \right)^T\!\mathrm C_{ij}^\nu
                \left(\frac{\partial \mathbf w}{\partial x_j} \right) d\Omega -
    \int_\Omega \frac{\alpha \psi^2}{\mu_0 T} d\Omega,
\end{equation}
where $\mathbf f^\nu_n$ and $f^S_n$ are the normal components to the boundary of the corresponding vectors.
For a smooth solution, \eqref{entropy_int_eq} takes the form of an equality.
The semi-positiveness of the matrices $C^\nu_{ij}$ and the non-negativeness of the coefficient $\alpha$
guarantee the volume integrals on the right-hand side of inequality \eqref{entropy_int_eq} are non-positive.
Therefore, the entropy in $\Omega$ can only increase through the flux across the boundary~$\Gamma$, \emph{i.e.}, the integral 
term on the right-hand side of \eqref{entropy_int_eq}. This is the same term appearing in the 
entropy stability analysis of the compressible Navier--Stokes. Therefore, the following two theorems presented in
\cite{dalcin2019conservative} and \cite{parsani2015entropy} bound the flux of entropy through a solid wall and, hence, the time derivative
of the entropy function given in \eqref{entropy_int_eq}.

The first theorem is a generalization of Theorem
2.2 from \cite{dalcin2019conservative} (see also Theorem 3.1 in \cite{parsani2015entropy})
from the compressible Navier--Stokes equations to the compressible GLM-MHD system \eqref{conservation_law_eq}.
\begin{thm}
\label{inviscid_entropy_theorem}
At a solid wall, the condition
\begin{equation}
  \label{velocity_bc}
\vec v = \vec v^\text{\,w}, \quad \vec v^\text{\,w} \cdot \vec n = 0, \quad \forall \vec x \in \Gamma, \quad t \ge 0,
\end{equation}
where $\vec v^\text{\,w}$ represents the wall velocity and $\vec n$ is the unit outward normal vector to the wall,
bounds the advective contribution to the time
derivative of the entropy in equation \eqref{entropy_int_eq}.
\end{thm}
\begin{pf}
The normal advective entropy flux $f^S_n = S v_n$ in equation \eqref{entropy_int_eq} equals zero at the wall. $\square$
\end{pf}
The next theorem is an extension of Theorem 3.2 from \cite{parsani2015entropy}
from the compressible Navier--Stokes equations to the GLM-MHD system \eqref{conservation_law_eq}.
\begin{thm}
\label{viscous_entropy_theorem}
The wall boundary condition (or prescription of a heat entropy flux, also known as heat entropy
transfer) 
\begin{equation}
  \label{thermal_bc}
g(t) = \kappa \frac{1}{T} \frac{\partial T}{\partial n}, \quad \forall \vec x \in \Gamma, \quad t \ge 0,
\end{equation}
where $\partial T/\partial n$ denotes the temperature gradient normal to the wall,
bounds the diffusive contribution to the time derivative of the entropy in \eqref{entropy_int_eq}.
\end{thm}
\begin{pf}
    Contraction of the normal diffusive flux, $\mathbf f^\nu_n$, with the entropy variables, $\mathbf w$, reads, 
    %\textcolor{red}{Is the first term correct? Why the 2nd and third terms are zero?}
\begin{align*}
    \mathbf w^T \mathbf f^\nu_n &=  \frac{\vec v \cdot \vec{\tau}_n}{T} - \frac{1}{T}
    \left( \vec v \cdot \vec{\tau}_n - \kappa \frac{\partial T}{\partial n} - 
    \frac{\mu_R}{\mu_0^2} \left(\left(\vec \nabla \times \vec B\right) \times \vec B\right)_n \right) + 
    \frac{\vec B}{\mu_0 T} \frac{\mu_R}{\mu_0}  \left(\vec \nabla \vec B - \left(\vec \nabla \vec B\right)^T\right)_n \\
    &= \kappa \frac{\partial T}{\partial n}\frac{1}{T} = g(t),
\end{align*}
    where the subscript $n$ denotes \review{the component normal} to the wall of the corresponding tensor,
    \review{and all the terms involving the magnetic field $\vec B$ cancel trivially from the vector calculus identity
    $\vec C \times (\vec \nabla \times \vec D)  = \vec C \cdot (\vec \nabla \vec D - (\vec \nabla \vec D)^T)$}.
    Therefore, the prescribed boundary entropy flux (or entropy transfer) function $g(t)$ 
    bounds the diffusive contribution to the time derivative of the entropy in \eqref{entropy_int_eq}. $\square$
\end{pf}

\begin{remark}
As outlined in \cite{parsani2015entropy}, this boundary condition at first 
glance appears ``uncommon". However, the scalar value $\kappa \frac{1}{T} \frac{\partial T}{\partial n}$
    accounts for the change in entropy at the boundary and matches the result of the thermodynamic (entropy) analysis
    of a generic system~\cite{bejan1995entropy}. In fact, the compressible Navier--Stokes model encapsulates correctly the
thermodynamic analysis at the boundary of the domain, $\Omega$.

It is essential to understand that the boundary term $\kappa \frac{1}{T} \frac{\partial T}{\partial n}$ 
is not a peculiar term only for a wall \cite{parsani2015entropy}. 
It is the contribution of a generic piece of boundary of $\Gamma$ to the entropy analysis. Thus, the development of any entropy 
conservative and entropy stable boundary conditions for the compressible
Navier--Stokes equations and the GLM-MHD system \eqref{conservation_law_eq} have to provide a 
    proper treatment for it.
\end{remark}

\subsection{Solid wall boundary conditions}
\label{continuous_bc_section}

In the previous subsection, we have shown that when considering a solid wall, the
entropy analysis of the GLM-MHD system \eqref{conservation_law_eq} yields entropy 
conservative and entropy stable boundary conditions, which are the same as those 
obtained for the compressible Navier--Stokes equations \cite{parsani2015entropy,dalcin2019conservative}. 
Thus, the entropy analysis of the GLM-MHD system \eqref{conservation_law_eq} does not provide information on if and how boundary conditions
must be specified for the magnetic field $\vec B$ and the generalized Lagrange 
multiplier $\psi$. In this work, we develop new discrete entropy 
conservative and entropy stable solid wall boundary conditions for the magnetic field, which mimic the 
solid wall boundary conditions for the resistive MHD equations \cite{davidson2002resistbook}.  

Magnetic wall boundary conditions come in different forms depending on the conductive properties of the wall.
 We refer the reader to the book of Davidson \cite{davidson2002resistbook} for a detailed discussion and derivation of magnetic boundary conditions for
electrically (perfectly) insulating walls, electrically conducting thin walls,
and electrically perfectly conducting walls.

For electrically insulating walls, the magnetic field is prescribed at a boundary $\Gamma^\text{ins}$ as
\begin{equation}
\label{insulating_cond}
\vec B = \Bextvec, \quad \forall \vec x \in \Gamma^\text{ins}, \quad t \ge 0,
\end{equation}
where $\Bextvec$ is the external magnetic field.
This condition follows from the assumption that no electric current flows from the fluid to the wall or within the wall.
Consequently, there is no jump in the magnetic field across the fluid-solid interface. 
Since no electric current flows through the wall, an additional condition on the external magnetic
field is $(\vec \nabla \times \Bextvec) \cdot \vec n = 0$.

For electrically conducting walls, the magnetic wall boundary condition involves the continuity of the wall-normal magnetic field
and the wall-tangent electric field at the fluid-solid interface. By the Ohm's and Amp\`ere's laws, and the thin wall approximation,
the magnetic boundary condition at a boundary $\Gamma^\text{cond}$ reads 
\begin{equation}
\label{conductive_cond}
\vec B \cdot \vec n = \Bextvec  \cdot \vec n, \quad 
\left(\vec \nabla \times \vec B\right) \times \vec n = c_d^{-1} \left(\Bextvec - \vec B\right)
, \quad \forall \vec x \in \Gamma^\text{cond},\quad t \ge 0,
\end{equation}
where $c_d = (\mu_R / \mu_\text{w}) d_\text{w} $ is a wall conductance coefficient,
$\mu_R$ is the fluid electrical resistivity introduced in \eqref{eq:diffusive-flux-mhd},
$\mu_\text{w}$ denotes the wall resistivity, and $d_\text{w}$ denotes the wall thickness.
The wall is assumed to be thin, that is, $d_\text{w} \ll L^*$, where $L^*$ is the characteristic length of the domain.
The wall conductivity properties are usually expressed in terms of the non-dimensional wall conductance parameter,
which is defined as
\begin{equation}
\label{wall_cond_param}
c = \frac{\mu_R \, d_\text{w}}{\mu_\text{w} \, L^*},
\end{equation}
and, therefore, this non-dimensional quantity expresses the ratio of wall and fluid conductances.
For perfectly conducting walls where $c \rightarrow \infty$, the magnetic boundary condition \eqref{conductive_cond} reduces to
\begin{equation}
\label{conductive_perf_cond}
\vec B \cdot \vec n = \Bextvec \cdot \vec n, \quad 
\left(\vec \nabla \times \vec B\right) \times \vec n = \vec 0
, \quad \forall \vec x \in \Gamma^\text{cond},\quad t \ge 0,
\end{equation}
and it is valid for walls of any shape and dimension.

For the ideal MHD equations (\emph{i.e.}, system \eqref{conservation_law_eq} with 
$\overset{\leftrightarrow}{\mathbf{f}^\nu} =\overset{\leftrightarrow}{\mathbf{0}}$), the wall-normal magnetic field should be zero, 
\emph{i.e.}, $\Bextvec  \cdot \vec n = 0$, and the condition on the tangential component
of the magnetic field is redundant \cite{goedbloed2004idealbook}. In the 
resistive MHD equations, the wall-normal magnetic field can be nonzero, 
but the total magnetic flux through any closed surface must be zero due to the
divergence-free constraint on the magnetic field. A nonzero wall-normal magnetic 
field forms a Hartmann boundary layer, which usually requires 
appropriate mesh refinement for an accurate numerical solution \cite{davidson2002resistbook}.

The GLM term $\psi$ does not require additional boundary conditions at a solid wall.
The rationale for this claim is twofold. First, there is no diffusion of $\psi$ in the system.
Second, $\psi$ acts as an artificial pressure controlling the divergence of the magnetic field, as described in \cite{IdealGLMMHD}.
%In analogy with the compressible Euler equations, where no pressure or density/temperature needs to be specified on the wall,
%$\psi$ does not require any special treatment at solid wall boundaries.

\subsection{One-dimensional formulation}
Following the approach of \cite{DGSEM_MHD_Part2},
we reproduce here the system of equations \eqref{conservation_law_eq} in a one-dimensional case.
Later on, and for the sake of simplicity, the numerical method will also be formulated in one space dimension.

In one dimension, the GLM-MHD system \eqref{conservation_law_eq} reads
\begin{equation}
\label{conserv_eq_1d}
    \partial_t \mathbf u +
    \partial_x \mathbf f^a (\mathbf u) +
    \mathbf \Upsilon (\mathbf u, \vec \nabla \mathbf u) =
    \partial_x \mathbf f^\nu (\mathbf u, \vec \nabla \mathbf u) +
    \mathbf r.
\end{equation}
Following the notation in \cite{DGSEM_MHD_Part2}, we omit the sub-index in the conservative fluxes,
$\mathbf f^a$ and $\mathbf f^\nu$, and in the coordinate, $x$, to simplify the notation and enhance 
readability. Thus, we use $\mathbf f^a$,  $\mathbf f^\nu$, and $x$ instead of 
$\mathbf f^a_1$, $\mathbf f^\nu_1$, and $x_1$, respectively.
The advective flux in the $x$ direction, split into Euler, MHD, and GLM parts, is given by
\begin{equation}
\label{invisced_flux_eq}
    \mathbf f^a =
    \begin{pmatrix}
    \rho v_1 \\
    \rho v_1^2 + p \\
    \rho v_1 v_2 \\
    \rho v_1 v_3 \\
    v_1 \left( \frac{1}{2}\rho|\!|\vec v|\!|^2 +
    \frac{\gamma p}{\gamma-1}\right) \\
    0 \\
    0 \\
    0 \\
    0
    \end{pmatrix}
    +
    \begin{pmatrix}
    0 \\
    \frac{1}{\mu_0}\left(\frac{1}{2}|\!|\vec B|\!|^2-B_1^2\right) \\
    - B_1 B_2 / \mu_0 \\
    - B_1 B_3 / \mu_0 \\
    \frac{1}{\mu_0} \left( v_1 |\!|\vec B|\!|^2 -
    B_1 \left( \vec v \cdot \vec B \right) \right)\\
    0 \\
    v_1 B_2 - v_2 B_1 \\
    v_1 B_3 - v_3 B_1 \\
    0
    \end{pmatrix}
    +
    \begin{pmatrix}
    0 \\
    0 \\
    0 \\
    0 \\
    \frac{c_h}{\mu_0} \psi B_1 \\
    c_h \psi \\
    0 \\
    0 \\
    c_h B_1
    \end{pmatrix}.
\end{equation}

The non-conservative term, $\mathbf \Upsilon = \mathbf \Upsilon^\text{MHD} + \mathbf \Upsilon^\text{GLM}$,
consists of the following two parts:
\begin{equation}
\label{Upsilon_eq}
    \mathbf \Upsilon^{\text{MHD}} =
    \frac{\partial B_1}{\partial x} \pmb{\phi}^\text{MHD} =
    \frac{\partial B_1}{\partial x}
    \begin{pmatrix}
        0 \\
        B_1/\mu_0 \\
        B_2/\mu_0 \\
        B_3/\mu_0 \\
        \vec v \cdot \vec B/\mu_0 \\
        v_1 \\
        v_2 \\
        v_3 \\
        0
    \end{pmatrix},\quad
    \mathbf \Upsilon^\text{GLM} =
    \frac{\partial \psi}{\partial x} \pmb{\phi}^\text{GLM} =
    \frac{\partial \psi}{\partial x}
    \begin{pmatrix}
        0 \\
        0 \\
        0 \\
        0 \\
        v_1 \psi/\mu_0 \\
        0 \\
        0 \\
        0 \\
        v_1
    \end{pmatrix}.
\end{equation}
%{\color{red} To delete: mistype in \cite{DGSEM_MHD_Part2} was fixed. SEND EMAIL TO THE AUTHORS NOTIFYING THE TYPO.}

The visco-resistive flux in the $x$ direction reads
\begin{equation}
\label{viscous_flux_eq}
    \mathbf f^\nu \left( \mathbf u, \vec \nabla \mathbf u \right) =
    \begin{pmatrix}
        0 \\
        \frac{4}{3} \mu_{NS} \px{v_1} \\
        \mu_{NS} \px{v_2} \\
        \mu_{NS} \px{v_3} \\
        \mu_{NS} \left( \frac{4}{3} v_1 \px{v_1} +
        v_2 \px{v_2} + v_3 \px{v_3} \right) + \kappa \px{T} +
        \frac{\mu_R}{\mu_0^2} \left( B_2 \px{B_2} + B_3 \px{B_3} \right) \\
        0 \\
        \frac{\mu_R}{\mu_0} \px{B_2} \\
        \frac{\mu_R}{\mu_0} \px{B_3} \\
        0
    \end{pmatrix}.
\end{equation}

%As in the previous subsection, reaction term is defined as $\mathbf r = (0,\vec 0,0,\vec 0,-\alpha\psi)^T$, with $\alpha \ge 0$.

\section{Entropy stable DG-SBP-SAT discretization}
\label{discret_section_label}
This section presents the new discrete entropy conservative and entropy stable solid wall boundary conditions for the
GLM-MHD system \eqref{conservation_law_eq}.
These conditions are compatible with the collocated entropy conservative and entropy stable discontinuous Galerkin (DG) discretization schemes with the summation-by-parts (SBP) property and coupled 
with simultaneously approximated terms (SAT), (\emph{i.e.}, DG-SBP-SAT
operators). The DG-SBP-SAT method was developed and proven to be entropy
stable for the GLM-MHD system in
unbounded computational domains in \cite{DGSEM_MHD_Part1}. Thus, this work extends the 
entropy stability proof to include solid wall boundary conditions which can also
be ``trivially" extended to modal DG-SBP operators.

First, we revisit the method's definition and entropy stability proof for nodes collocated inside the domain $\Omega$.
Then, we introduce the new solid wall boundary conditions and their original proof for the entropy conservation and 
stability for boundary nodes at a solid wall.
We follow the works presented in \cite{parsani2015entropy, dalcin2019conservative}, where discrete entropy-stable solid
wall boundary conditions are formulated, verified, and validated for the compressible Navier--Stokes equations.
Deviating from the original notation of \cite{DGSEM_MHD_Part1},
we introduce a nonsymmetric two-point flux function to address the nonconservative term in the MHD equations.
Technically, this approach arises from rearranging terms in the previously introduced expressions.
Such rearrangement significantly enhances the clarity of our presentation.

Following the approach outlined in \cite{DGSEM_MHD_Part2},
all definitions and proofs are presented for the one-dimensional (1D) case,
assuming that all variables vary solely along one spatial coordinate and are uniform in the other two.
While the entropy stability of the proposed discrete boundary conditions are
demonstrated in 1D, their formulation is inherently nodal. Therefore, they are applicable across various scenarios where DGSEM is used,
including unstructured curvilinear nonconforming meshes with $hp$-refinement 
\cite{friedrich2018entropy,fernandez2020entropy}.

\subsection{DG-SBP-SAT discretization in 1D}
To present the DG-SBP-SAT discretization, we rewrite the system of equations
\eqref{conserv_eq_1d} as a first-order system:
\begin{equation}
\label{splitted_conserv_eq}
\begin{split}
    & \partial_t \mathbf u + \partial_x \mathbf f^a (\mathbf u) + \mathbf \Upsilon (\mathbf u, \partial_x \mathbf u) = \partial_x \mathbf f^\nu (\mathbf u, \mathbf g) + \mathbf r(\mathbf u), \\
    & \mathbf g = \partial_x \mathbf w,
\end{split}
\end{equation}
where $\mathbf g$ denotes the spatial gradient of the entropy variables, $\mathbf w$.

In the case of a 1D problem, the computational mesh is composed of non-overlapping intervals. 
Each interval is mapped into a reference space, $\xi \in [-1,1]$. 
Within each element, all variables are represented using Lagrange interpolating polynomials of degree $N$, 
denoted as $\{l_k\}_{k=1}^N$, which are defined on $N+1$ Legendre--Gauss--Lobatto (LGL) nodes, $\{\xi_k\}_{k=0}^N$. 
The polynomial basis is continuous within each element but discontinuous at the interfaces between elements. 
A key aspect for ensuring the entropy conservation property of the scheme is that the undivided derivative matrix, 
$Q_{jk} = \omega_j D_{jk} = \omega_j l_k^\prime(\xi_j)$, has the SBP property \cite{carpenter2014entropy,carpenter2015entropy},
\emph{i.e.},
\begin{equation}
\label{SBP_def}
    Q + Q^T = \mathcal{B},
\end{equation}
where $\omega_j$ represents the reference space quadrature weight, $D_{jk}$ is
the $jk$ element of the one-dimensional (SBP) operator for the first
derivative, $l_k^\prime$ indicates the derivative of the $k$-th Lagrange
interpolating polynomial, and $\mathcal{B}$ is a boundary operator defined as
\begin{eqnarray}
\label{B_def}
    \mathcal{B} = \diag(-1,0,\dots,0,1).
\end{eqnarray}
Boundary conditions and conditions at the interfaces between mesh elements are incorporated into the semi-discretization
in a weak sense using the simultaneous approximation terms (SATs) approach;
see the original SAT works in \cite{Carpenter1994,Carpenter1999} and its
extension to the compressible Navier--Stokes in the context of entropy
stability \cite{parsani2015entropyInterfaces}. More details on 
the method and its extension to multidimensional cases, including the
$h/p$-refinement, can be found in \cite{DGSEM_MHD_Part1, parsani2021high, fernandez2020entropy}.

For the node $j$ in a mesh element, the spatial discretization leads to the
following semi-discrete form (\emph{i.e.}, an ordinary differential equation for the unknown vector ${\mathbf u}_j$),
\begin{equation}
\label{discrete_eq}
    J_j \omega_j \frac{d {\mathbf u}_j}{d t} = -\mathbf F_j^{a} + \mathbf F_j^{\nu} + J_j \omega_j \mathbf r_j,
\end{equation}
where $J_j$ is the determinant of the Jacobian of the transformation from the reference space to the physical space,
$\mathbf F^{a}_j$ represents the discretization of the advective flux divergence and non-conservative terms combined,
$\mathbf F^{\nu}_j$ represents the discretization of the diffusive flux divergence, and $\mathbf r_j = \mathbf r(\mathbf u_j)$.
Specifically, using a DG-SBP-SAT operator, the discretization of the divergence
of the convective flux and the
non-conservative term yields
\begin{equation}
\label{advect_discrete_def}
    \mathbf F_j^{a} = 2\sum_{k=0}^N Q_{jk} \mathbf f^*_{(j,k)} +
    \delta_{j0} \left( \mathbf f^*_{(0,0)} - {\mathbf f}^*_{(0,L)} - \text{DISS}^a_{(0,L)} \right) -
    \delta_{jN} \left( \mathbf f^*_{(N,N)} - {\mathbf f}^*_{(N,R)} + \text{DISS}^a_{(N,R)} \right).
\end{equation}
Here, the index $k$ ranging from $0$ to $N$ corresponds to the index of a node within the mesh element. 
The subscripts $L$ and $R$ refer to the nodes to the ``left" and the ``right"
interface of the ``left" and ``right" adjacent elements, respectively. These nodes are
collocated at the same physical location as those with index $k = 0$ and
$k= N$ within the element under consideration and are used to construct the SAT interface terms. %outside of the element but aligned with the $0$ and $N$ nodes, respectively. For mesh elements adjacent to the boundary, these $L$ or $R$ nodes belong to boundary (ghost) nodes. Different values of the variables at nodes adjacent to the interface on the left and right sides indicate discontinuities in the solution at the interfaces. The Kronecker delta function $\delta_{ij}$ precedes the SAT terms, which are active only at the mesh element interfaces.

At this point of the section, we depart from the notation used in \cite{DGSEM_MHD_Part1, DGSEM_MHD_Part2} and introduce a non-symmetric numerical two-point flux function, denoted as $\mathbf f^*_{(j,k)} = \mathbf f^*(\mathbf u_j, \mathbf u_k)$. 
This function is derived by rearranging different terms in the expressions from \cite{DGSEM_MHD_Part2} (specifically, Eq.~32). 
This step indirectly accounts for the non-conservative term $\mathbf \Upsilon$ and simplifies the expressions that will follow. 
The non-symmetric two-point flux is defined as
\begin{equation}
\label{nonsym_flux_def}
    \mathbf f^*(\mathbf u_i, \mathbf u_j) = \mathbf f^a(\mathbf u_i, \mathbf u_j) + \frac{1}{2} 
    \left( \pmb{\phi}_i^\text{MHD} B_{1,j} + \pmb{\phi}_{i}^\text{GLM} \psi_j \right),
\end{equation}
where $\mathbf f^a(\mathbf u_i, \mathbf u_j)$ represents any symmetric two-point flux function. 
In our simulations, we use the entropy conservative two-point flux function 
$\mathbf f^\text{EC}(\mathbf u_i, \mathbf u_j)$ designed in \cite{IdealGLMMHD}. 
This two-point flux function is reproduced in \ref{two_point_flux_app} for completeness.
The interface dissipation term $\text{DISS}^a_{(i,j)}$ in \eqref{advect_discrete_def} 
introduces an additional penalty on the solution jump at the mesh element interfaces.
The particular form of $\text{DISS}^a_{(i,j)}$ is specified in Section~\ref{discrete_advect_bc_section}.

%In the definition of $\hat{\mathbf f}^*_{(i, j)}$, there is an option to set up an additional penalty 
%on the solution jump at the mesh element interfaces,
%\begin{equation}
%    \hat{\mathbf f}^*_{(i, j)} = \mathbf{f}^*_{(i,j)} - \text{DISS}^a_{(i,j)},
%\end{equation}
%where $\text{DISS}^a_{(i,j)}$ is a dissipation operator that is symmetric with respect to its parameters. 
%The particular form of the dissipation term $\text{DISS}^a_{(i,j)}$ will be specified later in section~\ref{discrete_advect_bc_section}.

The diffusive term is discretized using a standard SBP differentiation operator
constructed at the $N+1$ LGL points. After some algebraic manipulations and
simplifications, the resulting expression for the approximation of the divergence of the diffusive
flux reads
\begin{equation}
\label{visc_discrete_def}
    \mathbf F_j^{\nu} = \sum_{k=0}^N Q_{jk} \mathbf f^\nu_k +
    \delta_{j0} \left( \mathbf f^\nu_0 - \hat{\mathbf f}^\nu_{(0,L)} + \text{DISS}^\nu_{(0,L)} \right) -
    \delta_{jN} \left( \mathbf f^\nu_N - \hat{\mathbf f}^\nu_{(N,R)} - \text{DISS}^\nu_{(N,R)} \right),
\end{equation}
where $\mathbf f^\nu_k = C^\nu(\mathbf u_k) \mathbf g_k$ represents the diffusive flux at node $k$.
The discrete gradient of entropy variables, $\mathbf g_k$, is determined using a consistent discretization approach,
\begin{equation}
\label{g_discrete_def}
    J_k \omega_k \mathbf g_k = \sum_{n=0}^N Q_{kn} \mathbf w_n +
    \delta_{kN} \left( \hat{\mathbf w}_{(N,R)} - \mathbf w_N \right) - 
    \delta_{k0} \left( \hat{\mathbf w}_{(0,L)} - \mathbf w_0 \right).
  \end{equation}
The diffusive numerical flux $\hat{\mathbf f}^\nu_{(i,j)}$ in \eqref{visc_discrete_def}
and the numerical entropy variables $\hat{\mathbf w}_{(i,j)}$ in \eqref{g_discrete_def}
are computed using the BR1 method \cite{BR1_1997},
\begin{equation}
\label{BR1_eq}
    \hat{\mathbf f}^\nu_{(i,j)} = \frac{1}{2} \left( \mathbf f^\nu_i + \mathbf f^\nu_j \right), \quad %+ \text{DISS}^\nu_{(i,j)}, \quad
    \hat{\mathbf w}_{(i,j)} = \frac{1}{2} \left( \mathbf w_i + \mathbf w_j \right).
\end{equation}
The visco-resistive interface dissipation $\text{DISS}^\nu_{(i,j)}$ is defined in Section~\ref{internal_conservation_section}.

\subsection{Entropy conservation for internal nodes}
\label{internal_conservation_section}

In \cite{DGSEM_MHD_Part1}, it is shown that the DG-SBP-SAT method described
previously,
does not produce or dissipate numerical entropy in internal nodes, \emph{i.e.},
it is entropy conservative when $\text{DISS}^a = \text{DISS}^\nu =0$, and $\alpha = 0$.\footnote{By internal nodes,
we mean nodes lying inside mesh elements or pairs of nodes lying on interfaces between adjacent mesh elements.}
This feature can be verified numerically by simulating problems with smooth
solutions and periodic computational domain.

The proof of entropy conservation of the interior DG-SBP-SAT operator relies on the generalized Tadmor's condition \cite{DGSEM_MHD_Part1}(Eq.~4.8) and \cite{DGSEM_MHD_Part2}(Eq.~51),
which accounts for the non-conservative contributions. In our notation, the generalized Tadmor's condition reads
\begin{equation}
    \label{gen_tadmors_cond}
    \mathbf w_k^T \mathbf f^*_{(k,j)} - \mathbf w_j^T \mathbf f^*_{(j,k)} = \Psi^*_k - \Psi^*_j, \quad
    \Psi^*_k = \mathbf w_k^T \mathbf f^*_{(k,k)} - f^S_k,
\end{equation}
where $f^S_k$ is the entropy flux at the node $k$.
%Direct substitution shows that
%\begin{equation}
%    \Psi^*_k = \mathbf w_k^T \mathbf f^a_k - f^S_k + \frac{1}{2}\theta B_1 =
%    \frac{v_1}{T} \left( p + \frac{|\!|\vec B|\!|^2}{2} \right) +
%    \frac{B_1}{T} \left(c_h \psi - \frac{(\vec v \cdot \vec B)}{2} \right).
%\end{equation}
For completeness, we repeat the proof of entropy conservation here. Using condition \eqref{gen_tadmors_cond} and the SBP property of the differentiation
operator encapsulated in \eqref{SBP_def}, we next show that the numerical scheme is entropy conservative within a single mesh element.
Mimicking the continuous analysis presented in Section
\ref{sec:eq-thermo}, we contract \eqref{discrete_eq} with the entropy
variables, $\mathbf w$, 
on a single mesh element and after some algebraic manipulations and
simplifications, we arrive at the following expression %(see \ref{appendix_num_ss} for a detailed derivation) 
\begin{equation}
    \label{entropy_discrete_eq}
    \begin{split}
    &\frac{d }{d t} \sum_{j=0}^N J_j \omega_j S_j +
    \sum_{j=0}^N J_j \omega_j \mathbf g_j^T \mathrm C^\nu_j \mathbf g_j +
    \sum_{j=0}^N J_j \omega_j \frac{\alpha \psi_j^2}{\mu_0 T} = \\
    &P^\text{a,cons}_{(0,L)} + P^\text{a,diss}_{(0,L)} + P^{\nu,\text{cons}}_{(0,L)} + P^{\nu,\text{diss}}_{(0,L)} +
     P^\text{a,cons}_{(N,R)} + P^\text{a,diss}_{(N,R)}
      + P^{\nu,\text{cons}}_{(N,R)} + P^{\nu,\text{diss}}_{(N,R)},
    \end{split}
\end{equation}
where $S_j$ is the entropy function at the node $j$. Therefore, the
discrete entropy within a mesh element changes solely due to the dissipation
(second term on the left-hand side),
the algebraic damping on $\psi$ (third term on the left-hand side),
and the entropy flow across the boundaries of the mesh element (all terms on the right-hand side). 
In the following expressions, we separate the advective and visco-resistive entropy conservative and dissipative contributions
at the left and right boundaries of the mesh element,
\begin{align}
    \label{p_a_cons_L}
    P^\text{a,cons}_{(0,L)} &= \left(\mathbf w_0^T \mathbf f^*_{(0,L)} - \Psi^*_0\right), \\
    \label{p_a_cons_R}
    P^\text{a,cons}_{(N,R)} &= -\left(\mathbf w_N^T \mathbf f^*_{(N,R)} - \Psi^*_N\right), \\
    \label{p_nu_cons_L}
    P^{\nu,\text{cons}}_{(0,L)} &= - \left(\hat{\mathbf w}_{(0,L)} - \mathbf w_0 \right)^T \mathbf f^\nu_0 - \mathbf w_0^T \hat{\mathbf f}^\nu_{(0,L)}, \\
    \label{p_nu_cons_R}
    P^{\nu,\text{cons}}_{(N,R)} &= \left(\hat{\mathbf w}_{(N,R)} - \mathbf w_N \right)^T \mathbf f^\nu_N + \mathbf w_N^T \hat{\mathbf f}^\nu_{(N,R)}, \\
    \label{p_a_diss_L}
    P^\text{a,diss}_{(0,L)} &= \mathbf w_0^T \text{DISS}^a_{(0,L)}, \\
    \label{p_a_diss_R}
    P^\text{a,diss}_{(N,R)} &= \mathbf w_N^T \text{DISS}^a_{(N,R)}, \\
    \label{p_nu_diss_L}
    P^{\nu,\text{diss}}_{(0,L)} &= \mathbf w_0^T \text{DISS}^\nu_{(0,L)}, \\
    \label{p_nu_diss_R}
    P^{\nu,\text{diss}}_{(N,R)} &= \mathbf w_N^T \text{DISS}^\nu_{(N,R)}.
\end{align}

Next, we demonstrate that the appropriate choice of SAT terms in
\eqref{p_a_cons_L} through \eqref{p_nu_diss_R} guarantees 
entropy dissipation at the interfaces between adjacent mesh elements.
The scalar equation describing the rate of change of the entropy function in
the whole domain is
obtained by adding contributions from \eqref{entropy_discrete_eq} over all the mesh elements.
\review{For any interface between two neighboring mesh elements, there exists a node $\nodeA$ with index $N$ in the left mesh element
and a node $\nodeB$ with index $0$ in the right mesh element that are spatially coincident.}
We will separately examine the advective and diffusive fluxes to show that for every $\nodeA$ and $\nodeB$ node pair,
the total entropy production is zero for conservative fluxes and non-positive
for entropy dissipative fluxes.

Terms \eqref{p_a_cons_L} and \eqref{p_a_cons_R} are the contributions to the 
time rate of change of the entropy function due to the advective SATs for the
\review{interface nodes $\nodeA$ and $\nodeB$}. Their cumulative contribution to
$dS/dt$ is zero, as can be shown by using the generalized Tadmor's condition
\eqref{gen_tadmors_cond}:
\begin{align*}
    \frac{dS^\text{a,cons}_{\nodeA+\nodeB}}{dt} =
    \left( \mathbf w_\nodeB^T {\mathbf f}^*_{(\nodeB,\nodeA)} - \Psi^*_\nodeB \right) -
    \left( \mathbf w_\nodeA^T {\mathbf f}^*_{(\nodeA,\nodeB)} - \Psi^*_\nodeA \right) =
    \left( \mathbf w_\nodeB^T \mathbf f^*_{(\nodeB,\nodeA)} - \mathbf w_\nodeA^T \mathbf f^*_{(\nodeA,\nodeB)} \right) -
    \left( \Psi^*_\nodeB - \Psi^*_\nodeA \right) = 0.
\end{align*}
The dissipation operator $\text{DISS}^a$
is antisymmetric in its arguments and can be expressed as a linear operator on
the jump in the entropy variables,
\begin{equation}
    \label{diss_a_definition}
    \text{DISS}^a_{(i,j)} = \hat D [\![\mathbf w]\!]_{(i,j)},
\end{equation}
where $\hat D$ is a symmetric positive semi-definite matrix. 
Thus, with the definition \eqref{diss_a_definition}, the (numerical)
contribution to the time rate of change of the entropy function due to
the terms \eqref{p_a_diss_L} and \eqref{p_a_diss_R} is given by
\begin{align*}
    \frac{dS^\text{a,diss}_{\nodeA+\nodeB}}{dt} =
    \mathbf w_\nodeB^T \text{DISS}^a_{(\nodeB,\nodeA)} + \mathbf w_\nodeA^T \text{DISS}^a_{(\nodeA,\nodeB)} =
    - [\![\mathbf w]\!]_{(\nodeA,\nodeB)}^T \hat D [\![\mathbf w]\!]_{(\nodeA,\nodeB)},
\end{align*}
which is negative due to the positive semi-definiteness of $\hat D$.

Similar to the above reasoning, using the expressions given in \eqref{BR1_eq}, 
we can show that the cumulative contribution to the time derivative of the
entropy function due to the entropy conservative visco-resistive
part of the flux expressed by \eqref{p_nu_cons_L} and \eqref{p_nu_cons_R} is zero,
\begin{equation}
    \label{cons_entropy_term}
    \frac{dS^{\nu,\text{cons}}_{\nodeA+\nodeB}}{dt} =
    \left( \hat{\mathbf w}_{(\nodeA,\nodeB)} - \mathbf w_\nodeA \right)^T \mathbf f^\nu_\nodeA -
    \left( \hat{\mathbf w}_{(\nodeB,\nodeA)} - \mathbf w_\nodeB \right)^T \mathbf f^\nu_\nodeB +
    \mathbf w_\nodeA^T \hat{\mathbf f}^\nu_{(\nodeA,\nodeB)} - \mathbf w_\nodeB^T \hat{\mathbf
    f}^\nu_{(\nodeB,\nodeA)} = 0.
\end{equation}

Finally, we consider the contribution of the visco-resistive operator $\text{DISS}^\nu$, which reads
\begin{equation}
    \label{visc_diss_interface_term}
    \text{DISS}^\nu_{(i,j)} = \hat L [\![\mathbf w]\!]_{(i,j)}, \quad
    \hat L = \beta \frac{C^\nu_{i} + C^\nu_{j}}{2}, \quad \beta \ge 0,
\end{equation}
where $\hat L$ is a scaled average matrix between the two symmetric positive semi-definite dissipation matrices
$C^\nu_i = C^\nu_{11}(\mathbf u_i)$ and $C^\nu_j = C^\nu_{11}(\mathbf u_j)$ defined in \eqref{entropy_diff_term}.
The nonnegative coefficient $\beta$ modulates the strength of the dissipation term.
Therefore, the contribution to the time rate of change of the entropy function
due to the (numerical) visco-resistive dissipative mechanisms is 
\begin{equation*}
    \frac{dS^{\nu,\text{diss}}_{\nodeA+\nodeB}}{dt} =
    \mathbf w_\nodeA^T \text{DISS}^\nu_{(\nodeA,\nodeB)} + \mathbf w_\nodeB^T \text{DISS}^\nu_{(\nodeB,\nodeA)} =
    - [\![\mathbf w]\!]^T_{(\nodeA,\nodeB)} \hat L [\![\mathbf w]\!]_{(\nodeA,\nodeB)},
\end{equation*}
which is nonpositive due to the definition of $\hat L$.

%This {\color{red}passage} finalizes the proof of entropy stability properties of the scheme for internal nodes.
In the following section, we introduce the discrete boundary conditions
for the advective and the visco-resistive terms,
which guarantees a physically consistent production of the discrete entropy
function in a bounded domain.

\section{Entropy conservative solid wall boundary conditions}
We consider boundary nodes within the computational domain corresponding to the internal fluid state and denote any of such nodes as
$\Mnode$.\footnote{At the boundary, we switch from the notation $L$, $R$ to $\Mnode$, $\Pnode$ to define internal and corresponding ghost external nodes.}
Each node $\Mnode$ lies (conceptually) on one side of the element boundary face and has a corresponding virtual or ghost node, denoted as $\Pnode$,
which spatially coincides with node $\Mnode$ but lies (conceptually) on the opposite side of the face. 
Our objective is to define manufactured quantities at the ghost node $\Pnode$ as a function of the discrete 
state and its gradient at node $\Mnode$ and the boundary conditions data 
such that i) the resulting boundary conditions imposition is entropy conservative or entropy stable and 
ii) the same discretization procedure used to treat internal interfaces between mesh elements
can be used to impose boundary conditions. The latter objective has also been one of the guiding principles 
in the construction of the entropy conservative and entropy stable solid wall boundary conditions
presented in \cite{dalcin2019conservative}.   

\review{Similar} to the approach proposed in \cite{parsani2015entropy,dalcin2019conservative} for
the compressible Navier--Stokes equations, we must ensure 
that the discrete imposition of the solid wall boundary conditions for the GLM-MHD system \eqref{conservation_law_eq} 
leads to an appropriate and physically consistent entropy production or dissipation.
%\textcolor{red}{Explain that we move from L and R to + and -} \textcolor{blue}{- VP: a footnode added}

\subsection{Advective numerical flux divergence operator}
\label{discrete_advect_bc_section}
%\textcolor{red}{MP: I believe the next sentence is wrong} \textcolor{blue}{- VP: updated}

Herein, boundary conditions are formulated in primitive variables, defined by letter $\mathbf v$. 
The following theorem defines the manufactured state $\mathbf v^a_{\Pnode}$
to be substituted into the advective numerical flux
$\mathbf f^*(\mathbf v_{\Mnode}, \mathbf v_{\Pnode})$ and advective dissipation
$\text{DISS}^a(\mathbf v_{\Mnode}, \mathbf v_{\Pnode})$
in expression \eqref{advect_discrete_def}.
We define the manufactured state in the proper and general three-dimensional setting,
while the proof is carried out under the one-directional assumption for simplicity.
In our notation, we omit the subscript $\Mnode$ on the individual components of the state
$\mathbf v_{\Mnode} = (\rho, \vec v, T, \vec B, \psi)^T$, to improve readability. 
\begin{thm}
    \label{advective_conservative_flux_theorem}
The advective flux penalty term, disregarding any contributions from the advective dissipation operator $\mathrm{DISS}^a$,
    is entropy-conservative if the vector boundary state is defined in primitive variables as
    \begin{equation}
        \label{inviscid_boundary_conditions}
        \mathbf v^a_{\Pnode} = \begin{bmatrix}
          \rho \\
          \vec v - 2 \left(\vec v \cdot \vec n\right) \vec n \\
          T \\
          \vec B - 2 \left(\vec B \cdot \vec n - \Bextvec \cdot \vec n\right) \vec n \\
          \psi
          \end{bmatrix}.
    \end{equation}
\end{thm}
\begin{pf}
Without loss of generality, we assume that node $\Mnode$ is assigned the index $N$ within its mesh element.
    In this context, the `right' node is node $\Pnode$ and the entropy
    produced at this node by the advective operator $\mathbf F^{a}$,
    disregarding the contribution from the advective dissipation operator $\text{DISS}^a$,
    corresponds to \eqref{p_a_cons_R}; direct substitution of \eqref{inviscid_boundary_conditions}
    into \eqref{p_a_cons_R} with the proper node index substitutions yields
    \begin{equation}
        \label{adv_bound_prod}
        \frac{dS^{a,\text{cons}}_{\Mnode}}{d t} = 
        -\left(\mathbf{w}_{\Mnode}^T \mathbf{f}^*(\mathbf u_{\Mnode},\mathbf u_{\Pnode}) - \Psi^*_{\Mnode} \right)= 0,
    \end{equation}
    which proves the theorem.  $\square$
\end{pf}

Now, we examine the effect of the dissipation operator $\text{DISS}^a$ applied to the solution jump at the boundary.
\review{In this work, we used} the local Lax-Friedrichs (LLF) dissipation, modified to act in entropy \review{variables},
\begin{equation}
    \label{LLF_definition}
    \text{DISS}^{a,\text{LLF}}_{(i,j)} = \frac{1}{2} |\lambda|_\text{max} \mathcal{H} [\![\mathbf w]\!]_{(i,j)},
\end{equation}
where $|\lambda|_\text{max}$ is the maximum absolute eigenvalue of the flux Jacobian considering both the $i$ and $j$ states.
%and the transformation matrix $\mathcal{H}$ is a numerical analog of the \review{entropy Jacobian}
%${\partial \mathbf u}/{\partial \mathbf w}$. %\ref{w_to_u_matrix_appendix}.
%\review{(underline in definition of $\mathcal{H}$ was removed)}.
%Our proof is quite general and allows for several possible definitions of $\mathcal{H}$.
%As one possibility, $\mathcal{H} = \mathcal{H}(\overline{\mathbf u})$ can be evaluated at
%some average state $\overline{\mathbf{u}} = \overline{\mathbf{u}}(\mathbf{u}_i,\mathbf{u}_j)$.
%Another option defines $\mathcal H = \mathcal H(\mathbf u_i, \mathbf u_j)$ to be
\review{The symmetric positive semi-definite matrix $\mathcal{H} = \mathcal{H}(\mathbf{u}_i, \mathbf{u}_j)$
is a numerical entropy Jacobian satisfying the relation
\begin{equation}
    \label{eq:llfmod_property}
    \left(\jump{\mathbf{u}}\right)_i = \left(\mathcal{H}\jump{\mathbf{w}}\right)_i, \quad i = 1, 2, 3, 4, 6, 7, 8, 9, \quad
    \left(\jump{\mathbf{u}}\right)_5 \approx \left(\mathcal{H}\jump{\mathbf{w}}\right)_5,
\end{equation}
following the construction recipes from \cite{derigs2017novel,IdealGLMMHD} to reconcile differences between LLF dissipation operators
formulated in conservative \review{vs.} entropy \review{variables}.
The expression for $\mathcal{H}$ is reproduced in Appendix~\ref{w_to_u_matrix_appendix}.
%In both cases, the entries of $\mathcal H$ result from some averaging procedure
%between states $i$ and $j$ denoted here with an overline.
The following theorem %substantiates the restriction on the average procedure,
guarantees entropy stability of the LLF dissipation operator at the boundary.}

\begin{thm}
    \label{LFF_boundary_theorem}
    The local Lax-Friedrich dissipation operator \eqref{LLF_definition}
    with the boundary state \eqref{inviscid_boundary_conditions} is 
    entropy stable at the boundary.% if in the definition of the matrix $\underline{\mathcal H}$, the
    %average state $\overline{\mathbf u}$ satisfies
    %\begin{equation}
    %    \label{average_restrictions}
    %    \overline{v}_n = 0, \quad \overline{B}_n = \Bext_n, \quad \overline{T} > 0, \quad \overline{\rho} > 0,
    %\end{equation}
    %where $\overline{v}_n$ and $\overline{B}_n$ are the normal components of the corresponding average vectors.
\end{thm}
\begin{pf}
Without loss of generality, we assume that the node $\Mnode$ has index $N$ inside the mesh element. 
In this context, the `right' node is denoted by $\Pnode$ and the
dissipative part of the advective entropy flux defined by \eqref{p_a_diss_R} 
relies on the computation of LLF operator \eqref{LLF_definition}. 
The jump in entropy variables at the boundary has only two non-zero components:
\begin{equation}
    \label{boundary_entropy_jump}
    [\![\mathbf w]\!]_{\Mnode,\Pnode} = 
    \left[0, -\frac{2 v_1}{T}, 0, 0, 0, \frac{2 \left(\Bext_1 - B_1\right)}{\mu_0 T}, 0, 0, 0\right]^T.
\end{equation}
Therefore, we represent \review{in \eqref{eq:diss_intermediate_calculations}} only the two columns of matrix $\mathcal{H}$,
interacting with the nonzero elements in \eqref{boundary_entropy_jump}.
\review{The ghost state \eqref{inviscid_boundary_conditions} results into the identities
$\avr{v_1} = 0$, $\avr{B_1} = \Bext_1$, and $\overline{p} = \rho R T$, simplifying expression \eqref{LLF_definition} at the boundary node}
%Under the assumption \eqref{average_restrictions}, the dissipative term equals
\begin{align}
    \label{eq:diss_intermediate_calculations}
    \text{DISS}^{a,\text{LLF}}_{\Mnode,\Pnode} =
    \frac{1}{2} |\lambda|_\text{max} {\mathcal{H}} [\![\mathbf w]\!]_{\Mnode,\Pnode} =
    \frac{1}{2} |\lambda|_\text{max}
    \begin{bmatrix}
        0 & 0 \\
        \review{\rho T} & 0 \\
        0 & 0 \\
        0 & 0 \\
        0 & \review{T} \Bext_1 \\
        0 & \review{T} \mu_0 \\
        0 & 0 \\
        0 & 0 \\
        0 & 0 \\
    \end{bmatrix}
    \begin{bmatrix}
        -2 v_1 / T \\
        2 (\Bext_1 - B_1) / \mu_0 T
    \end{bmatrix} \review{=
    |\lambda|_\text{max} \begin{bmatrix}
        0 \\
        - \rho v_1 \\
        0 \\
        0 \\
        \Bext_1 (\Bext_1 - B_1) / \mu_0 \\
        (\Bext_1 - B_1) \\
        0 \\
        0 \\
        0 \\
    \end{bmatrix}}.
\end{align}
At the boundary node, the contribution of the LLF dissipation operator to the time derivative of the global entropy function of the system
is expressed by using \eqref{LLF_definition} in
\eqref{p_a_diss_R} with corresponding substitution of indexes,
\begin{equation}
    \label{llf_entropy_contribution}
    \frac{dS^{a,\text{LLF}}_{\Mnode}}{dt} = \mathbf w^T_{\Mnode} \text{DISS}^{a,\text{LLF}}_{\Mnode,\Pnode} = -
    |\lambda|_\text{max} \review{\frac{1}{T}} \left( \rho v_1^2 + \frac{(\Bext_1 - B_1)^2}{\mu_0}\right).
\end{equation}
Positivity of density \review{$\rho$} and temperature \review{$T$}
guarantees nonpositivity of the expression \eqref{llf_entropy_contribution} and thus, the entropy stability of the LLF operator.  $\square$
\end{pf}

\subsection{Diffusive numerical flux divergence operator}
\label{discrete_diffuse_bc_section}

In the analysis of the visco-resistive term, we need to define not only the manufactured state
$\mathbf{v}^\nu_{\Pnode}$, but also the manufactured gradient of entropy variables $\mathbf g_{\Pnode}$ at the ghost node $\Pnode$,
both to be substituted in expressions \eqref{visc_discrete_def} and \eqref{g_discrete_def}.
Since our solid wall boundary conditions involve conditions on gradients of primitive variables,
we need a transformation procedure to go back and forth between gradients of primitive and entropy variables.
The procedure adopted in this work is as follow \cite{parsani2015entropy,dalcin2019conservative},
\begin{enumerate}
\item From the known gradient of entropy variables $\mathbf g_{\Mnode}$, we compute the gradient of primitive variables $\boldsymbol \theta_{\Mnode}$
  by multiplying $\mathbf g_{\Mnode}$ with the Jacobian matrix of the variable transformation
  $\mathbf w \rightarrow \mathbf v$ computed at the internal state $\mathbf v_{\Mnode}$,
  \begin{equation}
    \boldsymbol \theta_{\Mnode} = \left[\frac{\partial \mathbf v}{\partial \mathbf w}\right]_{\Mnode} \mathbf g_{\Mnode}.
  \end{equation}
  
\item From the gradient of primitive variables $\boldsymbol \theta_{\Mnode}$, we define the manufactured gradient
  of primitive variables $\boldsymbol \theta_{\Pnode}$ at the ghost node.
  Definitions for the components of $\boldsymbol \theta_{\Pnode}$ in terms of those of $\boldsymbol \theta_{\Mnode}$
  and local state $\mathbf v_{\Mnode}$ are specified in Theorems \ref{insul_theorem} and \ref{conduct_theorem} presented next.
  
\item From the gradient of primitive variables $\boldsymbol \theta_{\Pnode}$, we compute the gradient of entropy variables  $\mathbf g_{\Pnode}$
  by multiplying  $\boldsymbol \theta_{\Pnode}$ with the Jacobian matrix of the variable transformation
  $\mathbf v \rightarrow \mathbf w$ computed at the manufactured state vector $\mathbf v^\nu_{\Pnode}$,
  \begin{equation}
    \mathbf g_{\Pnode} = \left[\frac{\partial \mathbf w}{\partial \mathbf v}\right]_{\Pnode} \boldsymbol \theta_{\Pnode}.
  \end{equation}
\end{enumerate}

We first restrict our analysis to the case $\text{DISS}^\nu = 0$ and show the
entropy conservative properties of the manufactured boundary state.
Without loss of generality, we assume that the boundary node $\Mnode$ has index $N$
and lies within the rightmost mesh element, while the ghost node $\Pnode$ is to the right of node $\Mnode$.
Thus, with the particular form of boundary fluxes \eqref{BR1_eq}, the contribution of the 
diffusive dissipation term \eqref{p_nu_cons_R} to the 
time derivative of the discrete entropy function reads
\begin{equation}
    \label{diffusive_entropy_production}
    \frac{dS_{\Mnode}^{\nu,\text{cons}}}{dt} =
    \left( \hat{\mathbf w}_{\Mnode,\Pnode} - \mathbf w_{\Mnode} \right)^T \mathbf f^\nu_{\Mnode} + \mathbf w_{\Mnode}^T \hat{\mathbf f}^\nu_{\Mnode,\Pnode} =
    \frac{1}{2} \left( \mathbf w_{\Pnode}^T \mathbf f^\nu_{\Mnode} + \mathbf w_{\Mnode}^T \mathbf f^\nu_{\Pnode} \right).
    %- \mathbf v_{\Mnode}^T \mathbf \Mnode^\nu_{(\Mnode,\Pnode)} [\![\mathbf v]\!]_{(\Mnode,\Pnode)}.
  \end{equation}
Guided by Theorem \ref{viscous_entropy_theorem},
we define the manufactured state $\mathbf w_{\Pnode}^\nu$ and the gradient $\mathbf g_{\Pnode}^\nu$
such that the boundary diffusive flux, $\mathbf f^\nu_{\Pnode} = C^\nu_{\Pnode} \mathbf g_{\Pnode}^\nu$, satisfies the condition
\begin{equation}
  \label{diff_bound_prod}
  \frac{1}{2} \left( \mathbf w_{\Pnode}^T \mathbf f^\nu_{\Mnode} + \mathbf w_{\Mnode}^T \mathbf f^\nu_{\Pnode} \right) = g(t).
\end{equation}

We consider two cases.
The first one covers electrically insulating walls, cf.~\eqref{insulating_cond}. %which corresponds to Dirichlet boundary conditions for the magnetic field,
The second one covers the electrically conducting walls, cf.~\eqref{conductive_cond}. %which correspond to Dirichlet boundary conditions for the normal component of magnetic field and Robin and Neumann boundary conditions for the tangential component of magnetic field, 
Boundary conditions for the velocity and temperature follows \eqref{velocity_bc}, \eqref{thermal_bc}.
In the following theorems, we define manufactured states and gradients in the proper and general three-dimensional setting,
while the proofs are carried out under the one-directional assumption for simplicity.
In our notation, we omit the subscript $\Mnode$ on the individual components of the state $\mathbf v_{\Mnode} = (\rho, \vec v, T, \vec B, \psi)^T$,
and gradient $\boldsymbol \theta_{\Mnode}^T = ( \theta_{(\rho)}^T, \theta_{(\vec v)}^T, \theta_{(T)}^T, \theta_{(\vec B)}^T, \theta_{(\psi)}^T)$,
where $\theta_{(\cdot)}$ refers to the gradient components corresponding to a particular set of primitive variables.

The following theorem proposes the manufactured state and gradient for electrically insulating walls.
\begin{thm}
    \label{insul_theorem}
    For electrically insulating walls and
    disregarding the interface dissipation operator $\mathrm{DISS}^\nu$,
    the manufactured state and gradient
    \begin{equation}
        \label{Direchlet_bc}
        \mathbf v_{\Pnode}^{\nu,\mathrm{ins}} = \begin{bmatrix}
            \rho \\ - \vec v +  2 \vec v^\mathrm{\,w} \\ T \\ - \vec B + 2 \Bextvec \\ \psi
        \end{bmatrix}, \quad
        \boldsymbol \theta_{\Pnode}^\mathrm{ins} = \begin{bmatrix}
            -\theta_{(\rho)} \\ \theta_{(\vec v)} \\ -\theta_{(T)} + 2 \kappa^{-1}Tg(t) \vec n \\ \theta_{(\vec B)} \\ -\theta_{(\psi)}
        \end{bmatrix}
    \end{equation}
    guarantee that the penalty terms for the diffusive flux,
    together with the penalty on the gradient of the entropy variables, are
    \begin{itemize}
    \item entropy conservative if the wall is adiabatic, that is $g(t) = 0$,
    \item entropy stable in the presence of a heat entropy ﬂux, that is $g(t) \ne 0$, where $g(t)$ is a given $L^2$ function.
    \end{itemize}
\end{thm}
\begin{pf}
Direct substitution of \eqref{Direchlet_bc} into \eqref{diffusive_entropy_production} satisfies the condition \eqref{diff_bound_prod},
    If the heat entropy flux $g(t) = 0$, the contribution to the discrete entropy is zero, and the scheme is entropy conservative.
    If the heat entropy flux $g(t)$ is a given $L^2$ function, only a bounded amount of discrete entropy enters the domain
    at node $\Mnode$ during any time interval $[0,t]$. Thus, summing up all the contributions from all boundary nodes yields a bounded 
    contribution to the time derivative of the discrete entropy function in the system at any time $t$.  $\square$
\end{pf}

The following theorem proposes the manufactured state and gradient for the electrically conducting walls.
An additional motivation for the proposed manufactured state for the magnetic component given
in this theorem is presented in \ref{discrete_bc_motivation}. 
\begin{thm}
    \label{conduct_theorem}
    For electrically conducting walls and
    disregarding interface dissipation operator $\mathrm{DISS}^\nu$,
    the manufactured state and gradient
    \begin{equation}
        \label{Neumann_bc}
        \mathbf v_{\Pnode}^{\nu,\mathrm{cond}} = \begin{bmatrix}
            \rho \\ - \vec v + 2 \vec v^\mathrm{\,w} \\ T \\ \vec B \\ \psi
        \end{bmatrix}, \quad
        \boldsymbol \theta_{\Pnode}^\mathrm{cond} = \begin{bmatrix}
            -\theta_{(\rho)} \\ \theta_{(\vec v)} \\ -\theta_{(T)} + 2 \kappa^{-1}Tg(t)\vec n \\ 
            \theta_{(\vec B)}^{\ T} + 2 c_d^{-1} (\Bextvec - \vec B)^T \vec n \\ -\theta_{(\psi)}
        \end{bmatrix},
    \end{equation}
    guarantee that the penalty terms for the visco-resistive flux,
    together with the penalty on the gradient of the entropy variables, are
    \begin{itemize}
    \item entropy conservative if the wall is adiabatic, that is $g(t) = 0$,
    \item entropy stable in the presence of a heat entropy ﬂux, that is $g(t) \ne 0$, where $g(t)$ is a given $L^2$ function.
    \end{itemize}
\end{thm}
\begin{pf}
    Direct substitution of \eqref{Neumann_bc} into \eqref{diffusive_entropy_production} satisfies the condition \eqref{diff_bound_prod}.
    If the heat entropy flux is zero, \emph{i.e.}, $g(t) = 0$, the contribution to the time derivative of the discrete entropy function 
    is zero, and the scheme is entropy conservative.
    If the heat entropy flux, $g(t)$, is a nonzero function in $L^{2}$, only a bounded amount of discrete entropy enters the domain
    at node $\Mnode$ during any time interval $[0,t]$. Therefore, summing up all the contributions from all boundary nodes leads to a
    bounded contribution to the time derivative of the discrete entropy function in the system at any time $t$.  $\square$
\end{pf}

Finally, we examine the visco-resistive interface dissipation $\text{DISS}^\nu$ defined in \eqref{visc_diss_interface_term}.
The following theorem proves that this operator leads to discrete entropy dissipation. 
\begin{thm}
    \label{L_boundary_stability}
    The visco-resistive dissipation operator $\mathrm{DISS}^\nu$, defined in \eqref{visc_diss_interface_term},
    with either boundary state and gradient \eqref{Direchlet_bc} or \eqref{Neumann_bc},
    is entropy dissipative at boundary nodes.
\end{thm}
\begin{pf}
    Without loss of generality, we assume that the node $\Mnode$ has index $N$
    and lies within the rightmost mesh element, while the ghost node $\Pnode$ is to the right of node $\Mnode$.
    Under this assumption, and with the corresponding substitution of indices, the
    visco-resistive dissipation operator $\text{DISS}^\nu$ contributes to the discrete entropy 
    through the term \eqref{p_nu_diss_R}.

    For the insulating wall case,  direct substitution of the boundary state and gradient
    \eqref{Direchlet_bc} into \eqref{p_nu_diss_R} simplifies to
    \begin{equation}
        \label{Insul_diss_value}
        \begin{split}
            \frac{d}{dt}S_{\Mnode}^{\nu,\text{diss,ins}} &= \mathbf w_{\Mnode}^T \text{DISS}^{\nu,\text{ins}}_{\Mnode,\Pnode} = \\
        &-\frac{2 \beta}{T} \left[\mu_\text{NS}\left(\frac{4}{3} v_1^2 +
        \left(v_2^\text{w} - v_2\right)^2  + \left(v_3^\text{w} - v_3\right)^2 \right) +
        \frac{\mu_R}{\mu_0^2} \left(\left(\Bext_2 - B_2\right)^2 + \left(\Bext_3 - B_3\right)^2\right) \right],
        \end{split}
    \end{equation}
    which is non positive for a positive temperature, $T$, and positive values of $\beta$, $\mu_\text{NS}$, and $\mu_\text{R}$.

    For the conducting wall case,  direct substitution of boundary state and gradient
    \eqref{Neumann_bc} into \eqref{p_nu_diss_R} simplifies to 
    \begin{equation}
        \label{cond_diss_value}
        \frac{d}{dt}S_{\Mnode}^{\nu,\text{diss,cond}} = \mathbf w_{\Mnode}^T \text{DISS}^{\nu,\text{cond}}_{\Mnode,\Pnode} =
        - \frac{2 \beta \mu_{NS}}{T} \left(\frac{4}{3} v_1^2 + \left(v_2^\text{w} - v_2\right)^2  +
        \left(v_3^\text{w} - v_3\right)^2 \right),
    \end{equation}
    which is also non positive under the same conditions for $T$, $\beta$, and $\mu_\text{NS}$. $\square$
\end{pf}

\begin{remark}
  The visco-resistive flux does not depend on $\rho$, $\psi$, and their gradients. Therefore, it is not strictly necessary to
  specify them in the manufactured states and gradients \eqref{Direchlet_bc} or \eqref{Neumann_bc}.
  Nonetheless, the complete specification may be helpful in the case of using unconventional diffusive fluxes arising from alternative physical models.
  See \cite{svard2018new,sayyari2021entropy} for an example of the use of alternative viscous fluxes
  including mass diffusion under an unifying mechanism for the mass/momentum/energy diffusive processes.
\end{remark}

\section{Numerical results}
\label{results_section_label}
In this section, we preset numerical results to demonstrate the accuracy and
reliability of the developed boundary conditions for the DG-SBP-SAT discretization
operators.

The unstructured grid solver SSDC \cite{parsani2020high}
used herein has been developed at KAUST. The solver is built
on top of the Portable and Extensible Toolkit for Scientific Computing (PETSc) \cite{petsc-user-ref},
its mesh topology abstraction (DMPlex) \cite{knepley2009mesh} and scalable ordinary differential equation
(ODE)/differential algebraic equations (DAE) solver library \cite{abhyankar2018petsc},
and the Message Passing Interface (MPI).
The SSDC solver is also based on various algorithms for interior operators, 
\cite{carpenter2016towards,parsani2015entropyInterfaces,carpenter2014entropy,fernandez2020entropy,fernandez2020entropy-nasatm}
boundary nodes treatments that preserve the accuracy of the interior operators \cite{parsani2015entropy,dalcin2019conservative},
and time integration techniques
\cite{ranocha2019relaxation,ALJAHDALI2022111333,ranocha2021optimized}
It uses a transformation from computational to physical space that satisfies both the entropy conservation and the
geometric conservation law at the semi-discrete level \cite{fisher2012high,nolasco2020optimized}.
The computational meshes are generated with Gmsh \cite{geuzaine2009gmsh}.
Figures are created using Matplotlib \cite{Matplotlib} and ParaView \cite{ParaView}.

SSDC solves the governing equations \eqref{conservation_law_eq} in non-dimensional form.
In the following, dimensionless quantities are denoted with a prime,
and (dimensional) reference quantities are denoted with an asterisk.
The independent reference quantities are defined as follows,
\begin{equation}
    \label{eq:reference-scales}
    \begin{split}
        L^* \vec x' = \vec x, \quad
        U^* \vec v' = \vec v, \quad
        \rho^* \rho' = \rho, \quad
        T^* T' = T, \quad
        B^* \vec B' = \vec B, \\
        R^* R' = R, \quad
        \mu_{NS}^* \mu_{NS}' = \mu_{NS}, \quad
        \kappa^* \kappa' = \kappa, \quad
        \mu_{R}^* \mu_{R}' = \mu_{R}, \quad
        %\mu_0^* \mu_0' = \mu_0,
    \end{split}
\end{equation}
where $L^*$ is the reference length scale, $U^*$ is the reference velocity, $\rho^*$ if the reference density,
$T^*$ is the reference temperature, $B^*$ is the reference magnetic field magnitude,
$R^*$ is the reference gas constant, $\mu_{NS}^*$ is the reference dynamic viscosity,
$\kappa^*$ is the reference thermal conductivity, and $\mu_R^*$ is the reference electric resistivity of the system.
In all examples in this section, the properties of the gas are assumed constant and equal to the respective reference values
such that $R' = \mu'_{NS} = \kappa' = \mu'_R = 1$.
Substitution of \eqref{eq:reference-scales} into the governing equations \eqref{conservation_law_eq} results into the following
six non-dimensional numbers
\begin{align}
    \label{eq:non-dim-number-first}
    \gamma &= \frac{C_p}{C_v}, & \text{(Heat capacity ratio)} \\
    \Ma &= \frac{U^*}{\sqrt{\gamma R^* T^*}}, & \text{(Gasdynamic Mach number)} \\
    \Rh &= \frac{L^* U^* \rho^*}{\mu_{NS}^*}, & \text{(Hydrodynamic Reynolds number)} \\
    \Pr &= \frac{C_p \mu_{NS}^*}{\kappa^*}, & \text{(Prandtl number)} \\
    \Mm &= \frac{U^* \sqrt{\mu_0 \rho^*}}{B^*}, & \text{(Magnetic or Alfv\'en Mach number)} \\
    \label{eq:non-dim-number-last}
    \Rm &= \frac{\mu_0 U^* L^*}{\mu_R}, & \text{(Magnetic Reynolds number)}
\end{align}
where the heat capacities at constant pressure ($C_p$) and volume ($C_v$) are assumed to be constant.
Other non-dimensional quantities can be defined with the independent reference quantities from \eqref{eq:reference-scales}, including
\begin{equation}
    t' = \frac{L^*}{U^*}t, \quad
    E' = \frac{E}{\rho^* \gamma C_v T^*}, \quad
    F' = \frac{L^*}{\rho^* {U^*}^2} F, \quad
    \psi' = \frac{\psi}{B^*}, \quad
    c'_h = \frac{c_h}{U^*}, \quad
    \alpha' = \frac{L^*}{U^*}\alpha,
\end{equation}
where $F$ is an external forcing required later for the specification of one of our proposed test cases.

Unless otherwise stated, we define the damping parameter $\alpha'$ as a multiple of the hyperbolic divergence cleaning speed $c_h'$,
in particular we use $\alpha' = c_h' / 0.18$ following recommendations from \cite{IdealGLMMHD,DEDNER2002645}.
The value of $c_h'$ is chosen on a case-by-case basis so that it does not negatively impact time step size.
Unless otherwise stated, the \review{modified} local Lax--Friedrichs dissipation \review{\eqref{LLF_definition}} and
the visco-resistive dissipation from \review{\eqref{visc_diss_interface_term}} are applied to mesh element interfaces.

For integration in time, we use the fourth-order accurate Dormand--Prince method \cite{DORMAND198019}
adjusted with the adaptive time-stepping technique based on the digital signal processing
\cite{soderlind2003digital,SODERLIND2006225,ranocha2023error}.

%\textcolor{red}{
%In addition to the presented in this article solid wall boundary conditions,
%some test cases in this section requires inlet and outlet conditions.
%We specified them in an occasion way to make the problems solves and we do not pretend
%they are correct in all aspects.
%}

\subsection{Convergence study}
\label{pipeflow_testcase_section}

\begin{figure}
    \centering
    \begin{subfigure}[b]{0.4\linewidth}
        \centering
        \includegraphics[width=1\textwidth]{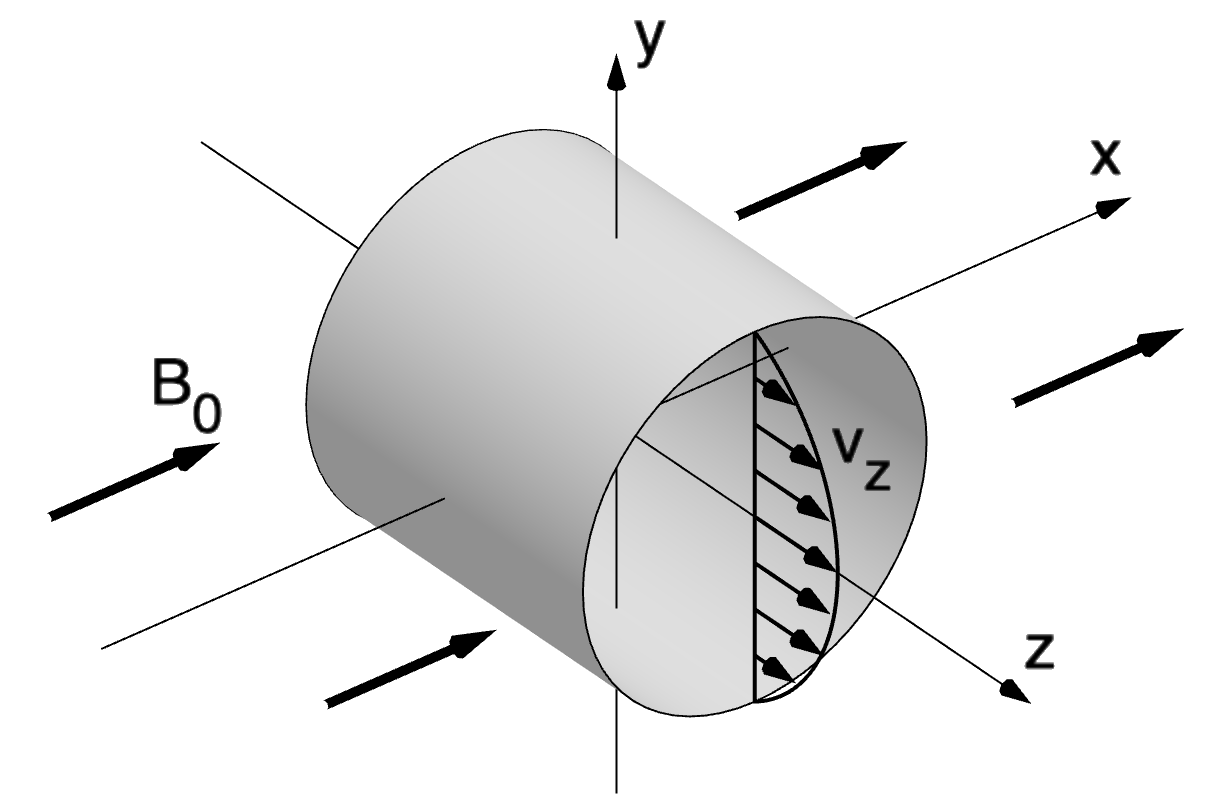}
    \end{subfigure}
 ~
    \begin{subfigure}[b]{0.27\linewidth}
        \centering
        \includegraphics[width=0.6\textwidth]{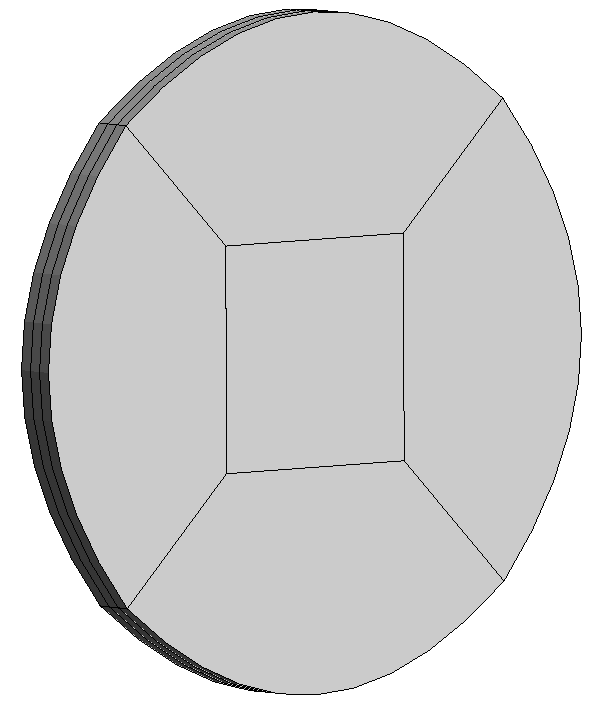}
    \end{subfigure}
 ~
    \begin{subfigure}[b]{0.27\linewidth}
        \centering
        \includegraphics[width=0.6\linewidth]{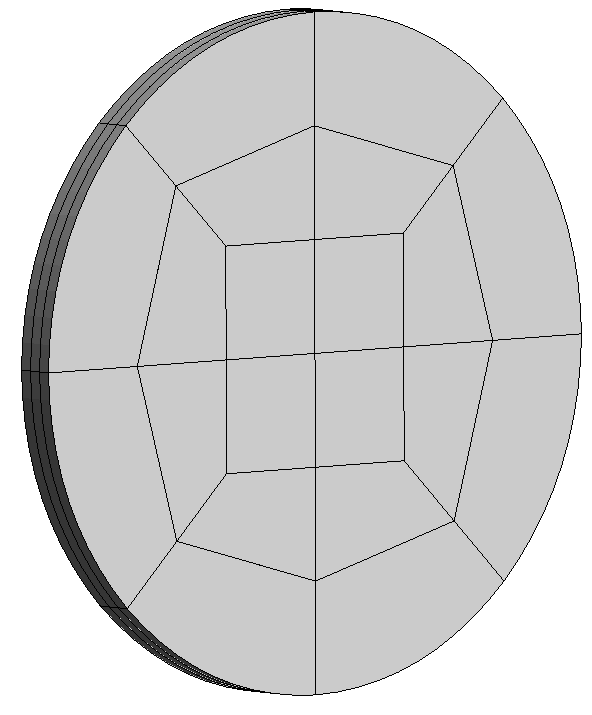}
    \end{subfigure}
 ~
    \caption{Convergence study: schematic, base hexahedral mesh, and the first level of mesh refinement for the pipe flow.}
    \label{fig:pipeflow-setup}
\end{figure}

In this section, we verify the accuracy of the proposed discrete boundary conditions. 
As a test problem, we consider the flow of an electrically conducting fluid in a circular pipe under an
external transverse uniform magnetic field and an external uniform force in the axial direction driving the fluid movement.
In the incompressible limit, this problem has a steady-state analytical solution expressed in terms of
modified Bessel functions of the first kind \cite{ihara1967flow}. The solution is
reproduced in \ref{solution_appendix_label} for completeness.

The flow domain is oriented such that the \( z \)-axis coincides with the longitudinal axis of a pipe of radius $a$ 
and the \( x \)-axis is aligned with an external magnetic field of strength \( B_0 \) (see Figure~\ref{fig:pipeflow-setup}).
The analytical solution takes the form
\[
 \vec{v} = [0, 0, u(x,y)], \quad \vec{B} = [B_0, 0, b(x,y)],
\]
and it is compatible with the boundary conditions defined in Section \ref{continuous_bc_section}.
Up to scaling factors depending on the external force and the Reynolds number,
the solution essentially depends on two non-dimensional parameters:
the wall conductance parameter $c$ defined in \eqref{wall_cond_param} with $L^*=a$ and 
the Hartmann number $\Ha$; see \ref{solution_appendix_label} for additional details.

\begin{figure}
    \centering
    \begin{subfigure}[b]{0.22\linewidth}
        \centering
        \includegraphics[width=\linewidth]{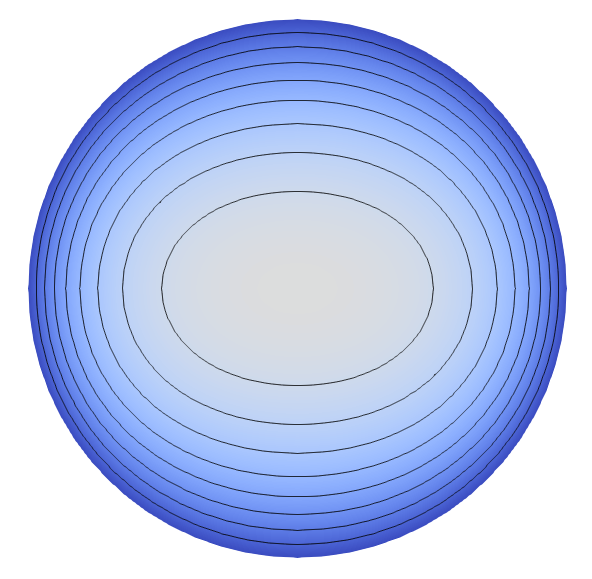}
        \caption{$u$, $c=0$.}
    \end{subfigure}
 ~
    \begin{subfigure}[b]{0.22\linewidth}
        \centering
        \includegraphics[width=\linewidth]{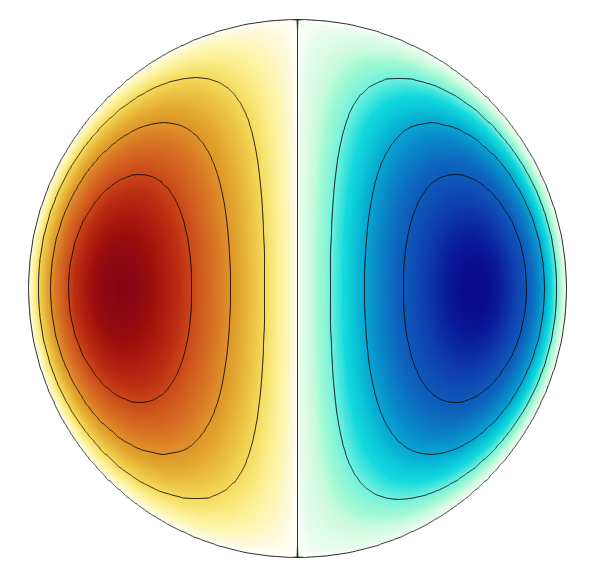}
        \caption{$b$, $c=0$.}
    \end{subfigure}
 ~
    \begin{subfigure}[b]{0.22\linewidth}
        \centering
        \includegraphics[width=\linewidth]{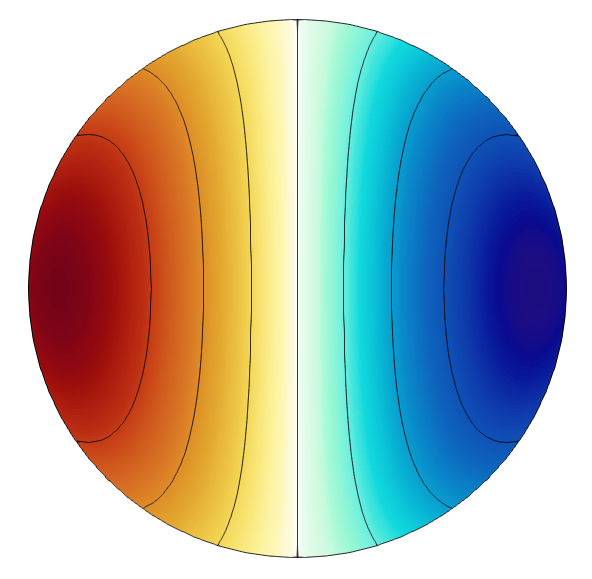}
        \caption{$b$, $c=1$.}
    \end{subfigure}
 ~
    \begin{subfigure}[b]{0.22\linewidth}
        \centering
        \includegraphics[width=\linewidth]{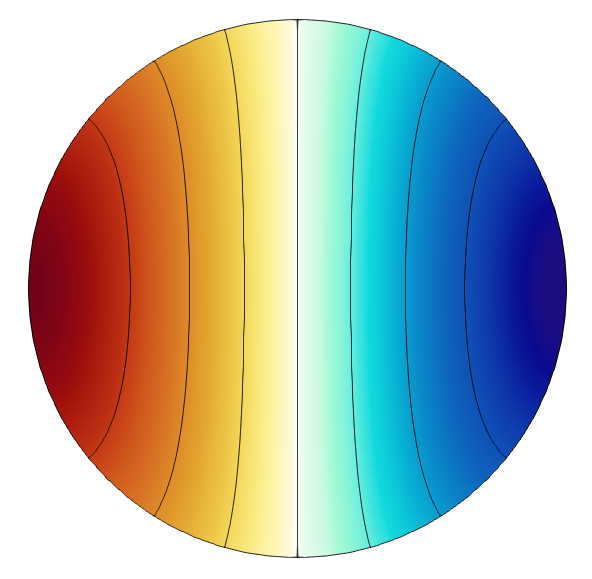}
        \caption{$b$, $c=\infty$.}
    \end{subfigure}
 ~
    \caption{Axial velocity, $u'$, and magnetic field, $b'$, in a pipe at three different values of
 the wall conductance parameter, $c=0,1,\infty$, at $Ha=5$. Due to the slight dependence of $u'$ on $c$ at $Ha=5$,
 only the plot for $u'$ at $c=0$ is shown.}
    \label{Ha5_pic}
\end{figure}

We study three cases for the wall conductance parameter:
\( c = 0 \), corresponding to an electrically insulating wall,
\( c = 1 \), corresponding to an electrically conducting wall, and
\( c = \infty \), corresponding to an electrically perfectly conducting wall.
The first case tests the electrical insulating properties of the discrete wall boundary conditions \eqref{Direchlet_bc},
while the other two test the electric conductivity properties of the discrete wall boundary conditions \eqref{Neumann_bc}.
In all cases, the Hartmann number is set to \( \Ha = 5 \). At this value of the Hartmann number,
electromagnetic forces are strong enough to deflect the axial velocity distribution from an axisymmetric  parabolic profile,
while the Hartman boundary layer is still thick enough to allow for the use of relatively coarse meshes as needed in a convergence study.

The resulting axial velocity and magnetic field distributions are plotted in Figure \ref{Ha5_pic}.
The conducting fluid moving in the magnetic field produces an electric current that flows along the iso-lines of the axial magnetic field, $b(x,y)$.
This electric current interacts with the magnetic field to produce a Lorentz force, which affects the velocity profile and exhibits anisotropy in the direction of $x$ vs. $y$.
Constant velocity contours expand from the center and compress towards the boundary in the \( x \)-direction along the external magnetic field.
In the case of an electrically insulating wall, the iso-lines of the magnetic field do not intersect the boundary,
which is consistent with no electric current flowing through the wall.
In the case of a perfectly conducting wall, the iso-lines of the magnetic field intersect the boundary at right angles, which is consistent with the electric
current flowing through the wall only in the wall-normal direction.

To obtain an accurate numerical solution for the incompressible equations using a compressible solver,
a small enough Mach number is required to achieve convergence.
An incompressible hydrodynamic Poiseuille flow solution and a uniform magnetic field are used as initial conditions.
The Mach number is progressively decreased in three stages: \( \Ma = 10^{-2}, 10^{-3}, 10^{-4} \).
At each stage, the solution is evolved to a steady state.
A Reynolds number $\Rh=1$ guarantees that this test will validate both viscous and advective parts of the boundary conditions. 
The hyperbolic divergence cleaning speed \( c'_h \) is of the same order as the maximum characteristic speed in the system but less than it 
to not affect the size of time step, $c'_h = \Ma^{-1}$.
We assume the wall to be adiabatic. Therefore, extra effort is needed to remove the heat generated by viscous and resistive mechanisms.
To this end, and to accelerate the decay of perturbations in the pressure field, we
applied to the energy equation damping mechanism similar to the one used for the generalized Lagrange multiplier $\psi$,
and then we redefine the damping term in \eqref{conservation_law_eq} to the following non-dimensional form
\begin{equation}
    \label{r_pipe_eq}
 \mathbf r' = [0, \vec 0, -\alpha'_1 \rho'_0 \gamma^{-1} (T' - T'_0), \vec 0, -\alpha'_2 \psi']^T,
\end{equation}
where $T'_0$ and $\rho'_0$ are undisturbed, reference non-dimensional temperature and density in the flow.
Following the recommendations in \cite{IdealGLMMHD,DEDNER2002645},
we set $\alpha'_1 = \text{Ma}^{-1}/0.18$, $\alpha'_2 = c'_h/0.18$ in all cases.
All the other parameters are defined from the Hartman number $\Ha$ and wall conductance parameter $c$ as detailed in \ref{solution_appendix_label}.

\begin{figure}
    \centering
    \includegraphics[width=\linewidth]{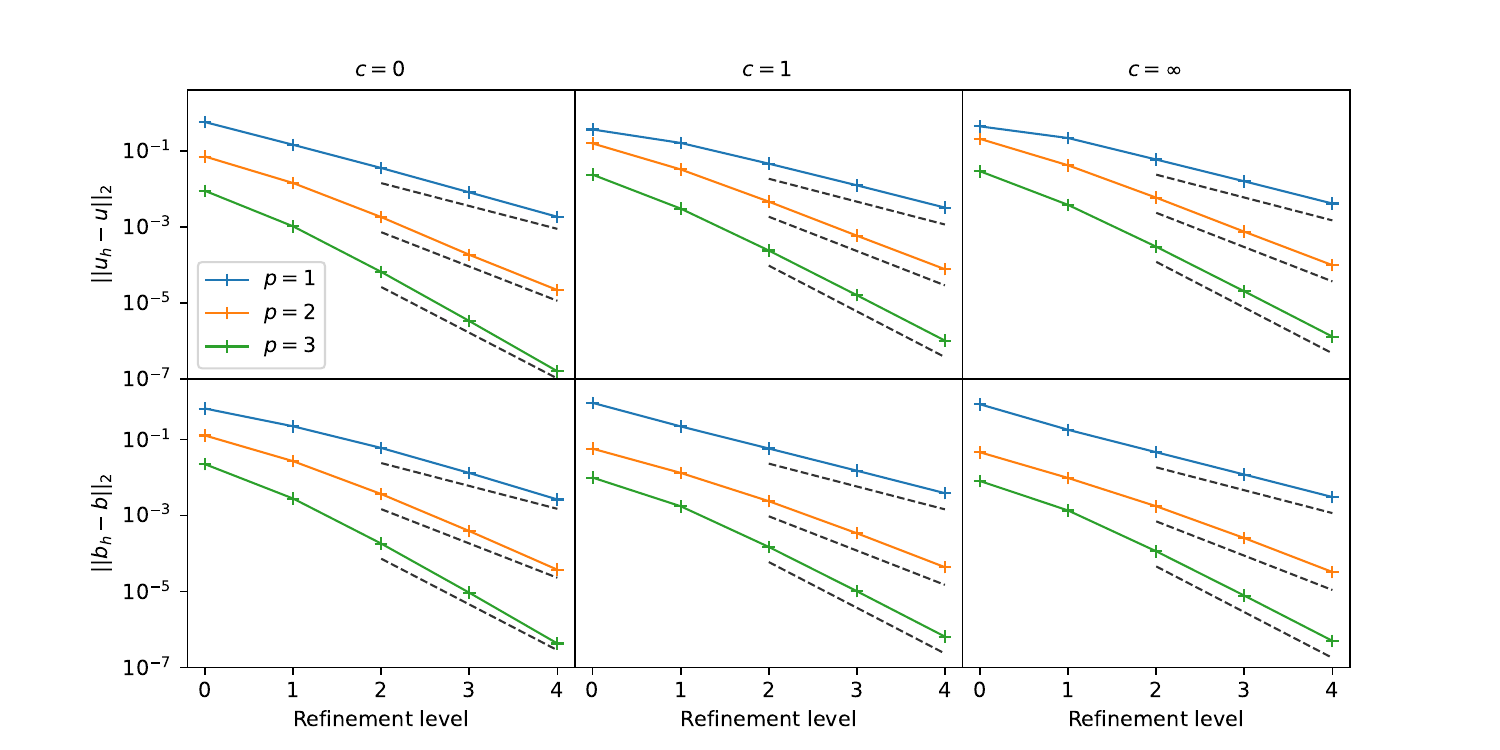}
    \caption{Convergence study for pipe flow with solution polynomial degrees \( p = 1, 2, 3 \);
 wall conductance parameter \( c = 0, 1, \infty \); \( \Ha = 5 \).
    \( L_2 \) errors in axial velocity \( u \) and axial magnetic field \( b \) compared
 to the analytical solution are shown with the ``expected" convergence rates
    \( o(h^{p+1}) \) (dashed lines). At a generic refinement level \( l \), the mesh has \( 5 \cdot 4^l \) elements
 in cross-section; in the axial direction number of elements is constant and equals $3$.
 The raw convergence data used for these plots is available in Table \ref{converg_table}.}
    \label{conv_plot}
\end{figure}

For the convergence study, we compute the discrete $L_2$ norm or the error between the analytical and numerical solutions
%For this convergence study, the error between the analytical and numerical solutions in the discrete $L_2$ norm is defined as
\begin{equation}
 |\!|e|\!|_{L_2} = |\!|\Omega|\!|^{-1} \left[\sum_{k=1}^{K} J_k \omega_k e_{k}^2 \right]^{1/2},
\end{equation}
where the index $k$ spans all nodes in the computational domain,
$||\Omega||$ indicates the volume of the computational domain,
$J_k$ is the determinant of the Jacobian of the transformation
from physical to computational space,
$\omega_k$ is the quadrature weight in computational space,
and $e_k$ is the difference between the analytical and numerical solutions at the node $k$.

The convergence study is performed using a sequence of nested three-dimensional butterfly-type hexahedral meshes.
The base mesh and the first level of refinement are shown in Figure~\ref{fig:pipeflow-setup}.
The length along the axial direction is equal to the $1/8$ of the pipe radius and is discretized with three mesh elements for all levels of refinement.
%The mesh resolution in the axial direction is fine enough not to affect the convergence rate at the parameters under the study.
The base mesh has five elements in the cross-section.
Each level of uniform mesh refinement splits each hexahedral cell into four new cells in the cross-section.
The meshes have curved boundary element faces to accurately approximate the cylindrical wall geometry using a polynomial degree that matches the solution polynomial order.
In all cases, we use absolute and relative tolerance of $10^{-11}$ for the adaptive time-stepping scheme to guarantee
the local truncation error in time is negligible compared to the space discretization error.

The results for the second, third, and fourth-order approximation of the numerical solution,
corresponding to $p = 1, 2, 3$ polynomial degrees and $N = 2, 3, 4$ LGL nodes in each coordinate direction,
are shown in Figure \ref{conv_plot}, where each column corresponds to different values of $c$.
For each value of $c$, the plot shows the $L_2$ error of the axial velocity and axial magnetic field -- these are the only two
solution components deviating from a uniform value in the analytical solution.
In all cases, the convergence rate reaches the ``expected" value of $p + 1$ (dashed lines).
The various data points from Figure \ref{conv_plot} are also available in the appendix in Table \ref{converg_table}.

\subsection{Entropy conservation study}
\begin{figure}[H]
\centering
\begin{subfigure}[b]{0.30\linewidth}
    \centering
    \includegraphics[width=1\linewidth]{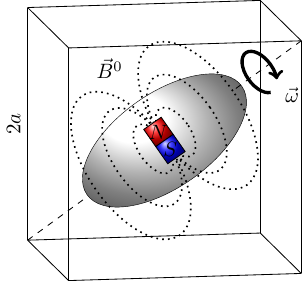}
\end{subfigure}
\qquad
\begin{subfigure}[b]{0.30\linewidth}
    \centering
    \includegraphics[width=1\linewidth]{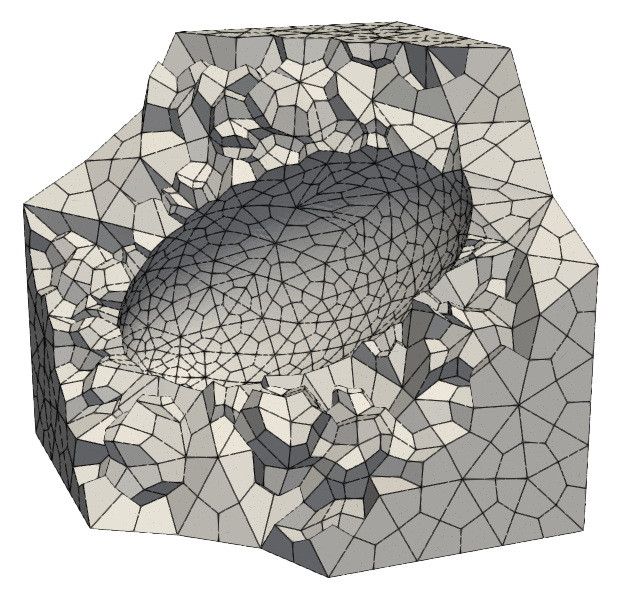}
\end{subfigure}
\caption{Discrete entropy conservation test: schematic of the domain and computational mesh for the spinning spheroid.}
\label{fig:spheroid-setup}
\end{figure}

This test case is designed to verify the entropy conservative properties of the spatial discretization algorithm
and our proposed hydrodynamic, thermal, and magnetic boundary conditions.
As shown in Figure~\ref{fig:spheroid-setup}, the computational domain has an exterior cubic shape.
Within the cube's interior, the plasma is set in motion by the rotation of an spheroidal solid body around its major axis,
which is aligned with one of the long diagonals of the cube.
An external magnetic field is generated by a magnetic dipole located within the spheroid.
The dipole axis matches the minor axis of the spheroid and rotates with it;
this way, the external magnetic field varies with time.

In Cartesian coordinates the plasma fills the interior of the cube of edge length $2a$,
$(x,y,z) \in [-a,a]^3$, and the exterior of the spheroid,
\begin{equation}
 \frac{\xi_1^2}{r_a^2} + \frac{\xi_2^2 + \xi_3^2}{r_b^2} \ge 1, \quad
    \xi_1 = \frac{x+y+z}{\sqrt{3}}, \quad \xi_2 = \frac{x-z}{\sqrt{2}}, \quad \xi_3 = \frac{x -2y + z}{\sqrt{6}},
\end{equation}
where the main axes of the spheroid are aligned with the orthonormal local coordinates $\xi_i$.
The polar and equatorial radiuses of the spheroid are $r_a = a$, $r_b = \frac{a}{2}$.
The reference length scale is equal to the half length of the cube edge, $L^*=a$.

The spheroid rotates around the semi-major axis $\xi_1$ with constant angular velocity $\vec \omega$,
\begin{equation}
  \vec \omega = \omega \frac{\partial \vec r}{\partial \xi_1}, 
%  \label{spheroid_surface_velocity}
\end{equation}
such that the spheroid surface velocity $\vec v^s$
at any point $\vec r$ is tangent to the surface of the spheroid, that is
\begin{equation}
  \vec v^s(\vec r) = \vec \omega  \times \vec r.
  \label{spheroid_surface_velocity}
\end{equation}
The reference velocity is equal to the maximum velocity attained at the equator of the spheroid, $\xi_1 = 0$, and equal to $U^* = \omega r_b$.

The external magnetic field $\Bextvec$ is generated by a magnetic dipole with magnetic moment $\vec m(t)$,
aligned with the semi-minor axis of the spheroid $\xi_2$ at $t=0$, and rotating together with the spheroid. More precisely,
\begin{equation}
    \label{spheroid_bext}
 \Bextvec(\vec r, t) = \frac{\mu_0}{4\pi} \left( \frac{3 \vec r (\vec r \cdot \vec m)}{r^5} - \frac{\vec m}{r^3}\right), \quad
 \vec m = \vec m(t) = m_0 \left(\frac{\partial \vec r}{\partial \xi_2} \cos{\omega t} + \frac{\partial \vec r}{\partial \xi_3} \sin{\omega t}\right),
\end{equation}
where $m_0$ is the strength of the magnetic dipole.
The reference magnetic field magnitude is defined as the maximum magnitude of the magnetic field attained on the surface of the spheroid
at $\vec \xi_\star = (0,r_b,0), t=0$, and equal to $B^* = |\!|\Bextvec(\vec \xi^\star, 0)|\!| = \mu_0 m_0 / 2\pi r_b^3$.

The reference density and temperature are set to the initial uniform density and temperature, $\rho^* = \rho$, $T^* = T$ at $t=0$.
The properties of the plasma are constant and defined through the non-dimensional numbers $\Ma = \Mm = \Rh = \Rm = \Pr = 1$, $\gamma = 1.1$.
The non-dimensional hyperbolic divergence cleaning speed is set to $c'_h = 1$.
The non-dimensional initial conditions at $t=0$ are set to $\rho' = T' = 1$, $\vec v' = 0$, $\vec B' = \Bextvec / B^*$, $\psi' = 0$.
Boundary conditions include zero velocity on the cube surface and non-uniform velocity $\vec v^s$ on the spheroid surface.
All walls are assumed adiabatic, that is, $g(t)=0$.

Two cases are considered to cover our proposed magnetic boundary conditions.
In the first case, the spheroid is electrically insulating, while the cube walls are electrically conducting with wall conductance parameter $c=1$.
In the second case, the spheroid surface is electrically conducting with $c=1$, while the cube walls are electrically insulating.
In both cases, the external magnetic field on all boundaries is $\Bextvec$ as defined in \eqref{spheroid_bext}.
Note that the two cases considered verify each of the discrete conditions \eqref{Direchlet_bc} and \eqref{Neumann_bc} in general circumstances, including:
curvilinear surface with non-zero velocity, non-homogeneous wall-normal magnetic field, not-steady state solution and boundary conditions,
the velocity exhibits a large discontinuity on the spheroid surface at $t=0$.
We do not explicitly test for a perfectly conducting wall with $c=\infty$.
This case corresponds to the homogeneous variant of \eqref{conductive_cond} where the right-hand side is zero
and is already covered by testing the magnetic boundary conditions with any positive value for wall conductance parameter, in particular, $c=1$.

The computational domain is discretized with a curvilinear computational mesh with 1,474 tetrahedral elements; each split into four hexahedrons (see Figure~\ref{fig:spheroid-setup}).
The mesh geometry and the solution are represented with polynomials of degree $p=3$.
Due to the nodal character of the entropy stability proofs of our boundary conditions,
the polynomial degree $p$ does not affect the entropy conservation properties. Therefore, the results for $p=3$ we present here extend to any other values of $p$.

At the discrete level, the spheroid surface velocity $\vec v^s$ is not exactly perpendicular to the (discrete) surface normal.
Therefore, in order to satisfy condition \eqref{velocity_bc}, values of $\vec v^s$ obtained from expression
\eqref{spheroid_surface_velocity} cannot be used directly to specify the wall velocity $\vec v^\text{w}$ in \eqref{Direchlet_bc} and \eqref{Neumann_bc}.
Instead, the spheroid surface velocity values from \eqref{spheroid_surface_velocity}
must be projected out of the discrete wall normal direction, that is
\begin{equation}
 \vec v^\text{w}_j = \vec v^s_j - (\vec v^s_j \cdot \vec n_j) \vec n_j,
\end{equation}
where the index $j$ goes through all the nodes on the spheroid's surface, and $\vec n_j$ is the unit normal at any such node.

Starting from the element-wise entropy conservation balance stated in \eqref{entropy_discrete_eq},
we turn off the interface dissipation mechanisms, $\text{DISS}^a = 0$ and $\text{DISS}^\nu = 0$,
as well as the damping term in the divergence cleaning mechanism, $\alpha = 0$.
Without these dissipation and damping terms, and after adding up the contributions of all the mesh elements,
the total rate of change of discrete entropy exactly balances the discrete entropy production due to dissipation, that is
\begin{equation}
 \label{eq:dSdt-DT-spheroid}
 \frac{dS}{dt} + DT = 0, \quad
 \frac{dS}{dt} = \sum_k J_k \omega_k \frac{d S_k}{d t} =
 \sum_k J_k \omega_k  \mathbf w_k^T \frac{d \mathbf u_k}{d t} , \quad
 DT = \sum_k J_k \omega_k \mathbf g_k^T \mathrm C^\nu_k \mathbf g_k,
\end{equation}
where index $k$ goes through all nodes in the mesh
and $\frac{d \mathbf u_k}{d t}$ corresponds to the time derivative in the left-hand side of \eqref{discrete_eq}.

\begin{figure}[H]
    \centering
    \begin{subfigure}[b]{0.45\linewidth}
        \centering
        \includegraphics[width=1\linewidth]{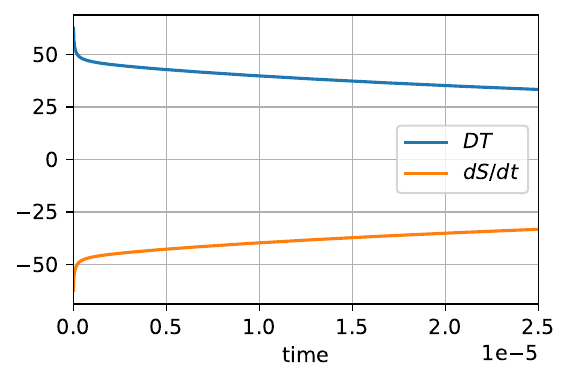}
        \caption{Insulating spheroid ($c=0$).}
        %\label{entropy_cons_c0}
    \end{subfigure}
    \begin{subfigure}[b]{0.45\linewidth}
        \centering
        \includegraphics[width=1\linewidth]{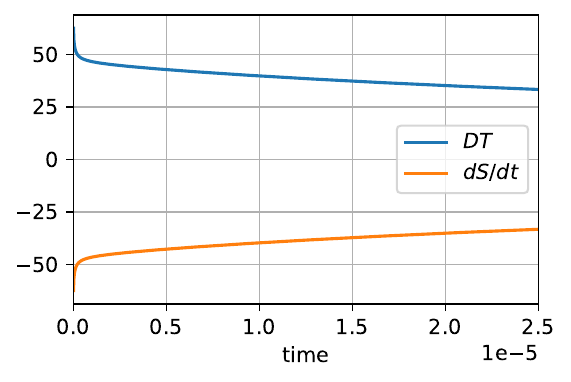}
        \caption{Conducting spheroid ($c=1$).}
        %\label{entropy_cons_c1}
    \end{subfigure}
    \begin{subfigure}[b]{0.45\linewidth}
        \centering
        \includegraphics[width=1\linewidth]{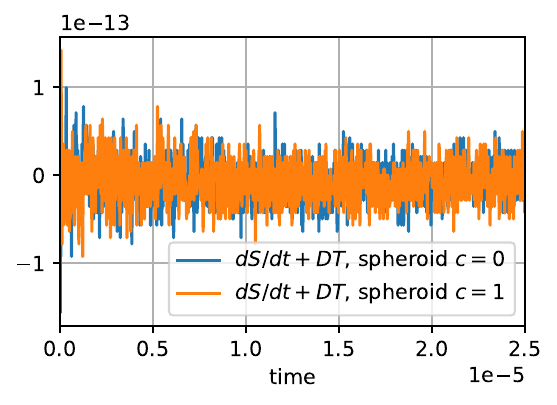}
        \caption{Instantaneous entropy balance.}
        %\label{entropy_cons_balance}
    \end{subfigure}
    \caption{Discrete entropy conservation test: time derivative of the entropy function, $dS/dt$, 
    dissipation term, $DT$, and instantaneous entropy production balance, $dS/dt+DT$, for electrically 
    insulating ($c=0$) and electrically conducting ($c=1$) surface of the spheroid.}
    \label{fig:dSdt-spheroid}
\end{figure}

The numerical verification of the entropy production balance \eqref{eq:dSdt-DT-spheroid} is shown in Figure \ref{fig:dSdt-spheroid}.
\review{Due to the small Reynolds numbers selected for this numerical conservation study, the time-step size is quite small ($\Delta t \approx 1.9\times 10^{-8}$). 
Thus, we show some results only for the first $1,300$ time-steps.}
The time rate of change of the total discrete entropy, $dS/dt$, is always negative,
while the entropy production visco-resistive dissipation $DT$ is always positive.
Furthermore, $dS/dt$ balances $DT$ up to values close to machine epsilon for IEEE double precision floating point arithmetic.

\enlargethispage{5mm}

\subsection{Fluid flow in microchannel}
\label{microchannel_section}

In this more practical example, we investigate a microfluidics application involving pumping and mixing in a microchannel.
MHD techniques offer significant advantages for microfluidic applications,
inducing bi-directional fluid motion without mechanical moving parts while
requiring low voltage drops, which is compatible with biological samples and related applications.
When an external magnetic field is applied to the fluid,
a pass-through electric current generates a Lorentz force that drives the fluid motion.

\begin{figure}[H]
\centering
\begin{subfigure}[b]{0.6\linewidth}
    \centering
    \includegraphics[width=\linewidth]{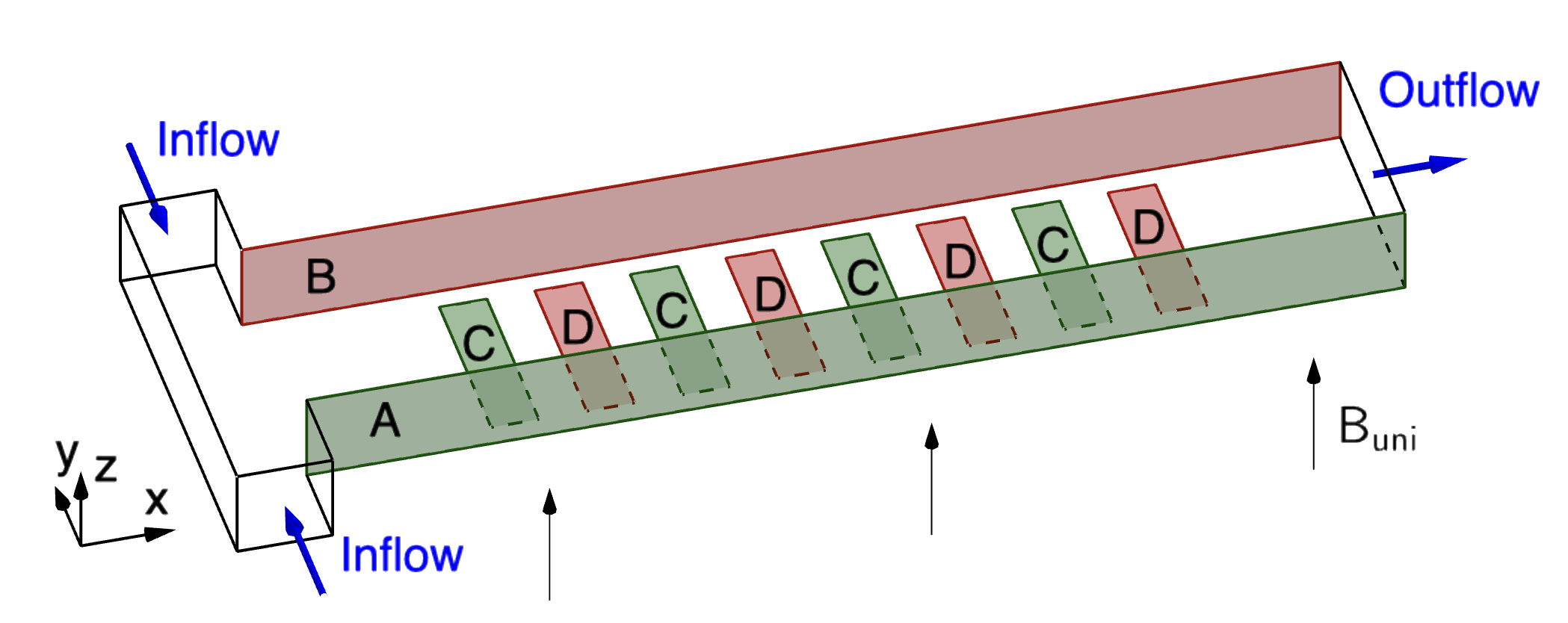}
    %\caption{}
    %\label{fig:channel_schematic}
\end{subfigure}
~
\begin{subfigure}[b]{0.3\linewidth}
    \centering
    \includegraphics[width=\linewidth]{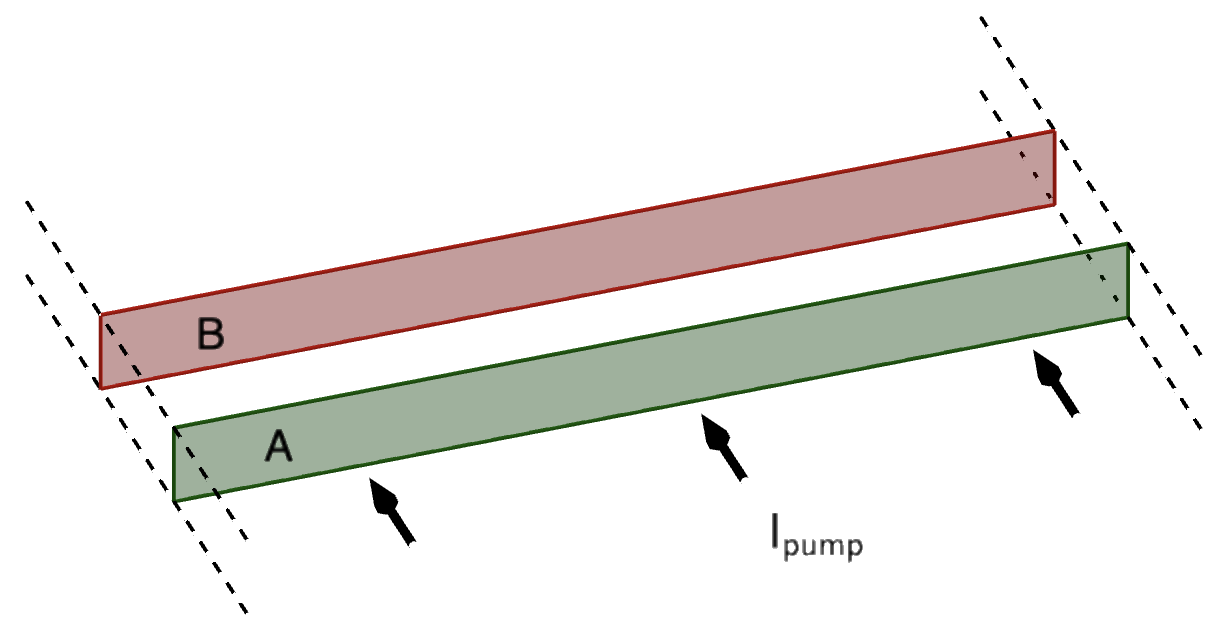}
    \includegraphics[width=0.9\linewidth]{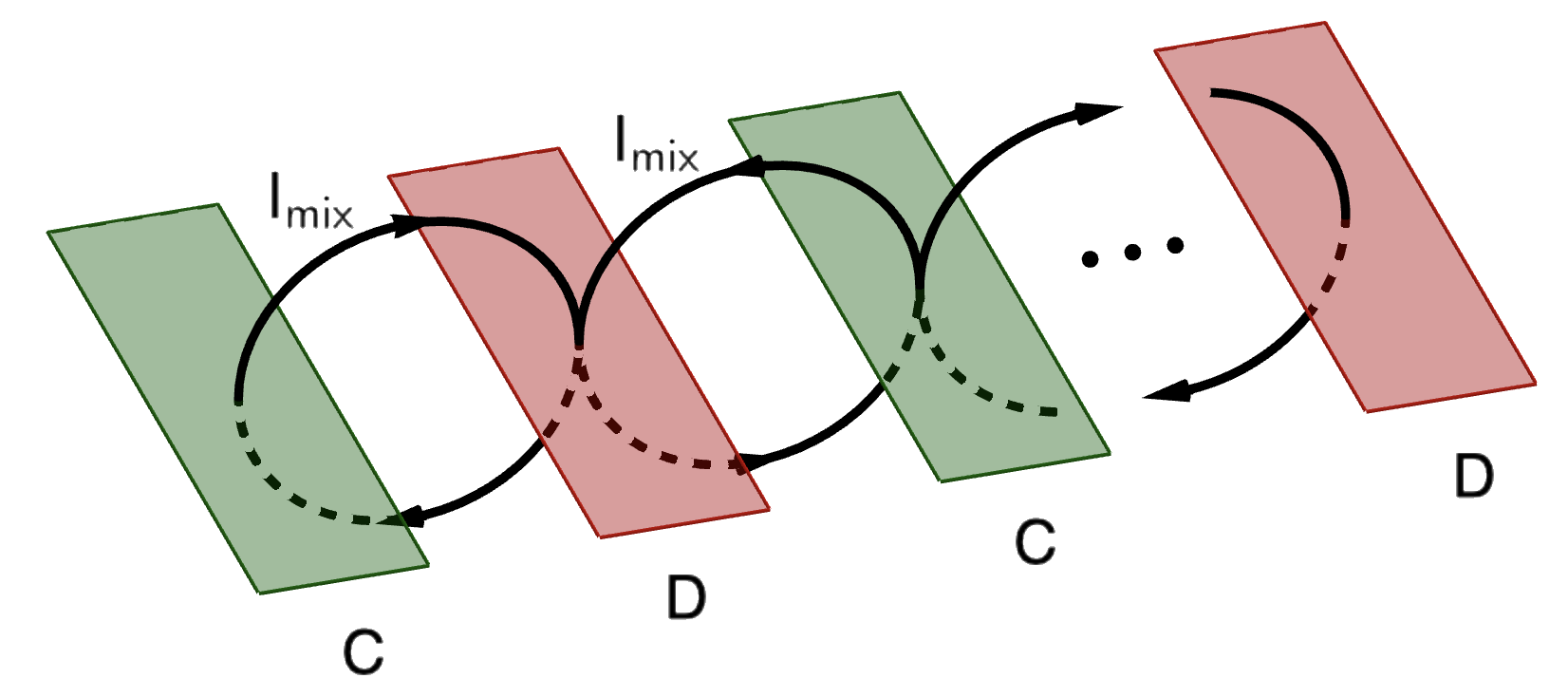}
    %\caption{}
    %\label{fig:channel_schematic}
\end{subfigure}
~
\caption{Microchannel flow test: schematic of the domain, uniform magnetic field $\vec B_\text{uni}$,
and currents $I_\text{pump}$ and $I_\text{mix}$.}
\label{fig:channel_schematic}
\end{figure}

A detailed description of the geometry used in this work can be found in \cite{lakim2014design} (Figure~3).
For completeness, we also report a schematic of the setup; see Figure~\ref{fig:channel_schematic}.
The main channel is a cuboid with a rectangular base and a rectangular cross-section.
The fluid comes in through symmetric side square channels at one end of the channel
and goes out through the other end of the main channel.
The channel features two sets of electrodes, denoted A-B and C-D.
The voltage drop between electrodes A and B, located at side walls $y=const.$, and the accompanying electric current
causes the fluid to flow along the channel (in $x$-direction) in an external magnetic field (in $z$-direction).
The voltage drop between another set of electrodes, \emph{i.e.}, electrodes C and D, located on the bottom of the channel at $z=0$,
causes the fluid to circulate in loops parallel to the $(x,y)$-plane.

The computational mesh includes 17,024 equal-sized cubic elements;
the numerical solution is approximated with polynomials of degree \( p = 5 \) within each element.
The solid wall boundaries of the channel are stationary ($\vec v^\text{\,w} = 0$) and adiabatic ($g(t) = 0$).
Initial conditions are set to $\mathbf v_\text{init} = [\rho^0,\vec 0, T^0, \Bextvec, 0]^T$,
where $\rho^0$, $T^0$ and $\Bextvec$ are the reference density, temperature, and external magnetic field, respectively.

At the inlet, we impose \review{a} uniform normal velocity \review{profile}, an external magnetic field $\Bextvec$ and zero normal gradients.
The inlet velocity magnitude is defined such that the inlet mass flow rate at time step $n+1$ equals
the outlet mass flow rate at the previous time step $n$. The initial mass flow rate is zero.
With the notation from Sections \ref{discrete_advect_bc_section} and \ref{discrete_diffuse_bc_section},
the inlet manufacture state and gradients are defined as
\begin{equation}
    \label{chanel_inlet_boundary}
 \mathbf v_{\Pnode}^\text{inlet} = \left[\rho^0 \left(\frac{p_{\Mnode}}{p^0}\right)^{\frac{1}{\gamma}}, -V \vec n,
 T^0 \left(\frac{p_{\Mnode}}{p^0}\right)^{\frac{\gamma-1}{\gamma}}, \Bextvec, 0\right]^T, \quad
    \theta_{\Pnode}^\text{inlet} = \mathbf \theta_{\Mnode} - (\mathbf \theta_{\Mnode} \cdot \vec n) \, \vec n,
\end{equation}
where $V = \rho_{\Pnode}^{-1} A_\text{inlet}^{-1} \dot m$ is the inlet velocity directed against the outward normal vector $\vec n$,
$A_\text{inlet}$ is the total inlet area, and
$\dot m = \int_\text{outlet} \rho (\vec v \cdot \vec n) da$ is the outlet mass flow rate at the previous time step.
At the ghost inlet boundary node, the density $\rho_{\Pnode}$ and temperature $T_{\Pnode}$ are defined such that
the pressure equals that of the internal state, that is $p_{\Pnode} = p_{\Mnode}$,
and the specific entropy equals that of the reference state, that is $s_{\Pnode} = s^0(\rho^0, T^0)$.
% = \frac{R}{\gamma-1}\ln{\frac{T^0}{T^*}} - R\ln{\frac{\rho^0}{\rho^*}}$.

At the outlet, we impose pressure, external magnetic field, and zero normal gradients,
\begin{equation}
 \mathbf v_{\Pnode}^\text{outlet} = \left[\rho_{\Mnode} \left(\frac{p^0}{p_{\Mnode}}\right)^{\frac{1}{\gamma}}, \vec v_{\Mnode},
 T_{\Mnode} \left(\frac{p^0}{p_{\Mnode}}\right)^{\frac{\gamma-1}{\gamma}}, \Bextvec, 0\right]^T, \quad
    \theta_{\Pnode}^\text{outlet} = \mathbf \theta_{\Mnode} - (\mathbf \theta_{\Mnode} \cdot \vec n) \, \vec n,
\end{equation}
where $p^0 = R \rho^0 T^0$ is the reference pressure.
At the ghost outlet boundary node, the density $\rho_{\Pnode}$ and temperature $T_{\Pnode}$ are defined such that
the pressure equals the reference pressure, that is $p_{\Pnode} = p^0$,
and the specific entropy equals that of the internal state, that is $s_{\Pnode} = s_{\Mnode}$.

We encode electric currents passing through the conducting electrodes in the external magnetic field.
The external magnetic field consists of three additive parts,
$\Bextvec = \vec B_\text{uni} + \vec B_{I_\text{pump}} + \vec B_{I_\text{mix}}$.
Here, \( \vec B_{\text{uni}} = (0, 0, B_\text{uni}) \) is a uniform background magnetic field directed along $z$-axis (see Figure~\ref{fig:channel_schematic}).
\( \vec B_{I_{\text{pump}}} \) represents the magnetic field from a uniform electric current $I_{\text{pump}}$
flowing in positive $y$-direction through an infinite straight conductor with rectangular cross-section, passing
electrodes A and B, as shown in Figure~\ref{fig:channel_schematic}.
The analytical expression for \( \vec B_{I_{\text{pump}}} \) can be found in \cite{NTMDT} and in \ref{rectangular_wire_appendix}.
\( \vec B_{I_{\text{mix}}} \) represents the magnetic field from a series of circular current loops $I_\text{mix}$,
crossing neighboring electrodes in the direction from C to D within the fluid, as shown in Figure~\ref{fig:channel_schematic}.
The analytical expression for \( \vec B_{I_{\text{mix}}} \) can be found in \cite{offaxisloop,jackson1999classical} and in \ref{current_loop_appendix}.
We assume patches A, B, C, and D to be perfectly electrically conducting, while other walls are perfectly electrically insulating.

Non-dimensional variables are defined through the thickness of the channel in the $z$-direction, $L^* = \qty{0.2}{\milli\meter}$ \cite{lakim2014design},
the reference density and temperature, $\rho^* = \rho^0, T^* = T^0$, and the magnitude of the uniform external magnetic field, $B^* = B_\text{uni}$.
The velocity scale $U^*$ is approximately equal to the bulk velocity.
The properties of the fluid are constant and defined such that $\Ma = 0.01$ and $\Rh = 1.5$,
which corresponds to a near-incompressible viscous flow.
Because of incompressibility, the velocity and magnetic field do not significantly depend on the heat capacity ratio
and thermal conductivity of the fluid, therefore we set $\gamma = 1.4$, $\Pr=0.72$.
The magnetic Mach and Reynolds numbers are set to $\Mm = 0.3$ and $\Rm = 1.5$.
The non-dimensional currents are $I'_{\text{pump}} = \mu_0 I_{\text{pump}} /B^*L^* = 0.72$, $I'_\text{mix} = \mu_0 I_\text{mix}/B^*L^* = 1$.
The hyperbolic divergence cleaning speed is $c'_h = \Ma^{-1}$.

\begin{figure}
    \centering
    \begin{subfigure}[b]{0.46\linewidth}
        \centering
        \includegraphics[width=\linewidth]{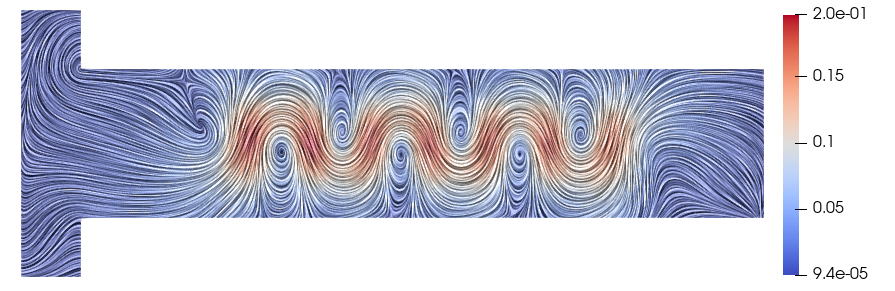}
        \caption{tangential magnetic field magnitude and lines}
        \label{channel_mag_lic}
    \end{subfigure}
 ~
    \begin{subfigure}[b]{0.46\linewidth}
        \centering
        \includegraphics[width=\linewidth]{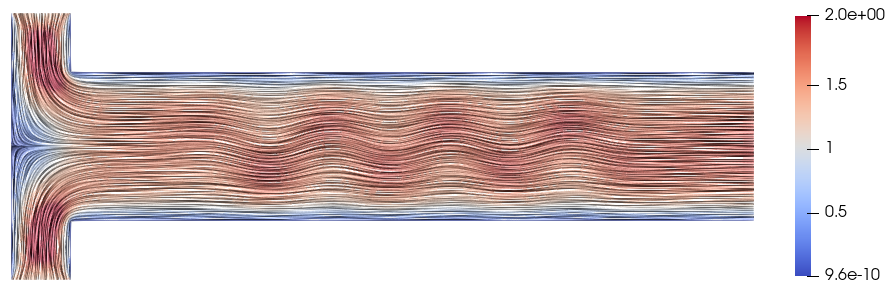}
        \caption{tangential velocity magnitude and lines}
        \label{channel_vel_lic}
    \end{subfigure}
    
    \begin{subfigure}[b]{0.46\linewidth}
        \raggedright
        \includegraphics[width=0.85\linewidth]{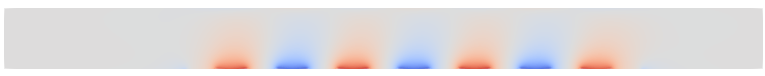}
        \caption{$B'_y$ on the mid $x$-$z$ plane}
        \label{channel_mag_y}
    \end{subfigure}
 ~
    \begin{subfigure}[b]{0.46\linewidth}
        \raggedright
        \includegraphics[width=0.85\linewidth]{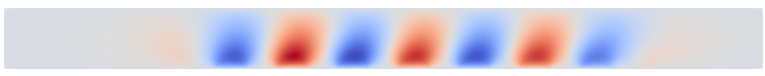}
        \caption{$v'_y$ on the mid $x$-$z$ plane}
        \label{channel_vel_y}
    \end{subfigure}

    \begin{subfigure}[b]{0.46\linewidth}
        \centering
        \includegraphics[width=\linewidth]{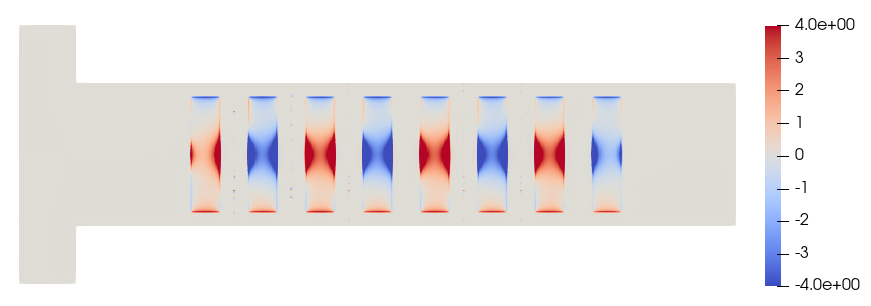}
        \caption{$(\vec \nabla' \times \vec B')_z$ on the bottom boundary, $z'=0$}
        \label{channel_current_z}
    \end{subfigure}
 ~
    \begin{subfigure}[b]{0.46\linewidth}
        \centering
        \includegraphics[width=\linewidth]{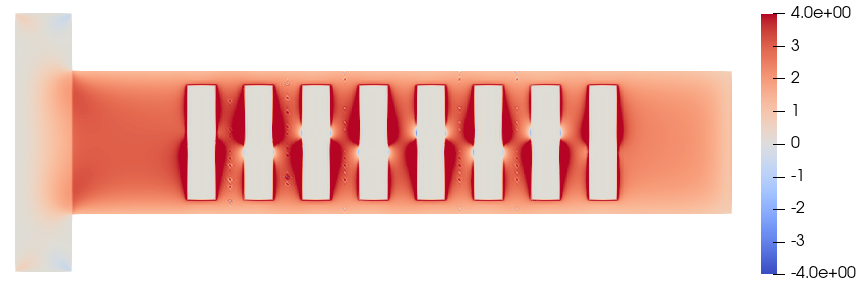}
        \caption{$(\vec \nabla' \times \vec B')_y$ on the bottom boundary, $z'=0$}
        \label{channel_current_y}
    \end{subfigure}
    \caption{Microchannel flow: steady state solution on mid $x$-$y$ plane (a, b), mid $x$-$z$ plane (c, d), and bottom surface (e, f).}
    %; (a), (b) -- visualization of the velocity and magnetic fields in the middle $x$-$y$ plane;
 %(c), (d) -- visualization of the velocity and magnetic field in middle $x$-$z$ plane;
 %(e), (f) -- visualization of $\vec \nabla \times \vec B$ proportional to electric current}
    \label{channel_pic}
\end{figure}

After an initial transient period, the system reaches a steady state shown in Figure~\ref{channel_pic}.
We chose the problem parameters to quantitatively reproduce a sinusoidal flow pattern, as described in \cite{lakim2014design}.
This sinusoidal flow pattern results from the combined effects of the pumping and mixing magnetic forces.
According to the design of the boundary conditions, normal to the wall electric current is close to zero
at the electrically insulating walls, as shown in Figure~\ref{channel_current_z} where wall normal component of $\vec \nabla \times \vec B$,
proportional to the electric current, is represented at the bottom boundary of the channel.
At the conducting patches, normal electric current diverges from zero and concentrates near the perimeter of the patches.
The total electric current through the conducting patches corresponds to the value imposed through the external magnetic field $\Bextvec$:
it is positive for patches C and negative for patches D. The tangential electric current is close to zero
at the perfectly conducting patches, as shown in Figure~\ref{channel_current_y} where one of the tangential components of
$\vec \nabla \times \vec B$ is represented at the same surface as in Figure~\ref{channel_current_z}.

\subsection{Self-field magnetoplasmadynamic thruster}
\label{thruster_section}

\begin{figure}
    \centering
    \begin{subfigure}[b]{0.4\linewidth}
        \centering
        \includegraphics[width=\linewidth,trim=3cm 3cm 3cm 2.8cm,clip]{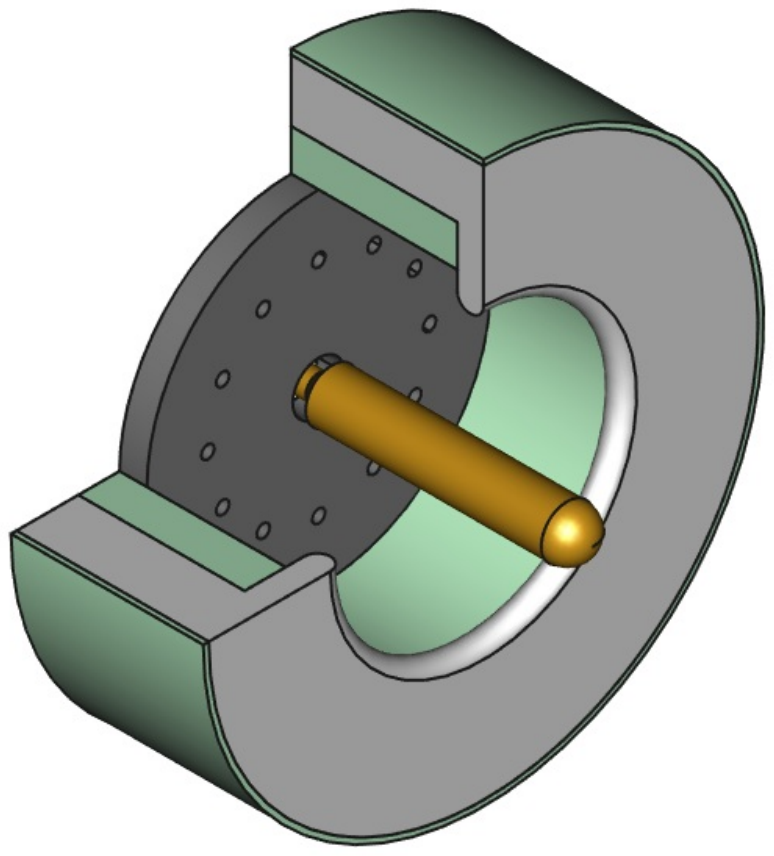}
        \includegraphics[width=\linewidth,trim=4cm 3cm 5cm 3cm,clip]{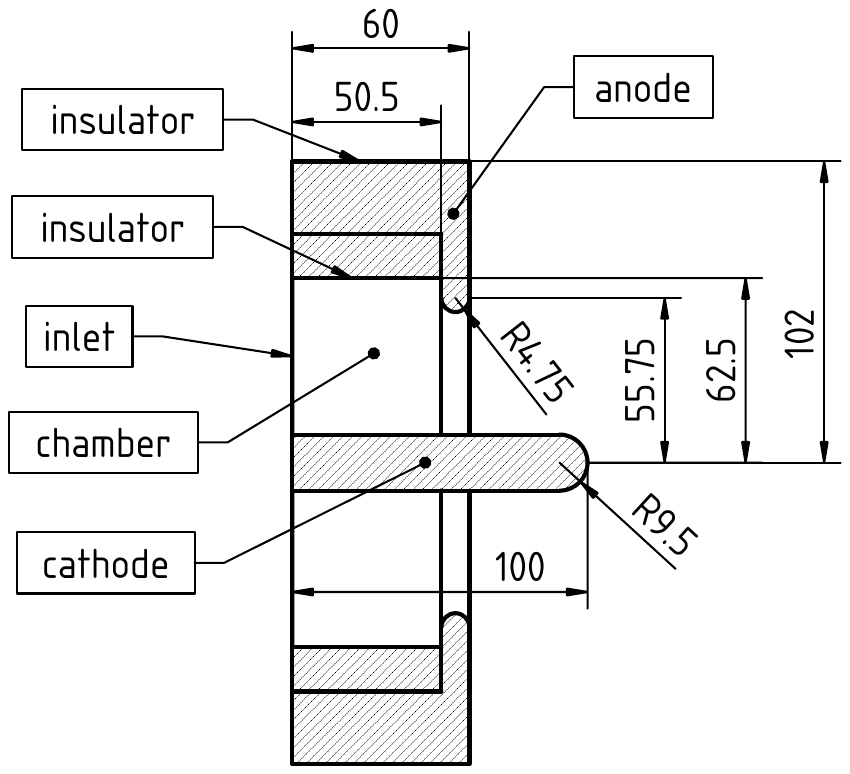}
        %\caption{}
        %\label{fig:thruster-schematic}
    \end{subfigure}
 ~
    \begin{subfigure}[b]{0.56\linewidth}
        \centering
        \includegraphics[width=\linewidth]{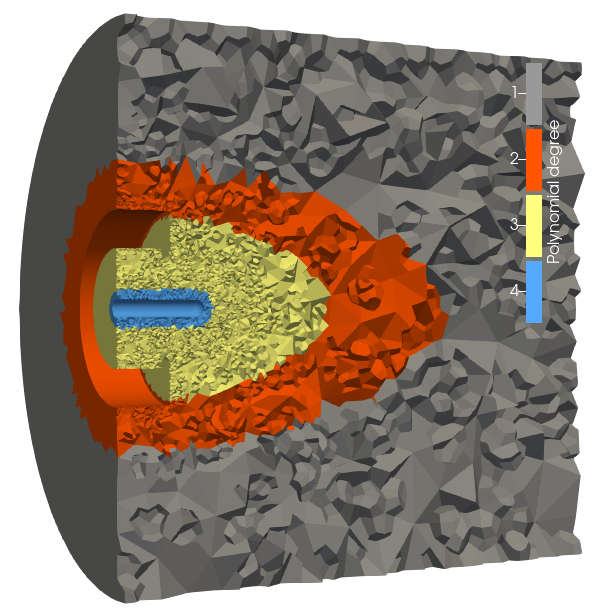}
        %\caption{}
        %\label{fig:thruster-grid}
    \end{subfigure}
    \caption{Princeton full-scale benchmark thruster: two-dimensional and three-dimensional schematics, computational domain, and grid. 
    Sizes are given in \unit{\milli\meter}. The solution is approximated with polynomials of varying degree from $p=\numrange{1}{4}$.}
    \label{fig:thruster-setup}
\end{figure}

In this final test case, we focus on the compressible flow of ionized gas in a self-field magnetoplasmadynamic (MPD) thruster, 
which is a type of electric propulsion device
that utilizes the principles of magnetohydrodynamics (MHD) to generate thrust by accelerating an ionized gas.
MPD thrusters are known for their compact design and efficiency, making them suitable for space applications.
For our simulation, we chose the Princeton full-scale benchmark thruster (PFSBT), which has been extensively studied 
experimentally \cite{burton1983measured} and numerically \cite{sankaran2005comparison,mayigue2018numerical}. 
The geometry and specifications of the PFSBT are well-documented in the literature, including
\cite{burton1983measured,sankaran2005comparison,mayigue2018numerical}. The thruster geometry features
a cylindrical cathode and an annular anode, as shown in Figure~\ref{fig:thruster-setup}. 
Argon is used as a propellant. The gas is injected in a chamber with insulating walls, partially through the 12 holes in the backplate
and partially through a narrow channel at the base of the cathode. An electric arc ionizes the gas within a few millimeters downstream of the inlet.
In its interaction with the self-induced magnetic field, the electric current flowing through the plasma accelerates the fluid by the Lorenz force.

Although our simulation does not aim to replicate the exact experimental regime, 
it provides valuable insights into the fundamental behaviors of the thruster operation under simplified assumptions. 
Specifically, we assume an ideal gas model, constant viscosity, thermal conductivity, and electrical resistivity for the plasma. 
The operating conditions for the MPD thruster include temperatures exceeding \qty{d4}{\kelvin}, where internal molecular degrees of freedom become significant.
In our model, we set the specific heat ratio $\gamma=\num{1.1}$, as follows from \cite{sankaran2005comparison} (Fig.~2).
The argon gas constant is $R_\text{Ar} = \qty{208.13}{\joule\per\kilo\gram\per\kelvin}$.

Boundary conditions include zero velocity on the surface of the thruster and zero heat flux through it.
The cathode and anode are perfectly electrically conducting, while the other parts of the thruster are perfectly electrically insulating.
At the internal electrically insulating surface bounding the thruster chamber and also at the inlet, the
external magnetic field is defined by the discharge current $I_d$ flowing through the cathode.
In cylindrical coordinates $(r,\theta,z)$, the external magnetic field is 
\begin{equation}
    \label{eq:thruster-B-ext-1}
 \Bext_r = \Bext_z = 0, \quad \Bext_\theta = \frac{I_d \mu_0}{4 \pi r}, \quad 
    \forall \vec x \in \Gamma^\text{inlet} \cup \Gamma^\text{int-ins}, \quad t > 0, 
\end{equation}
which is tangent to all the thruster surfaces. Thus, the normal magnetic field is zero at the conducting surfaces.
At perfectly conducting surfaces, tangent components of the external magnetic field do not influence the solution and, therefore, can be set to zero.
At the external insulating surface of the thruster, the external magnetic field can be assumed to be zero;
the electric current flowing through the bulk of the anode approximately shields the magnetic field induced by currents in the cathode.   
\begin{equation}
    \label{eq:thruster-B-ext-0}
 \Bextvec = 0, \quad \forall \vec x \in \Gamma^\text{cond} \cup \Gamma^\text{ext-ins}, \quad t > 0.
\end{equation}
In the simulation, the discharge current was set to $I_d = \qty{5}{\kilo\ampere}$ which is around \qty{30}{\percent}
of the operational discharge current equals \qty{16}{\kilo\ampere}.

Following \cite{sankaran2005comparison,mayigue2018numerical}, we further simplify the scheme of the propellent injection.
The propellent is injected through the whole backplate of the chamber at a temperature high enough for the gas to be already ionized.
The normal inlet velocity $V$ is defined from the constant mass flux at the inlet
resulting in a total mass flow rate of \review{$\dot m = \qty{2.4}{\gram\per\second}$}, \review{which is 40\% of} the operational mass flow rate of the thruster.
In the notation of Section \ref{discrete_advect_bc_section} and \ref{discrete_diffuse_bc_section}, the inlet boundary conditions are
\begin{equation}
    \label{thruster_inlet_bc}
 \mathbf v_{\Pnode}^\text{inlet} = \left[\rho^0 \left(\frac{p_{\Mnode}}{p^0}\right)^\frac{1}{\gamma}, -V \vec n, 
 T^0 \left(\frac{p_{\Mnode}}{p^0}\right)^\frac{\gamma-1}{\gamma}, \Bextvec, 0 \right]^T, \quad 
    \theta_{\Pnode}^\text{inlet} = \theta_{\Mnode} - (\theta_{\Mnode} \cdot \vec n) \, \vec n,
\end{equation}
where the velocity $V = \rho_{\Pnode}^{-1} A_\text{inlet}^{-1} \dot m$ is directed against the outward-pointing normal vector $\vec n$ and
$A_\text{inlet} \approx \qty{120}{\square\centi\meter}$ is the inlet area.
Density and temperature $\rho_{\Pnode}$ and $T_{\Pnode}$ are defined such that inlet pressure is that of the internal state, $p_{\Pnode} = p_{\Mnode}$,
and specific entropy is that of the reference state, that is, $s_{\Pnode} = s^0(\rho^0,T^0)$.
From a consideration of results from the previous simulations in \cite{mayigue2018numerical},
the reference density and temperature are set to $\rho^0 = \qty{0.1}{\gram\per\meter\cubed}$, $T^0 = \qty{d5}{\kelvin}$.
The reference pressure $p^0 = R \rho^0 T^0$.
The second expression in \eqref{thruster_inlet_bc} amounts to imposing zero normal gradients for all variables.
All other boundaries of the computational domain are treated as outlets.

The computational domain has a cylindrical shape with a diameter and length equal to six cathode lengths.
It is large enough for the far-field outlet boundaries to not noticeably to affect the solution, and
such that all significant physical phenomena are captured within the region of interest close to the thruster.
Outlet boundary conditions include zero magnetic field and external pressure \review{$p^\text{out} = 0.1p^0$,} approximating a vacuum condition.
In the notation of Section \ref{discrete_advect_bc_section} and \ref{discrete_diffuse_bc_section}, outlet boundary conditions are 
\review{
\begin{equation}
    \label{thruster_outlet_bc}
 \mathbf v_{\Pnode}^\text{outlet} = \left[\rho_{\Mnode}, \vec V^\text{out},
 T_{\Mnode} \left(\frac{p^\text{out}}{p_{\Mnode}}\right), \vec 0, 0 \right]^T, \quad
    \theta_{\Pnode}^\text{inlet} = \theta_{\Mnode} - (\theta_{\Mnode} \cdot \vec n) \, \vec n,
\end{equation}
}
where $\vec V^\text{out} = \vec v_{\Mnode} - 2 \min{(\vec n \cdot \vec v_{\Mnode}, 0)} \, \vec n$ is the velocity of the gas 
with a simple limiter to prevent backflow if $\vec n \cdot v_{\Mnode} < 0$.
\review{The density $\rho_{\Pnode}$ and the temperature $T_{\Pnode}$ are defined such that the pressure satisfies $p_{\Pnode} = p^\text{out}$.}

The problem is non-dimensionalized with the reference density and temperature, $\rho^* = \rho^0$, $T^* = T^0$, which are constant.
Reference length $L^* = \qty{0.1}{\meter}$ equals the length of the cathode and approximately equal to the anode annular diameter.
Reference velocity \review{$V^* = {\rho^*}^{-1} A_\text{inlet}^{-1} \dot m = \qty{2}{\kilo\meter\per\second}$} corresponds to the velocity of the gas at the inlet.
The reference magnetic field magnitude $B^* = \mu_0 I_d/2\pi r_c = \qty{0.1}{\tesla}$ is the maximum magnitude of the magnetic field
in the domain attained at the base of the cathode of radius $r_c = \qty{0.95}{\centi\meter}$.
\review{Based on these reference values, the Mach number is $\Ma = 0.4$ and the Alfv\'en Mach number is $\Mm=0.5$.}
\review{The} viscosity, electrical resistivity, and thermal conductivity are defined such that $\Rh=500$, $\Rm=50$, and $\Pr=2/3$.
The magnetic Reynolds number $\Rm=50$ is chosen to prevent convection of the
magnetic field lines outside of the thruster chamber to qualitatively reproduce \cite{sankaran2005comparison,mayigue2018numerical}.
\review{The} hyperbolic divergence cleaning speed \review{is set to} $c'_h = \Mm^{-1}$.

The unstructured computational mesh consists of 113,643 tetrahedral elements, each subdivided into four equally sized hexahedral elements, 
as shown in Figure~\ref{fig:thruster-setup}. The mesh is refined in the region of interest so that the maximum edge length of the tetrahedral elements 
near the cathode is $r_c/2$, while in the far-field zone, it grows progressively up to $5 r_c$.
An approximately steady-state solution was obtained using polynomials degrees gradually varying
from $p=1$ in the under-resolved far-field region and up to $p=4$ in the region close to the cathode boundary. 
The polynomial degree used in the different regions is indicated with varying colors in Figure~\ref{fig:thruster-setup}. 
An entropy-stable algorithm with h/p-nonconforming refinement was used, following the formulation from \cite{fernandez2020entropy}. 
The solution preserves axisymmetry, including $v_\theta = 0$ and $B_z = B_r = 0$.
The dimensionless divergence of the magnetic field in the $L_2$ norm remains under $|\!|\nabla' \cdot \vec B'|\!|_2 < \num{d-2}$.

\begin{figure}
    \centering
    \begin{subfigure}[b]{0.45\linewidth}
        \centering
        \includegraphics[width=\linewidth]{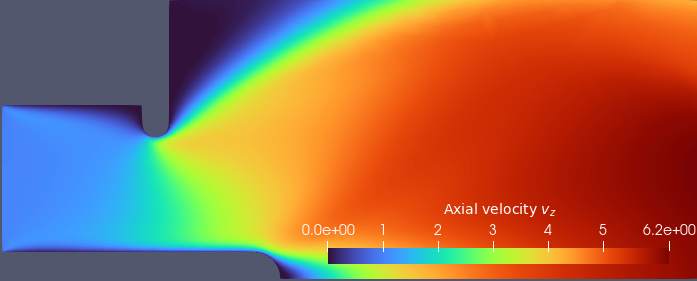}
        %\caption{}%Dimensionless axial velocity $v'_z$}
        %abel{fig:thruster-velocity}
    \end{subfigure}
 ~
    \begin{subfigure}[b]{0.45\linewidth}
        \centering
        \includegraphics[width=\linewidth]{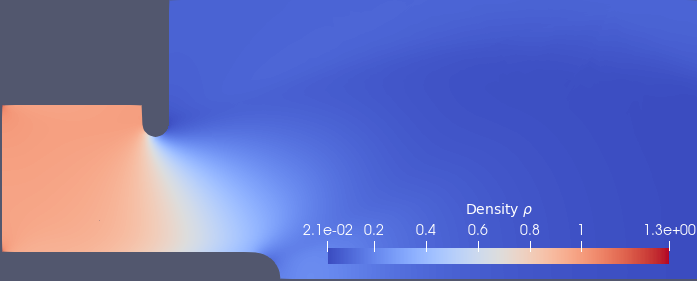}
        %\caption{}%Dimensionless density $\rho'$}
        %\label{fig:thruster-density}
    \end{subfigure}

    \begin{subfigure}[b]{0.45\linewidth}
        \centering
        \includegraphics[width=\linewidth]{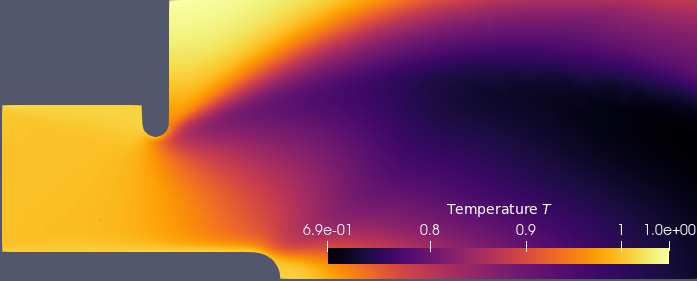}
        %\caption{}%Dimensionless temperature $T'$}
        %\label{fig:thruster-temperature}
    \end{subfigure}
 ~
    \begin{subfigure}[b]{0.45\linewidth}
        \centering
        \includegraphics[width=\linewidth]{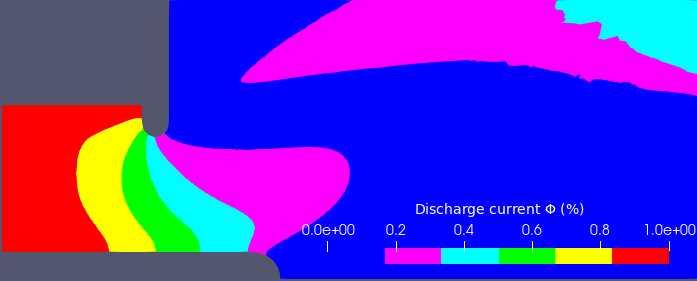}
        %\caption{}%Discharge current referred to the total current $I_d=\qty{5}{\kilo\ampere}$}
        %\label{fig:thruster-current}
    \end{subfigure}
    \caption{Princeton full-scale benchmark thruster: steady state solution for \review{$\Ma=0.4$, $\Mm=0.5$}.
    %discharge current $I_d=\qty{5}{\kilo\ampere}$.}
    }
    \label{fig:thruster-contours}
\end{figure}

\review{In agreement with experimental observations and previous numerical simulations from \cite{sankaran2005comparison},
the plasma accelerates as it passes through the thruster chamber, as shown in Figure~\ref{fig:thruster-contours}.}
Due to the axial symmetry of the solution, the electric current has a stream function $\Phi = r B_\theta$ shown in Figure~\ref{fig:thruster-contours}.
Electric current streamlines lie on iso-surfaces of $\Phi$. Consistently with the continuous analysis, \review{the electric current streamlines
which} originate at the cathode surface end at the anode surface. The gas tends to convect electric current streamlines downstream of the thruster chamber,
while the electrical resistivity keeps such spread within the nearby region.

\section{Conclusions}
\label{conclusions_section_label}
We have used entropy stability and the summation-by-parts framework to derive 
entropy conservative and entropy stable wall boundary conditions for the resistive magnetohydrodynamic 
equations in the presence of an adiabatic wall or a wall with a prescribed
 heat entropy flow, addressing three scenarios: electrically 
insulating walls, thin walls with finite conductivity, 
and perfectly conducting walls.
A point-wise entropy-stable numerical procedure has been presented for weakly 
enforcing these solid wall boundary conditions at the semi-discrete level, 
combining a penalty flux approach with a simultaneous approximation-term technique. 
The resulting semi-discrete operator mimics the behavior at the continuous level, 
and the proposed non-linear boundary treatment provides a mechanism for ensuring 
non-linear stability in the $L_2$ norm of the continuous and semi-discretized 
resistive magnetohydrodynamic equations.

The design order properties of the scheme are validated in the context of laminar 
flow of an electrically conducting fluid in a circular pipe, under an external transverse
uniform magnetic field and an external uniform force in the axial direction driving the
fluid movement.
The entropy conservation properties of the proposed numerical technique have been assessed
through a numerical experiment involving the flow around as stationary rotating spheroid.
The error in the entropy function balance showed an excellent agreement with
the theory for a mix of electrically insulating and electrically conducting walls.
The pumping and mixing of conducting fluid in a microchannel have been showcased as
an application to microfluidics.
The numerical results show very good agreement with the results available
from the literature.

The robustness of the complete semi-discrete operator (\emph{i.e.}, the entropy-stable 
interior operator coupled with the new boundary treatment) has been demonstrated 
for a self-field magnetoplasmadynamic thruster with a considerable pressure drop. 
The numerical simulation has been performed with a $p$-refinement approach, 
leading to a numerical solution that shows excellent agreement with the results 
available from the literature. This test has been successfully computed without 
the need to introduce artificial dissipation or filtering to stabilize the computations.
Although the robustness and efficacy of the techniques presented in this work
 have been validated using discontinuous spectral collocation operators on 
unstructured grids, the new boundary conditions can be applied to a comprehensive
 class of spatial discretizations. They are compatible with any diagonal-norm 
summation-by-parts spatial operator, including finite element, finite difference,
finite volume, nodal and modal discontinuous Galerkin, and flux reconstruction 
schemes.

\section*{Acknowledgements}
The research reported in this paper was funded by King Abdullah University of Science and Technology.
We are thankful for the computing resources of the Supercomputing Laboratory at King Abdullah University 
of Science and Technology and, in particular, Shaheen III supercomputer.

%{\color{red} Refer Shaheen III here with the link to the site because there are no article to cite yet.}

\appendix

\section{Entropy stable two-point flux function}
\label{iflux_appendix_label}

\label{two_point_flux_app}

The entropy consistent two point flux function for the discrete GLM-MHD system \eqref{discrete_eq} 
derived in \cite{IdealGLMMHD} is
\begin{equation}
\label{two_point_flux}
    \mathbf f^\text{EC}_{(i, j)} =
    \begin{pmatrix}
      f^\text{EC}_{1}\\
      \vdots\\
      f^\text{EC}_{5}\\
      \vdots\\
      f^\text{EC}_{9}
    \end{pmatrix}
    =
    \begin{pmatrix}
        \rho^\text{ln}\{\!\{v_1\}\!\} \\
        \rho^\text{ln}\{\!\{v_1\}\!\}^2 + \overline{p} + 
        \frac{1}{2\mu_0} \Bigl( \{\!\{B_1^2\}\!\} + \{\!\{B_2^2\}\!\} + \{\!\{B_3^2\}\!\} \Bigr) - \frac{1}{\mu_0} \{\!\{B_1\}\!\}^2 \\
        \rho^\text{ln}\{\!\{v_1\}\!\}\{\!\{v_2\}\!\} - \frac{1}{\mu_0} \{\!\{B_1\}\!\}\{\!\{B_2\}\!\} \\
        \rho^\text{ln}\{\!\{v_1\}\!\}\{\!\{v_3\}\!\} - \frac{1}{\mu_0} \{\!\{B_1\}\!\}\{\!\{B_3\}\!\} \\
        f^\text{EC}_{5} \\
        c_h \{\!\{\psi\}\!\} \\
        \{\!\{v_1\}\!\}\{\!\{B_2\}\!\} - \{\!\{v_2\}\!\}\{\!\{B_1\}\!\} \\
        \{\!\{v_1\}\!\}\{\!\{B_3\}\!\} - \{\!\{v_3\}\!\}\{\!\{B_1\}\!\} \\
        c_h \{\!\{B_1\}\!\}
    \end{pmatrix},
\end{equation}
where
\begin{align*}
  f^\text{EC}_{5} = f^\text{EC}_{1} & \Biggl[\frac{1}{(\gamma - 1)\beta^\text{ln}} - 
                   \frac{1}{2} \Bigl(\{\!\{v_1^2\}\!\} + \{\!\{v_2^2\}\!\} + \{\!\{v_3^2\}\!\}\Bigr) \Biggr] 
                  + f^\text{EC}_{2} \{\!\{v_1\}\!\} + f^\text{EC}_{3} \{\!\{v_2\}\!\} + f^\text{EC}_{4} \{\!\{v_3\}\!\} \\
                   + \frac{1}{\mu_0} & \Biggl[ f^\text{EC}_{6} \{\!\{B_1\}\!\} + f^\text{EC}_{7} \{\!\{B_2\}\!\} + f^\text{EC}_{8} \{\!\{B_3\}\!\}
                   - \frac{1}{2} \Bigl(\{\!\{v_1 B_1^2\}\!\} + \{\!\{v_1 B_2^2\}\!\} + \{\!\{v_1 B_3^2\}\!\}\Bigr) \\
                  & + \{\!\{v_1 B_1\}\!\}\{\!\{B_1\}\!\} + \{\!\{v_2 B_2\}\!\}\{\!\{B_1\}\!\} + 
                    \{\!\{v_3 B_3\}\!\}\{\!\{B_1\}\!\} +  f^\text{EC}_{9} \{\!\{\psi\}\!\} - c_h \{\!\{B_1 \psi\}\!\} \Biggr]
\end{align*}
and 
\begin{equation*}
    \overline{p} = \frac{\{\!\{\rho\}\!\}}{\{\!\{\beta\}\!\}}, \quad \beta = \frac{\rho}{p}.
\end{equation*}

The previous expressions use of the arithmetic and logarithmic means, respectively, defined as follow,
\begin{align*}
  \{\!\{a\}\!\} :=\frac{1}{2}(a_i + a_j), \quad 
  a^\text{ln} := \frac{a_j - a_i}{\ln{a_j} - \ln{a_i}}.
\end{align*}

%\section{Jacobian of the transformation from entropy to conservative variables}
\section{\review{Numerical entropy Jacobian used in the definition of LLF dissipation}}
\label{w_to_u_matrix_appendix}
\review{The numerical entropy Jacobian $\mathcal{H}$ designed for GLM-MHD system in \cite{IdealGLMMHD} completes the definition of LLF operator \eqref{LLF_definition}. In the following, we reproduce the expression for $\mathcal{H}$ reusing definitions from \ref{iflux_appendix_label}.}

Consider a state of the GLM-MHD system defined in primitive variables as $\mathbf v = (\rho, u, v, w, T, B_1, B_2, B_3, \psi)$. Then,
\review{
\begin{align*}
    \label{H_matrix}
    &\mathcal{H} = \mathcal{H}(\mathbf v_i, \mathbf v_j) = \\
    &\frac{1}{R}\begin{bmatrix}
\logavr{\rho} & \logavr{\rho} \avr{u} & \logavr{\rho} \avr{v} & \logavr{\rho} \avr{w} & \overline{E} & 0 & 0 & 0 & 0 \\
\logavr{\rho} \avr{u} & \logavr{\rho} \avr{u}^2 + \overline{p} & \logavr{\rho} \avr{u} \avr{v} & \logavr{\rho} \avr{u} \avr{w} & (\overline{E} + \overline{p}) \avr{u} &0&0&0&0\\
\logavr{\rho} \avr{v} & \logavr{\rho} \avr{v} \avr{u} & \logavr{\rho} \avr{v}^2 + \overline{p} & \logavr{\rho} \avr{v} \avr{w} & (\overline{E} + \overline{p}) \avr{v} &0&0&0&0\\
\logavr{\rho} \avr{w} & \logavr{\rho} \avr{w} \avr{u} & \logavr{\rho} \avr{w} \avr{v} & \logavr{\rho} \avr{w}^2 + \overline{p} & (\overline{E} + \overline{p}) \avr{w} &0&0&0&0\\
\overline{E}&(\overline{E}+\overline{p})\avr{u}&(\overline{E}+\overline{p})\avr{v}&(\overline{E}+\overline{p})\avr{w}&R\mathcal{H}_{5,5}&\tau\avr{B_1}&\tau\avr{B_2}&\tau\avr{B_3}&\tau\avr{\psi}\\
0      & 0            & 0            & 0            & \tau \avr{B_1}             & \tau \mu_0 & 0         & 0         & 0      \\
0      & 0            & 0            & 0            & \tau \avr{B_2}             & 0         & \tau \mu_0 & 0         & 0      \\
0      & 0            & 0            & 0            & \tau \avr{B_3}             & 0         & 0         & \tau \mu_0 & 0      \\
0      & 0            & 0            & 0            & \tau \avr{\psi}            & 0         & 0         & 0         & \tau \mu_0      \\
    \end{bmatrix}, \\
    &\overline{E} = \frac{1}{\gamma-1}\frac{\logavr{\rho}}{\logavr{\beta}} + \logavr{\rho} \left(\avr{u}^2+\avr{v}^2+\avr{w}^2\right) - \logavr{\rho} \frac{\avr{u^2+v^2+w^2}}{2}, \quad
    \tau = {\avr{\beta}}^{-1}, \\
    &R\mathcal{H}_{5,5} = \frac{\avr{u}^2+\avr{v}^2+\avr{w}^2}{2}\left(\overline{E}+\overline{p}\right) +
    \tau \left(\overline{E}+\frac{\avr{B_1}^2+\avr{B_2}^2+\avr{B_3}^2}{\mu_0}+\frac{\avr{\psi}^2}{\mu_0}\right).
    %&R\mathcal{H}_{5,5} = \frac{1}{\logavr{\rho}}\left(\frac{\left(\logavr{p}\right)^2}{\gamma-1}+\overline{E}^2\right) +
    %\overline{p}\left(\avr{u}^2+\avr{v}^2+\avr{w}^2\right) + \frac{\tau}{\mu_0} \left(\avr{B_1}^2 + \avr{B_2}^2 + \avr{B_3}^2 + \avr{\psi}^2\right).
\end{align*}
}

\begin{comment}
It is important for the definition of entropy stable matrix dissipation operator \eqref{Roe_upwind_eq}, that 
transformation matrix $\mathcal H$ is diagonal in the basis of right eigenvalues of the Flux Jacobian, 
defined in \eqref{right_ev_def}:
\begin{equation}
    \mathcal H = \mathcal R \mathcal Z \mathcal R^{-1}.
\end{equation}
Diagonal matrix $\mathcal Z$ is defined to be 
\begin{equation}
    \mathcal Z = \text{diag} \left( \frac{1}{2\gamma\rho}, \frac{T}{2\rho^2}, \frac{1}{2\gamma\rho}, \frac{T}{2}, 
    \frac{(\gamma-1)\rho}{\gamma}, \frac{T}{2}, \frac{1}{2\gamma\rho}, \frac{T}{2\rho^2}, \frac{1}{2\gamma\rho} \right).
\end{equation}
\end{comment}

\section{Motivation for the discrete boundary condition for electrically conducting walls}
\label{discrete_bc_motivation}
According to the expression for the diffusive flux \eqref{eq:diffusive-flux-mhd},
the magnetic boundary conditions for electrically conducting walls \eqref{conductive_cond}
define the normal resistive flux at the boundary, that is
\begin{equation}
    \label{resistive_flux_boundary}
    \overset{\leftrightarrow}{\mathbf f^\nu}_{(\vec B)} \cdot \vec n =
    \frac{\mu_R}{\mu_0} \left(\vec \nabla \vec B - \vec \nabla \vec B^T\right) \cdot \vec n \equiv
    \frac{\mu_R}{\mu_0} \left[ (\vec \nabla \times \vec B) \times \vec n \right]^T=
    \frac{\mu_R}{\mu_0} c_d^{-1} \left(\Bextvec_t - \vec B_t\right)^T,
    \quad \vec x \in \Gamma^\text{cond}, \quad t > 0,
\end{equation}
where $\Bextvec_t$ and $\vec B_t$ denote the tangential components of the corresponding vectors.
Note that the right hand side of \eqref{resistive_flux_boundary} is expressed here explicitly in terms of tangent components of the magnetic field;
this form is equivalent to the conditions in \eqref{conductive_cond} due to the continuity of the normal component of the magnetic field across the boundary.

The discrete boundary conditions \eqref{Neumann_bc} were designed to mimic \eqref{resistive_flux_boundary} on a discrete level.
The numerical resistive flux $\hat{\mathbf f}^\nu_{(\vec B)}(\mathbf u_{\Mnode}, \mathbf u_{\Pnode})$ defined in \eqref{BR1_eq}
is a linear function of the gradients, therefore it can be expressed as
\begin{equation}
    \label{eq:hat_theta}
 \hat{\mathbf f}^\nu_{(\vec B)}(\mathbf u_{\Mnode}, \mathbf u_{\Pnode}) =
 \frac{\mu_R}{\mu_0}\left(\hat{\theta}_{(\vec B)} - \hat{\theta}_{(\vec B)}^T\right),
\quad \hat{\theta}_{(\vec B)} = \frac{1}{2}\left(\theta_{(\vec B)}^{(-)} + \theta_{(\vec B)}^{(+)}\right),
\end{equation}
where $\hat{\theta}_{(\vec B)}$ is a numerical magnetic boundary gradient.
Direct substitution of the boundary gradients from \eqref{Neumann_bc} into the expression
for the numerical boundary gradient $\hat{\theta}_{(\vec B)}$ in \eqref{eq:hat_theta} results in
\begin{equation}
    \label{theta_b_boudnary}
 \hat{\theta}_{(\vec B)} =
 \frac{1}{2}\left(\theta_{(\vec B)} + \theta_{(\vec B)}^{T}\right)^{(-)} + c_d^{-1} \left(\Bextvec - \vec B\right)^T \vec n.
\end{equation}
As the resistive flux depends on the antisymmetric part of the gradients,
only the second term of \eqref{theta_b_boudnary} contributes to the normal numerical resistive flux at the boundary,
which therefore simplifies to
\begin{equation}
    \label{finial_expression}
 \hat{\mathbf f}^\nu_{(\vec B)}(\mathbf u_{\Mnode}, \mathbf u_{\Pnode}) \cdot \vec n=
 \frac{\mu_R}{\mu_0} c_d^{-1} \left( \left(\Bextvec - \vec B\right)^T \vec n - \vec n^T \left(\Bextvec - \vec B\right)\right) \cdot \vec n =
 \frac{\mu_R}{\mu_0} c_d^{-1} \left(\Bextvec_t - \vec B_t\right)^T.
\end{equation}
%where $\Bextvec_t$ and $\vec B_t$ denote the tangential components of the corresponding vectors.
%{\color{red}Note that boundary conditions for the normal component of the magnetic field $\vec B_n$
%are enforced through the advective boundary conditions \eqref{inviscid_boundary_conditions}.}
The discrete expression \eqref{finial_expression} mimics the continuous expression \eqref{resistive_flux_boundary}
and substantiates the particular form of the boundary gradient \eqref{Neumann_bc} for electrically conducting walls.
In particular, for the case of perfectly conducting boundary ($c_d^{-1}=0$), the numerical resistive flux through the boundary is identically zero.

\section{Analytical solution for MHD pipe flow}
\label{solution_appendix_label}
In this section, we present an analytical solution to the flow of a conducting fluid in a straight pipe of circular cross-section under uniform traverse magnetic field.
This analytical solution is derived in \cite{ihara1967flow} and used in section \ref{pipeflow_testcase_section} to validate
the accuracy of the solid wall boundary conditions we propose in this work.
\review{A} uniform external force in the axial direction causes the fluid motion, and the flow is assumed incompressible.  

Within the pipe, the solution satisfies the incompressible visco-resistive MHD equations,
\begin{equation}
    \begin{split}
    & \rho (\vec v \cdot \vec \nabla) \vec v = - \nabla \vec p + \eta {\vec \nabla}^2 \vec v + 
    \frac{1}{\mu_0} \left(\vec \nabla \times \vec B \right) \times \vec B + F, \quad \vec \nabla \cdot \vec v = 0,\\
    & \sigma \mu_0 \vec \nabla \times \left( \vec v \times \vec B \right) - 
    \vec \nabla \times \left(\vec \nabla \times \vec B\right) = 0, \quad \vec \nabla \cdot \vec B = 0,
    \end{split}
\end{equation}
where the density $\rho$, viscosity $\eta$, and conductivity $\sigma$ of the fluid are constant, and
$F$ is an external force.

We set a Cartesian coordinate system such that  the
$x$-axis is parallel to the uniform traverse external magnetic field of strength $\Bext$, and the
$z$-axis coincides with the longitudinal axis of the pipe of radius $a$.
Assuming laminar flow, the velocity and magnetic field of the solution take the form
\begin{equation}
    \vec v = (0, 0, u), \quad \vec B = (\Bext, 0, b),
\end{equation}
where functions $u=u(r,\theta)$ and $b=b(r,\theta)$ are defined in the cross-section in non-dimensional polar coordinates 
$r = a^{-1} \sqrt{x^2 + y^2}$ and $\theta$, and satisfies the following boundary conditions,
\begin{equation}
  u = 0, \quad b + c \frac{\partial b}{\partial r} = 0, \quad \text{at }  r=1,
  \label{bcpipeflow}
\end{equation}
where $c=\tau_w \sigma_0 / a \sigma$ is a dimensionless wall conductance parameter
depending on the thickness of the wall $\tau_w$ and conductivity of the wall material $\sigma_0$. 
% Any value from $c=0$ for electrically insulating walls to $c=\infty$ for perfectly conducting walls is allowed.
For electrically insulating walls ($c=0$), the boundary condition for $b$ \eqref{bcpipeflow} reduces to the Dirichlet condition $b|_{r=1} = 0$.
For perfectly conducting wall ($c=\infty$), \eqref{bcpipeflow} reduces to the Neumann condition $\partial_r b |_{r=1} = 0$.

%Two other dimensionless numbers governing the solution are the Hartmann number $\text{Ha}=a \Bext \sqrt{\sigma/\eta}$,
%and the Poiseuille number $P = -f a^2 / \eta v_0$, where $f$ is the external pressure gradient and $v_0$ is a reference velocity.
The second governing dimensionless parameter is the Hartmann number $\text{Ha}=a \Bext \sqrt{\sigma/\eta}$,
expressing the ratio of electromagnetic force to the viscous force. 
With the auxiliary constant $k=\text{Ha}/2$, the solution is given as a truncated series of the form
\begin{equation}
    \begin{split}
    \label{pipeflow_solution_eq}
    \frac{\eta}{a^2 F} u(r,\theta) = \tilde u(r,\theta) = - &\sum_{n=0}^{N_0} \frac{\epsilon_n}{k} M_n \frac{(-1)^n e^{kr\cos{\theta}} + e^{-kr\cos{\theta}}}{2} 
    I_n(kr) \cos(n\theta), \\
    \frac{\sqrt{\eta}}{a^2 F \mu_0 \sqrt{\sigma}} b(r,\theta) = \tilde b(r,\theta) =  + &\sum_{n=0}^{N_0} \frac{\epsilon_n}{k} M_n 
    \frac{(-1)^n e^{kr\cos{\theta}} - e^{-kr\cos{\theta}}}{2} I_n(kr) \cos(n\theta)-\frac{r}{2k} \cos(\theta),
    \end{split}
\end{equation}
where $I_n(k)$ are modified Bessel functions of the first kind. The coefficients $\epsilon_n$ and $M_n$ are defined as 
\begin{equation}
    \label{pipeflow_coeffitients_eq}
    \begin{split}
    \epsilon_n =& \begin{cases}
        {1}/{2},& \text{if } n = 0, \\
        1,      & \text{otherwise}, 
    \end{cases} \\
    M_n =& \sum_{s=0}^{S_0} A_s \frac{I_{n+2s+1}(k) + I_{n-2s-1}(k)}{I_n(k)}.
    \end{split}
\end{equation}
Bessel functions of negative order can be computed based on their symmetry, $I_n = I_{-n}$.
The coefficients $A_s$ are computed from the solution of the linear system of equations $B_{ms}A_s=C_m, 0\le s,m \le S_1$, 
where $S_1 \le S_0$. 
For finite values of $c$, we have
\begin{equation}
    \label{pipeflow_matrix_coeffitients_eq}
    \begin{split}
    B_{ms} &= \delta_{ms} + c k \sum_{n=0}^{N_1} (-1)^{n+1} \frac{\epsilon_n I_n^\prime(k)}{I_n(k)} 
    \left[ I_{n+2s+1}(k) + I_{n-2s-1}(k) \right] \left[ I_{n+2m+1}(k) + I_{n-2m-1}(k) \right],\\
    C_m &= \begin{cases}
        -\frac{1 + c}{2},& \text{if}\ m = 0, \\
        0,           & \text{otherwise},
    \end{cases}
    \end{split}
\end{equation}
where $I_n^\prime(k)$ is a derivative of the Bessel function, computed as
\begin{equation}
    I_n^\prime(k) = \frac{I_{n-1}(k) + I_{n+1}(k)}{2}.
\end{equation}
In case of large values of $c$, or even in the extreme $c=\infty$, the entries of the matrix $B_{ms}$ and the right hand side $C_m$ should be divided by $c$
to get alternative expressions with better numerical conditioning and thus more amenable to computation.

The analytical solution is expressed in non-dimensional variables as defined in \eqref{eq:reference-scales}--\eqref{eq:non-dim-number-last}.
The reference quantities are defined such that the non-zero components of the solution ($\rho'$, $T'$, $v'_z$, $B'_x$, $B'_z$)
have maximum magnitudes of order one. To this end, the reference length is defined as the pipe radius, that is $L^* = a$.
Due to incompressibility and temperature fluctuations damping mechanism, described in section \ref{pipeflow_testcase_section},
density and temperature are constant and arbitrary, therefore
the reference density and temperature can be set as $\rho^*=\rho$ and $T^*=T$ respectively.
The reference magnetic field magnitude is set to that of the external magnetic field, that is, $B^* = \Bext$.
The reference velocity is set to the velocity at the center line, that is $U^* = u(0,0)$.
The other reference quantities are defined though dimensionless numbers.
%The Mach number $\Ma$ and the Reynolds number $\Rh$ do not not significantly affect the solution 
%if they are small enough to not give rise of the turbulence and keep flow approximately incompressible.
%They are specified in Section \ref{pipeflow_testcase_section}.
The dimensionless external force $F'$ and the magnetic Mach and Reynolds numbers $\Mm$ and $\Rm$ are chosen 
such that the solution satisfies $u(0,0) = b(\frac{1}{2},\pi) = 1$
\footnote{At the coordinate $(r,\theta)=(\frac{1}{2},\pi)$ the axial magnetic field is close to its maximum value for all values of $c$.}.
%, and the Hartmann number has the prescribed value.
%The second condition comes from the requirement that $B'_z$ should be of magnitude one. 
%We enforce it to be equal one at the coordinate $(\frac{1}{2}, \pi)$ because at this coordinate $B_z$ has sufficient positive value for all $c$.
%With the prescribed values of wall conductance parameter $c$ and Hartmann number $\Ha$ we compute $\tilde u_0 = \tilde u(0,0)$ and $\tilde \Bext = \tilde b(\frac{1}{2},\pi)$ and then
To that end, we set
\begin{equation}
    F' = \frac{1}{\Rh \tilde u_\star}, \quad \Mm = \sqrt{\frac{\Rh \tilde u_\star}{\Ha \tilde b_\star}}, \quad \Rm = \Ha \frac{\tilde u_\star}{\tilde b_\star}. 
\end{equation}
where $\tilde u_\star = \tilde u(0,0)$ and $\tilde b_\star = \tilde b(\frac{1}{2},\pi)$.
As the solution corresponds to an incompressible flow regime, the Prandtl number and the heat capacity ratio are chosen arbitrarily and set to
$\gamma = \frac{7}{5}$ and $\Pr = \frac{3}{4}$.
With all the previous definition, the non-dimentional solution takes the following final form
\begin{equation}
    \rho' = 1, \quad \vec v' = (0, 0, \tilde u_\star^{-1} \tilde u), \quad T' = 1, \quad
    \vec B' = (1, 0, \tilde b_\star^{-1} \tilde b).
\end{equation}
%where $\Po = \Rh F'$ is the Poiseuille number.

The finite summation limits $N_0, S_0, N_1$, and $S_1$ in~\eqref{pipeflow_solution_eq}--\eqref{pipeflow_matrix_coeffitients_eq}
truncate the infinite series. These summation limits should be set to large enough values to achieve the required accuracy.
In practice, we set $N_0= N_1 = 30$ and $S_0 = S_1 = 10$.
At~$\Ha=5$, and for all the values of $c$ we tested, these summation limit values provide a maximum relative error
less than $10^{-15}$, which is close to machine epsilon for IEEE double precision floating point arithmetic.

\section{Raw data from convergence study}
Table \ref{converg_table} reproduces the raw data of the convergence study presented in Section \ref{pipeflow_testcase_section}.

\begin{sidewaystable}
    \caption{
        \label{converg_table} Convergence study for the pipe flow for solution polynomial degrees $p=1,2,3$; $Ha=5$; 
    wall conductance parameter $c=0,1,\infty$. The $L_2$ norm of the error in the axial velocity, $u$, and the axial magnetic field, $b$,
    are shown with the corresponding convergence rates. The number of mesh elements in the cross section at the refinement level $l$ is $5\cdot 4^l$.
    In axial direction number of mesh elements is constant and equal $3$.}
    %\small+
    \review{
    \begin{tabular}{ccccccccccccc}\hline
        \multirow{2}{*}{Grid} & \multicolumn{4}{l}{$c=0$} & \multicolumn{4}{l}{$c=1$} & \multicolumn{4}{l}{$c=\infty$} \\
        \cmidrule(lr){2-5} \cmidrule(lr){6-9} \cmidrule(lr){10-13}
        & $|\!|u-u_h|\!|_2$ & Rate & $|\!|b-b_h|\!|_2$ & Rate & $|\!|u-u_h|\!|_2$ & Rate & $|\!|b-b_h|\!|_2$ & Rate & $|\!|u-u_h|\!|_2$ & Rate & $|\!|b-b_h|\!|_2$ & Rate \\
        \hline
        $p=1$\\
        \hline
        0 & 5.71e$-$01 & $-$     & 6.57e$-$01 & $-$     & 3.66e$-$01 & $-$     & 9.28e$-$01 & $-$     & 4.41e$-$01 & $-$     & 8.49e$-$01 & $-$     \\
        1 & 1.44e$-$01 & $-$1.99 & 2.24e$-$01 & $-$1.55 & 1.62e$-$01 & $-$1.18 & 2.22e$-$01 & $-$2.06 & 2.17e$-$01 & $-$1.02 & 1.81e$-$01 & $-$2.23 \\
        2 & 3.53e$-$02 & $-$2.03 & 6.02e$-$02 & $-$1.89 & 4.60e$-$02 & $-$1.81 & 5.83e$-$02 & $-$1.93 & 5.91e$-$02 & $-$1.88 & 4.65e$-$02 & $-$1.96 \\
        3 & 8.07e$-$03 & $-$2.13 & 1.32e$-$02 & $-$2.19 & 1.24e$-$02 & $-$1.89 & 1.52e$-$02 & $-$1.94 & 1.58e$-$02 & $-$1.90 & 1.20e$-$02 & $-$1.95 \\
        4 & 1.83e$-$03 & $-$2.14 & 2.65e$-$03 & $-$2.32 & 3.20e$-$03 & $-$1.95 & 3.89e$-$03 & $-$1.96 & 4.11e$-$03 & $-$1.94 & 3.08e$-$03 & $-$1.96 \\
        \hline
        $p=2$\\
        \hline
        0 & 7.07e$-$02 & $-$     & 1.28e$-$01 & $-$     & 1.57e$-$01 & $-$     & 5.77e$-$02 & $-$     & 2.06e$-$01 & $-$     & 4.63e$-$02 & $-$     \\
        1 & 1.42e$-$02 & $-$2.32 & 2.70e$-$02 & $-$2.25 & 3.24e$-$02 & $-$2.27 & 1.33e$-$02 & $-$2.12 & 4.15e$-$02 & $-$2.31 & 9.75e$-$03 & $-$2.25 \\
        2 & 1.80e$-$03 & $-$2.98 & 3.68e$-$03 & $-$2.87 & 4.61e$-$03 & $-$2.81 & 2.40e$-$03 & $-$2.47 & 5.85e$-$03 & $-$2.83 & 1.76e$-$03 & $-$2.47 \\
        3 & 1.84e$-$04 & $-$3.30 & 3.91e$-$04 & $-$3.24 & 5.85e$-$04 & $-$2.98 & 3.43e$-$04 & $-$2.81 & 7.43e$-$04 & $-$2.98 & 2.55e$-$04 & $-$2.78 \\
        4 & 2.16e$-$05 & $-$3.09 & 3.75e$-$05 & $-$3.38 & 7.66e$-$05 & $-$2.93 & 4.38e$-$05 & $-$2.97 & 9.78e$-$05 & $-$2.93 & 3.26e$-$05 & $-$2.97 \\
        \hline
        $p=3$\\
        \hline
        0 & 8.84e$-$03 & $-$     & 2.25e$-$02 & $-$     & 2.36e$-$02 & $-$     & 1.00e$-$02 & $-$     & 2.90e$-$02 & $-$     & 8.11e$-$03 & $-$     \\
        1 & 1.03e$-$03 & $-$3.11 & 2.80e$-$03 & $-$3.01 & 3.02e$-$03 & $-$2.97 & 1.76e$-$03 & $-$2.51 & 3.77e$-$03 & $-$2.94 & 1.36e$-$03 & $-$2.58 \\
        2 & 6.49e$-$05 & $-$3.98 & 1.83e$-$04 & $-$3.94 & 2.38e$-$04 & $-$3.66 & 1.49e$-$04 & $-$3.56 & 2.99e$-$04 & $-$3.66 & 1.15e$-$04 & $-$3.57 \\
        3 & 3.31e$-$06 & $-$4.30 & 9.30e$-$06 & $-$4.30 & 1.59e$-$05 & $-$3.91 & 1.02e$-$05 & $-$3.87 & 2.00e$-$05 & $-$3.90 & 7.82e$-$06 & $-$3.87 \\
        4 & 1.57e$-$07 & $-$4.40 & 4.28e$-$07 & $-$4.44 & 1.01e$-$06 & $-$3.97 & 6.51e$-$07 & $-$3.97 & 1.29e$-$06 & $-$3.96 & 5.03e$-$07 & $-$3.96 \\
        \hline
    \end{tabular}
    } % review
\end{sidewaystable} 

\pagebreak

\section{Magnetic field of a wire with rectangular cross-section}
\label{rectangular_wire_appendix}
Here, we reproduce an analytical expression for the magnetic field induced by
\review{a} uniform electric current flowing through a conducting medium of infinite length
with rectangular cross section of width $w$ and thickness $h$~\cite{NTMDT}.
This analytical solution can be used to approximate the solution of a finite domain of length $L$ with $w \ll L$ and $h \ll L$.
In Cartesian coordinates \review{$(x, y, z)$}, with the current of density $i$ flowing in positive $y$ direction through the wire of rectangular cross section
$(x, z) \in [-w/2,w/2] \times [-h, 0]$, the components of the induced magnetic field $\vec B$ are
\begin{align*}
    B_x &= \tilde B_x\left(x + \frac{w}{2}, z\right) - \tilde B_x\left(x - \frac{w}{2}, z\right), \\
    B_y & = 0, \\
    B_z &= \tilde B_z\left(x + \frac{w}{2}, z\right) - \tilde B_z\left(x - \frac{w}{2}, z\right),
\end{align*}
where
\begin{align*}
    \tilde B_x(p, z) &= \frac{\mu_0 i}{4 \pi} \left[p \ln{\left(1 + \frac{h^2 + 2 z h}{p^2 + z^2}\right)} +
    2 (z + h) \arctan{\left(\frac{p}{z + h}\right)} - 2 z \arctan{\left(\frac{p}{z}\right)}\right], \\
    \tilde B_z(p, z) &= \frac{\mu_0 i}{4 \pi} \left[z \ln{\left(1 + \frac{p^2}{z^2}\right)} -
    (z + h) \ln{\left(1 + \frac{p^2}{(z + h)^2}\right)} -
    2 p \left(\arctan{\left(\frac{z + h}{p}\right)} - \arctan{\left(\frac{z}{p}\right)}\right)\right].
\end{align*}

\section{Magnetic field of a circular current loop}
\label{current_loop_appendix}
Here, we reproduce the expressions for the magnetic field induced by a closed circular current loop\cite{offaxisloop,jackson1999classical}. 
In cylindrical coordinates $(r,\theta,z)$, a total current $I$ flows through the circular loop with coordinates $z=0, r=a$ in the direction of incrementing $\theta$.
For any point in space that is not on the current loop, the components of the magnetic field are
\begin{align*}
    B_r &= B_0 \frac{\gamma}{\pi \sqrt{Q}}
          \left(E(k) \frac{1 + \alpha^2 + \beta^2}{Q - 4 \alpha} - K(k)\right), \\
  \quad B_\theta &= 0, \\
    B_z &= B_0 \frac{1}{\pi \sqrt{Q}}
    \left(E(k) \frac{1 - \alpha^2 - \beta^2}{Q - 4 \alpha} + K(k)\right),
\end{align*}
where
\begin{equation*}
    \alpha = \frac{r}{a}, \quad 
    \beta = \frac{z}{a}, \quad 
    \gamma = \frac{z}{r}, \quad
    Q = \left(1 + \alpha\right)^2 + \beta^2, \quad 
    k = \sqrt{\frac{4 \alpha}{Q}}, \quad
    B_0 = \frac{I \mu_0}{2 a},
\end{equation*}
and $K(k)$ and $E(k)$ are the complete elliptic integral functions of the first and second kind, respectively. 

%\printbibliography

\newpage

\bibliographystyle{elsarticle-num} 
\bibliography{mhd.bib}

%\end{document}

%\begin{thebibliography}{00}


%\bibitem{}

%\end{thebibliography}
\end{document}